\documentclass{amsart}
\usepackage{amsmath,amsfonts,latexsym,amssymb,amsthm}
\usepackage[mathscr]{eucal}
\newcommand{\R}{\ensuremath{\mathbb{R}}}
\newcommand{\Kr}{\ensuremath{K({\R})}}
\newcommand{\Lr}{\ensuremath{L({\R})}}
\newcommand{\mouse}{\ensuremath{\mathcal{M}}}
\newcommand{\nouse}{\ensuremath{\mathcal{N}}}
\newcommand{\overmouse}{\ensuremath{\overline{\mouse}}}
\newcommand{\core}{{\mathfrak{C}}}
\newcommand{\F}{\mathfrak{F}}
\newcommand{\f}{\mathfrak{f}}
\newcommand{\g}{\mathfrak{g}}
\newcommand{\h}{\mathfrak{h}}

\newcommand{\overcore}{\ensuremath{\overline{\core}}}
\newcommand{\fmouse}{\ensuremath{\mathfrak{M}}}
\newcommand{\fnouse}{\ensuremath{\mathfrak{N}}}
\newcommand{\A}{\ensuremath{\mathcal{A}}}
\newcommand{\B}{\ensuremath{\mathcal{B}}}
\newcommand{\calH}{{\mathcal{H}}}
\newcommand{\Q}{\ensuremath{\mathbb{Q}}}
\newcommand{\Lng}{\ensuremath{\mathcal{L}}}
\newcommand{\Jgamma}{\text{$J_\gamma^{\mouse}(\underline{\R})$}}
\newcommand{\BK}{\textup{BK}}
\newcommand{\Hull}{\textup{Hull}}
\newcommand{\pow}{\mathcal{P}}
\renewcommand{\P}{\ensuremath{\mathbb{P}}}
\renewcommand{\iff}{\textup{ \,\ iff \,\ }}
\newcommand{\AD}{\textup{AD}}
\newcommand{\DC}{\textup{DC}}
\newcommand{\ZF}{\textup{ZF}}
\newcommand{\AC}{\textup{AC}}
\newcommand{\ZFC}{\textup{ZFC}}
\newcommand{\Rtheory}{\textup{R}}
\newcommand{\Rplus}{\Rtheory^+}
\newcommand{\PM}{\textup{PM}}
\newcommand{\OR}{\textup{OR}}
\newcommand{\cOR}{{\mathfrak{O}}}
\newcommand{\en}{\text{$\in$}}
\newcommand{\Det}{\textup{Det}}
\newcommand{\dom}{\textup{dom}}
\newcommand{\ran}{\textup{ran}}
\newcommand{\ul}{\underline}
\newcommand{\ov}{\overrightarrow}
\newcommand{\wht}{\widehat}

\newcommand{\eq}{\!=\!}
\newcommand{\Eq}{{\underline{E}}}
\newcommand{\shrp}{{\sharp}}
\newcommand{\abs}[1]{\left|{#1}\right|}

\newcommand{\boldface}[2]{%
\protect\raisebox{0pt}[0pt][0pt]{%
$\underset{\displaystyle\widetilde{}}{\boldsymbol{#1}}{_{#2}}$}\mbox{\hskip 1pt}}
\newcommand{\boldfacemu}[2]{%
\protect\raisebox{0pt}[0pt][0pt]{%
$\underset{\displaystyle\widetilde{}}{\boldsymbol{#1}}{_{#2}^\mu}$}\mbox{\hskip 1pt}}
\newcommand{\mapsigma}[1]{\xrightarrow[\text{ \ \ $\Sigma_{#1}$}]{}}
\newcommand{\maps}[1]{\xrightarrow{\textup{#1}}}

\newcommand{\la}{\left\langle}
\newcommand{\ra}{\right\rangle}
\newcommand{\updownmap}[2]{\xrightarrow[\text{ \ \ $\Sigma_{#2}$}]{\text{#1}}} 
\newcommand{\premousesystem}{\text{$\la {\la \mouse_\alpha \ra}_{\alpha \in
\textup{OR}} ,  \la \pi_{\alpha\beta}\colon \mouse_\alpha \rightarrow
\mouse_\beta \ra_{\alpha\leq\beta\in\textup{OR}}\ra$}}
\newcommand{\premouseiteration}[3]
{\text{$\la {\la {#1}_{#2} \ra}_{{#2} \in \textup{OR}} , \la
\pi_{{#2}{#3}}\colon {#1}_{#2} \updownmap{cofinal}{1}{#1}_{#3}
\ra_{{#2}\leq{#3}\in\textup{OR}}\ra$}}
\newcommand{\premouseiterationover}[3]
{\text{$\la {\la {#1}_{#2} \ra}_{{#2} \in \textup{OR}} , \la
\overline{\pi}_{{#2}{#3}}\colon {#1}_{#2} \updownmap{cofinal}{1}{#1}_{#3}
\ra_{{#2}\leq{#3}\in\textup{OR}}\ra$}}
\newcommand{\nousensystem}{\text{$\la {\la \nouse_\alpha \ra}_{\alpha \in
\text{OR}} ,  \la \pi_{\alpha\beta}^n\colon \nouse_\alpha
\updownmap{cofinal}{1}\nouse_\beta \ra_{\alpha\leq\beta\in\text{OR}}\ra$}}
\newcommand{\mousesystem}{\text{$\la {\la \mouse_\alpha \ra}_{\alpha \in
\textup{OR}} ,  \la \pi_{\alpha\beta}\colon \mouse_\alpha \mapsigma{{n+1}}
\mouse_\beta \ra_{\alpha\leq\beta\in\textup{OR}}\ra$}}
\newcommand{\oversystem}{\text{$\la {\la \overline{\mouse}_\alpha \ra}_{\alpha \in
\textup{OR}} ,  \la \overline{\pi}_{\alpha\beta}\colon \overline{\mouse}_\alpha
\updownmap{cofinal}{1}\overline{\mouse}_\beta \ra_{\alpha\leq\beta\in\textup{OR}}\ra$}}
\newcommand{\coresystem}{\text{$\la {\la \core_\alpha \ra}_{\alpha \in
\textup{OR}} ,  \la \pi_{\alpha\beta}\colon \core_\alpha \mapsigma{{n+1}}
\core_\beta \ra_{\alpha\leq\beta\in\textup{OR}}\ra$}}
\newcommand{\overcoresystem}{\text{$\la {\la \overcore_\alpha \ra}_{\alpha \in
\textup{OR}} ,  \la \overline{\pi}_{\alpha\beta}\colon \overcore_\alpha
\updownmap{}{1}\overcore_\beta \ra_{\alpha\leq\beta\in\textup{OR}}\ra$}}
\newcommand{\premouseiterationt}[4]
{\text{$\la {\la {#1} \ra}_{{#2} \in \textup{OR}} , \la
\pi_{{#2}{#3}}\colon {#1} \updownmap{cofinal}{1}{#4}
\ra_{{#2}\leq{#3}\in\textup{OR}}\ra$}}

\newcommand{\be}{\begin{enumerate}}
\newcommand{\ee}{\end{enumerate}}
\newcommand{\bi}{\begin{itemize}}
\newcommand{\ei}{\end{itemize}}

\newtheorem{theorem}{Theorem}[section]
\newtheorem{corollary}[theorem]{Corollary}
\newtheorem{proposition}[theorem]{Proposition}
\newtheorem{lemma}[theorem]{Lemma}
\newtheorem*{main1}{Theorem \ref{GenDJ}}

\newtheorem*{question}{Question}

\theoremstyle{definition}
\newtheorem{definition}[theorem]{Definition}

\newtheorem{remark}[theorem]{Remark}
\newtheorem*{claiim}{Claim}
\newtheorem{clam}{Claim}[theorem]

\theoremstyle{remark}
\newtheorem*{rmk}{Comment}

\begin{document}

\title[Scales and the fine structure of $K(\R)$. Part I]{Scales and the fine structure of $\boldsymbol{K(\pmb{\R})}$\\ {Part I: Acceptability above the reals}}

\author{Daniel W. Cunningham}
\address{Mathematics Department,
State University of New York,
College at Buffalo,\\
1300 Elmwood Avenue,
Buffalo, NY 14222, USA}
\email{cunnindw@math.buffalostate.edu}

\begin{abstract}
This article is Part I in a series of three papers devoted to determining the minimal complexity of scales in the inner model \ $\Kr$. \ Here, in Part I, we shall complete our development of a fine structure theory for \ $\Kr$ \ which is essential for our work in Parts II and III. In particular, we prove the following fundamental theorem which supports our analysis of scales in \ $\Kr$: \ If \ $\mouse$ \ is an iterable real premouse, then \ $\mouse$ \ is acceptable above the reals. This theorem  will be used in Parts II and III to solve the problem of finding scales of minimal complexity in \ $\Kr$.
\end{abstract}

\maketitle

\section{Introduction}
\label{intro}
In \cite{Crcm} we introduced the Real Core Model \ $\Kr$ \  and showed that \ $\Kr$ \ is an inner model containing the reals and definable scales beyond those in \ $\Lr$. \ To establish our results on the existence of scales, we defined iterable real premice (see subsection \ref{realpremice} below) and showed how the {\it basic\/} fine structural notions of Dodd-Jensen \cite{DJ} generalize to iterable ``premice above the reals.'' Consequently, we were able to prove the following theorem (see \cite[Theorem 4.4]{Crcm}).

\begin{theorem}[$\ZF+\DC$]\label{firstscales}  Suppose that \ $\mouse$ \  is an iterable real premouse satisfying \ $\AD$. \ Then \ $\Sigma_1(\mouse)$ \ has the scale property. 
\end{theorem}

Theorem \ref{firstscales}, together with its proof, allowed us to determine the extent of scales in \ $\Kr$ \ and to prove that \ $\ZF+\AD+V\eq\Kr$ \ implies \ $\DC$ \ (see \cite{Crcm}). By allowing for real parameters in the proof of Theorem \ref{firstscales}, one can derive the following corollary. 

\begin{corollary}[$\ZF+\DC$]\label{relscales}  Suppose that \ $\mouse$ \  is an iterable real premouse satisfying \ $\AD$. \ Then \ $\Sigma_1(\mouse,\R)$ \ has the scale property. 
\end{corollary}

\begin{remark}\label{Rparameter} Throughout this paper (as in Parts II and III) we always allow the {\it set\/} of reals \
$\R$ \ to appear as a constant in our relevant languages (see subsections \ref{Preliminaries} and \ref{abovereals}). Thus for \ $n\ge 1$, \ the pointclass \ $\Sigma_n(\mouse)$ \ is in fact equal to the pointclass \ $\Sigma_n(\mouse,\{\R\})$; \ however, the pointclass \ $\Sigma_n(\mouse)$ \ is not necessarily equal to the pointclass \ $\Sigma_n(\mouse,\R)$.
\end{remark}

Given an iterable real premouse \ $\mouse$ \ and an arbitrary \ $m\ge 1$, \ Theorem \ref{firstscales} provokes:

\begin{question}[Q]When does the boldface pointclass $\boldface{\Sigma}{m}(\mouse)$
 have the scale property?
\end{question}

Before we address this question,  recall that \ $\mathcal{M}=(M,\R,\kappa,\mu)$ \ is a real 1--mouse
(see subsection \ref{realonemice}),  if \ $\mouse$ \ is an iterable real premouse and \ $\mathcal{P}(\R\times\kappa)\cap
\boldface{\Sigma}{1}(M)\not\subseteq M$, \ where \ $M$ \ has the form \ $J_\alpha[\mu](\R)$ \ and \ $\kappa$ \ is the
``measurable cardinal'' in \ $\mouse$. \ 

Real 1--mice suffice to define the real core model and to prove the results in
\cite{Crcm} about \ $\Kr$; \ however, real 1--mice are not sufficient to construct scales of minimal complexity. Our
solution to the problem of identifying these scales in \ $\Kr$ \ requires the development of a {\sl full\/}
fine structure theory for \ $\Kr$. \ In the paper \cite{Cfsrm} we initiated this development by generalizing
Dodd-Jensen's notion of a mouse to that of a {\it real mouse\/}  (see subsection \ref{realmice} below). This is
accomplished by 
\bi
\item isolating the concept of {\it acceptability above the reals\/}\footnote{This concept extends the Dodd-Jensen notion of
acceptability to include the set of reals.}  (see Definition
\ref{acceptable}),
\item replacing \ $\Sigma_1$ \ with \ $\Sigma_n$, \ where \ $n$ \ is the smallest integer such that \ 
$\mathcal{P}(\R\times\kappa)\cap\boldface{\Sigma}{n+1}(\mouse)\not\subseteq M$, 
\item defining an iteration procedure stronger than premouse iteration.
\ei

We now give a quick definition of a weak real mouse \ $\mouse$ \ and the natural number \ $m(\mouse)$. \  Let
\ $\mouse$ \ be a real mouse. Assume that there is an integer \ $m\ge 1$ \ such that
\ $\mathcal{P}(\R)\cap\boldface{\Sigma}{m}(\mouse)\not\subseteq M$ \ and let \ $m(\mouse)=m$ \ be the least such integer. We
say that \ $\mouse$ \ is {\it weak\/}  if 
\be
\item[(1)] $\mouse$ \ is a proper initial segment of an iterable real premouse and 
\item[(2)] $\mouse$ \ realizes a \ $\Sigma_m$ \ type not realized in any proper initial segment of \ $\mouse$. 
\ee
In (2) a \ $\Sigma_m$ \ type is a non-empty subset \ $\Sigma$ \ of \[\{ \theta\in \Sigma_m\cup\Pi_m  : \text{$\theta$
is a formula of one free variable}\}\] and \ $\mouse$ \ is said to realize \ $\Sigma$ \ if there is an \ $a\in M$ \
such that \ $\mouse\models\theta(a)$ \ for all \ $\theta\in\Sigma$.

In Part II \cite{Part2}, we shall present a partial answer to question (Q). Using the fine structure of real mice
developed  here and in \cite{Cfsrm}, we shall prove in Part II  the following theorem on the existence of scales:
\begin{theorem}[$\ZF+\DC$]\label{newthm}  Suppose that  \ $\mouse$ \  is a weak real mouse satisfying \ $\AD$. \ Then \
$\boldface{\Sigma}{m}(\mouse)$ \ has the scale property when \ $m=m(\mouse)$. 
\end{theorem}
We note that the above theorem requires {\it only\/} the determinacy of sets of reals in \ $\mouse$.\footnote{A critical
property for proving that \ $\Kr$ \ satisfies \ $\AD$ \ under certain large cardinal hypothesis.} The proof in
\cite{Part2} of Theorem
\ref{newthm} relies heavily on the fine structure of real mice; in particular, the proof relies on the fact that real
mice are  acceptable above the reals.  One might ask: {\sl Do weak real mice exist?\/} In the current paper we shall
prove the following fundamental theorem concerning acceptability above the reals (see Definition \ref{acceptable}). This
theorem will allow us to show that weak real mice do exist.
\begin{main1} Suppose that \ $\mouse$ \ is an iterable real premouse.  Then \ $\mouse$ \ is acceptable above the reals.
\end{main1}

Theorem \ref{GenDJ} and its proof are key components in our examination of the structure of
\ $\Kr$ \ and our analysis of its scales. For example, this theorem is used in Part III \cite{Part3} to
show that weak real mice, in fact, do exist in
\ $\Kr$.\footnote{We are implicitly assuming that $\R^\shrp$ exists.} \ In addition, Theorems \ref{newthm} and
\ref{GenDJ} are essential for our work in Part III, because they are used to justify our analysis of scales in \
$\Kr$ \ at the ``end of a gap''\footnote{A gap is an interval of ordinals in which no new $\Sigma_1$  truths about
the reals occur (see \cite{Part3}). New pointclasses will occur at the end of a gap.} and they allow us to obtain
scales of minimal complexity in \ $\Kr$. \ Consequently, in Part III we will be able to 
\bi
\item give a comprehensive answer to question (Q) and 
\item give a complete description of those levels of the Levy hierarchy for \ $\Kr$ \ possessing the scale property. 
\ei
In short, for our fine structural analysis of \ $\Kr$ \ the essential property is {\sl acceptability above the reals.\/}
In addition, Theorem \ref{GenDJ} is a critical component in the proofs of the major theorems
already established in \cite{Cnot}.\footnote{\label{justification} Theorem \ref{GenDJ} was first stated (without
formal proof) in \cite{Cnot}. It was also stated there that a proof of this theorem would be presented here.}  For
example, using Theorem \ref{GenDJ} we show (see \cite[4.5 \& 4.23]{Cnot}) that 
\bi
\item  $\ZF + \AD +
\exists X\subseteq\R\,[X\notin\Kr]$ \ implies the existence of an inner model of \ $\ZF + \AD + \DC$ \ containing a measurable
cardinal above its \ $\Theta$
\item $\ZF+\AD+\lnot\DC_{\R}$ \ implies that \ $\R^\dag$ \ (dagger) exists.
\ei
\begin{remark}\label{Justine}
In \cite{Part2} and \cite{Part3} we will be assuming the axiom of determinacy in order to produce scales in \
$\Kr$. \ Since \ $\AD$ \ implies that there is no well-ordering of the reals \ $\R$, \ we must not appeal to the
axiom of choice ($\AC$) in our study of \ $\Kr$. \ Theorem \ref{GenDJ} is a generalization of Lemma 5.21 of
Dodd-Jensen \cite{DJ}. The Dodd-Jensen proof of Lemma 5.21 exploits the axiom of choice in two different ways.
First, the Dodd-Jensen proof uses the fact that (Dodd-Jensen) premice satisfy the axiom of choice and secondly,
their proof presumes that \ $\AC$ \ holds in the universe. \ $\AD$, \ on the other hand, implies that ``premice
above the reals $\R$'' (see Definition \ref{premice}) fail to satisfy the axiom of choice and it also implies that \
$\AC$ \ fails in the universe. Thus it is critical that we do not inadvertently apply the axiom of choice in our proof
of Theorem \ref{GenDJ}. In fact, we shall present a proof of this theorem that relies on no principles of
choice.
\end{remark}

The current paper is organized into four sections. Section 1 offers an introduction and identifies our basic notation (see subsection \ref{Preliminaries} below). Section \ref{ZFcards} focuses on showing that certain {\sl relevant\/} cardinality calculations hold without the axiom of choice. These calculations will be used in our proof of Theorem \ref{GenDJ}. In Section \ref{finestructure} we present an overview together with some new results concerning the fundamental notions presented in \cite{Crcm} and \cite{Cfsrm} which will be used here and in Parts II \& III. Finally, Section \ref{realaccept} is devoted to the proof of Theorem \ref{GenDJ}.  

\subsection{Preliminaries and notation}\label{Preliminaries}

Let \ $\omega$ \ be the set of all natural numbers.  \ $\R = {^\omega \omega}$ \ is
the set of all functions from \ $\omega$ \ to \ $\omega$. \ We call \ $\R$ \ the set of
reals and regard \ $\R$ \ as a topological space by giving it the product topology,
using the discrete topology on \ $\omega$. \ For a set \  $A \subseteq \R$ \ we
associate a two person infinite game on \ $\omega$, \ with {\it payoff} \
$A$, \ denoted by \ $G_A$:
\[
\begin{aligned}[c]
{}&{\mathbf{I}}\phantom{{\mathbf{I}}} \qquad x(0) \qquad \phantom{x(1)} \qquad x(2) \qquad \phantom{x(3)}\quad \\
{}&{\mathbf{I}}\mathbf{I} \qquad \phantom{x(0)}\qquad x(1) \qquad \phantom{x(2)} \qquad x(3) \quad 
\end{aligned}
\begin{gathered}[c]
{\cdots}
\end{gathered}
\] 
in which player \ {\bf I} \ wins if \ $x \in A$, \ and \ {\bf II} \ wins if \ $x
\notin A$. \ We say that \ $A$ \ is {\it determined\/} if the corresponding game \ $G_A$ \
is determined, that is, either player {\bf I} or  {\bf II} has a winning
strategy (see \cite[p. 287]{Mosch}). 
The {\it axiom of determinacy\/} \ ($\AD$) \ is a regularity hypothesis about games on \
$\omega$ \ and states: \ $\forall A \subseteq \R \  ( A  \text{\ is
determined})$.

We work in \  $\ZF$ \ and state our additional
hypotheses as we need them.  We do this to keep a close watch on the use of determinacy in the
proofs of our main theorems. Variables \ $x, y, z, w \dots$ \ generally range over \ $\R$, \ while
\ $\alpha, \beta, \gamma, \delta \dots$ \ range over \ $\OR$, \ the class of ordinals.  For \ $x
\in \R$ \ and \ $i\in\omega$ \ we write \ $\lambda.nx(n+i)$ \ for the real \ $y$ \ such that \
$y(n)=x(n+i)$ \ for all \ $n$. \ We write \ $(x)_i$, \ or \ $x_i$ \ when the context is clear, \ for the real  \ $z$
\ such that \ $z(n)=x(\langle n, i \rangle)$, \ where \ $\langle \ , \  \rangle$ \ recursively codes a pair of integers
by a single integer.  If \ $0\leq j\leq\omega$ \ and \ $1\leq k\leq\omega$, \ then \ $\omega^j \times
(^\omega\omega)^k$ \ is recursively homeomorphic to \ $\R$, \ and we sometimes tacitly identify the two.  The cardinal \ $\Theta$ \ is the supremum of the ordinals which are the surjective image of \ $\R$. 

A {\it pointclass\/} is a set of subsets of \ $\R$ \ closed under recursive substitutions. A
boldface pointclass is a pointclass closed under continuous substitutions.  For a
pointclass \ $\Gamma$, \ we write \ ``$\Gamma$--$\AD$'' or ``$\Det(\Gamma)$'' \ to denote
the assertion that all games on \ $\omega$ \ with payoff in \ $\Gamma$ \ are
determined. For the notions of a scale and of the scale property as well as any other
notions from Descriptive Set Theory which we have not defined, we refer the reader
to Moschovakis \cite{Mosch}.

A proper class \ $M$ \ is called an {\it inner model\/} if and only if \ $M$ \ is a transitive  $\in$--model of \ $\ZF$ \ containing all the ordinals. We distinguish between the notations \ $L[A]$ and $L(A)$. \ 
The inner model \ $L(A)$ \ is defined to be the class of sets constructible {\it above\/} \ $A$, \ that is, \ one starts
with a set \ $A$ \ and iterates definability in the language of set theory. \ Thus, \ $L(A)$ \ is  the smallest
inner model \ $M$ \ such that \ $A\in M$. \ The inner model \ $L[A]$ \ is defined to be the class of sets constructible {\it
relative\/} to \ $A$, \ that is, \ one starts with the empty set and iterates definability in the language of set theory
augmented by the predicate \ $A$. \ Consequently, \ $L[A]$ \ is the smallest inner model \ $M$ \ such that \ $A\cap M\in
M$ \ (see page 34 of \cite{Kana}). Furthermore, one defines \ $L[A,B]$ \ to be the class of sets constructible {\it relative\/}
to \ $A$ \ and \ $B$, \ whereas \ $L[A](B)$ \ is defined as the class of sets constructible {\it relative\/} to \ $A$ \ and 
{\it above\/}
\ $B$. \ Thus, \ $A\cap L[A](B)\in L[A](B)$ \ and \ $B\in L[A](B)$.

Our general set theoretic notation is standard. Given a function \ $f$, \ we write \ $\dom(f) = \{x : \exists y (f(x) =y)\}$ \ and \ $\ran(f) = \{y : \exists x (f(x) = y)\}$. \ We shall write \ $\langle x_1,\dots,x_n\rangle$ \ to represent a finite sequence of elements.
For any set \ $X$, \ $(X)^{<\omega}$ \ is the set of all finite sequences of elements
of \ $X$,  \ $[X]^{<\omega}$ \ is the set of all finite subsets of \ $X$, \ and \
$\mathcal{P}(X)$ \ is the set of all subsets of \ $X.$ \ Given two finite sequences \
$s$ \ and \ $t$, \ the sequence \ $s{^\frown}t$ \ is the concatenation of \ $s$ \ to
\ $t$. \ Generally, \ $\mu$ \ will be a normal measure on \ $\mathcal{P}(\kappa)$, \ where
\ $\kappa$ \ is an ordinal. For any ordinals \ $\eta \leq \alpha$, \ 
${^\eta\alpha}\!\uparrow$ \ is the set of all strictly increasing \ $\eta$ \ sequences
from \ $\alpha$. \ $V_\alpha$ \ is the set of all sets of rank less than \ $\alpha$. \ 
We let \ $y = T_c(x)$ \ denote the formula ``$y$ \ is the transitive closure of \ $x$.'' \ 
 For a model \ $\mouse = (M,\in, \dots),$ \ we shall abuse standard notation slightly and write
 \ ${^\kappa M}=\{f\in M \ \vert \ f\colon \kappa \rightarrow M \}$. \ In addition, for a model (or inner model) \
$\mouse$ \ having only one ``measurable cardinal,'' we shall write \ $\kappa^{\mouse}$ \ to denote this cardinal
in \ $\mouse$. \ Similarly, when \ $\mouse$ \ has only one ``measure,'' we shall write \ $\mu^{\mouse}$ \ to denote
this measure. 

Given a model \ $\mouse = (M,c_1,c_2,\dots,c_m,A_1,A_2,\dots,A_N),$ \ where the \  $A_i$ \
are predicates and the \ $c_i$ \ are constants,  if \ $X\subseteq M$ \ then \
$\Sigma_n(\mouse,X)$ \ is the class of relations on \ $M$ \ definable over \ $\mouse$ \ by a
\ $\Sigma_n$ \ formula from parameters in \ $X \cup \lbrace c_1,c_2,\dots,c_m\rbrace$. \
$\Sigma_\omega(\mouse,X) =  \bigcup\limits_{n\in\omega} \Sigma_n(\mouse,X).$  We write
``$\Sigma_n(\mouse)$'' for \ $\Sigma_n(\mouse,\emptyset)$ \ and \
``$\boldface{\Sigma}{n}(\mouse)$'' for  the boldface class \ $\Sigma_n(\mouse,M).$ \ 
Similar conventions hold for \ $\Pi_n$ \ and \ $\Delta_n$ \ notations.  If \ $\mouse$ \ is
a substructure of \ $\nouse$ \ and \ $X\subseteq M \subseteq N$, \ then 
``$\mouse \prec_n^X \nouse$'' means that \ $\mouse \models \phi[a] \textup{ \ iff
 \ } \nouse \models \phi[a]$, \ for all \ $a \in (X)^{<\omega}$ \ and for all \
$\Sigma_n$ \ formulae \ $\phi$ \ (the formula \ $\phi$ \ is allowed constants taken
from \ $\lbrace c_1,c_2,\dots,c_m\rbrace$). \ We write \ ``$\mouse \prec_n \nouse$'' \ 
for \ ``$\mouse \prec_n^M \nouse$.'' \ In addition, for any two models \ $\mouse$ \ and \
$\nouse$, \  we write \ $\pi:\mouse\mapsigma{n}\nouse$ \ to indicate that the map \ $\pi$ \ is
a \ $\Sigma_n$--elementary embedding, that is, \ $\mouse \models \phi[a] \textup{ \ iff
 \ } \nouse \models \phi[\pi(a)]$, \ for all \ $a=\langle a_0,a_1,\dots\rangle \in (M)^{<\omega}$ \
and for all \ $\Sigma_n$ \ formulae \ $\phi$, \ where \ $0\le n\le\omega$ \ and \ 
$\pi(a)=\langle \pi(a_0),\pi(a_1),\dots\rangle$. \ We shall write \ $\mouse \cong \nouse$ \ to denote that the structures \
$\mouse$ \ and \ $\nouse$ \ are isomorphic.

\section{Inner models of $\ZF$ above the reals}\label{ZFcards}
In \cite{Part2} and \cite{Part3} we will be assuming the axiom of determinacy in order to produce scales in \
$\Kr$. \ Thus, it is imperative that we do not use the axiom of choice in our study of \ $\Kr$. \ Our fine structural
analysis of \ $\Kr$ \ is a generalization of the Dodd-Jensen analysis of the core model \ $K$. \ In their analysis,
Dodd and Jensen freely apply the axiom of choice. For example, they use the ``facts'' that (1) the successor of a
cardinal is a regular cardinal and (2) given \ $< \kappa$ \ many sets each of size \ $<  \kappa$, \ their union has 
size \ $< \kappa$, \ whenever \ $\kappa$ \ is a regular cardinal. \ We must assure that our analysis of \ $\Kr$ \ does
not inadvertently appeal to the axiom of choice.  Therefore, in this section we show that certain relevant
cardinality calculations hold without the axiom of choice (see Theorem \ref{card}). The results of this section will
be used in our proof of Theorem \ref{GenDJ} which asserts that iterable real premice are acceptable.

\begin{definition} Let \ $M$ \ be a transitive set or class.
\be
\item $M$ \ is a {\it inner class model\/} if and only if
\ $M$ is an $\in$--model of \ $\ZF$ \ with \ $\OR\subseteq M$.
\item $M$ \ is an {\it inner set model\/} if and only if
\ $M$ is an $\in$--model of \ $\ZF$ \ with \ $\OR^M\in\OR$.
\item $M$ \ is an {\it inner model\/} (of \ $\ZF$) \ if and only if \ $M$ \ is an inner class or set
model.  
\ee
\end{definition}

If an inner model \ $M$ \ satisfies the axiom of choice, then every set in \ $M$ \ can be well-ordered and thus, every
such set has cardinality in terms of the ordinals in \ $M$. \ In addition, \ $\AC$ \ implies that successor
cardinals are regular and it also implies results concerning the cardinality of a union of sets. 

If we assume \ $\ZF+\AD$, \ then any inner model \ $M$ \ which contains the set of reals does
not satisfy the axiom of choice. In fact, \ $\mathcal{P}(\lambda)$ \ cannot be well-ordered
in \ $M$ \ for any ordinal \ $\lambda\ge\omega$, \ and thus the cardinality of \ $\mathcal{P}(\lambda)$ \ in \ $M$ \
makes no apparent sense. In this paper however, it will be necessary to modify the standard definition of cardinality
in \ $M$. \ One benefit of this modified notion of cardinality is that the ``cardinality'' of every set will exist in
all of the inner models \ $M$ \ that we consider. We will also be able to establish some important cardinal
inequalities in \ $M$. \ We will take advantage of these cardinal inequalities in our proof of Theorem \ref{GenDJ}.
\begin{definition}\label{R-card} Let \  $M$ \ be an inner model containing \ $\R$.
\be
\item The \ $M$--{\it cardinality} \ of a set \ $a\in M$, \
denoted by \ $\abs{a}_M$, \ is the least ordinal \
$\lambda\in{\OR}$ \ such that \ $f\colon \lambda \times
{\R}\maps{onto} a$ \ for some \ $f\in M$, \ if such an \ $f$ \ exists. 
\item An ordinal \ $\lambda$ \ is called an \ $M$--{\it cardinal\/} \ if \
$\lambda=\abs{\lambda}_M$.
\item For an ordinal \ $\lambda$, \ the least \ $M$--{\it cardinal\/} \ greater than \ $\lambda$ \
is denoted by \ $\lambda^+_M$. 
\item An \ $M$--{\it cardinal\/} \ $\lambda$ \ is called $M$--{\it regular\/} if there is no function \ $f\in M$ \ such
that \ $f\colon \xi \times {\R}\maps{cofinal} \lambda$ \ where \ $\xi<\lambda$.\footnote{A map \ $f\colon \xi \times {\R}\to \lambda$ \ is cofinal if for all \ $\beta\in\lambda$ \ there is an \ $\alpha\in\xi$
\ and an \ $x\in\R$ \ such that \ $\beta\in f(\alpha,x)$.}
\ee
\end{definition}

\begin{rmk} We have chosen to use the cross product \ $\lambda\times{\R}$ \ in our definition of
\ $M$--{cardinality}. \  Since the cross product \ $0\times{\R}=\emptyset$ \ does not
involve the reals, we shall abuse standard notation and redefine \ $0\times{\R}={\R}$. \ Consequently, if \ $f\colon  {\R}\maps{onto} a$ \ for some \ $f\in M$, \ then \ $\abs{a}_M=0$. 
\end{rmk}

We note that if there is a surjection \ $F\colon \OR\times\R\to M$ \ and \ $F$ \ is definable
over \ $M$ \  such that the restriction \ $F\restriction(\lambda\times\R)$ \ is in \ $M$ \ for all ordinals \
$\lambda$, \ then the  \ $M$--cardinality of every set in \ $M$ \ exists. In particular, the \
$\Lr$--cardinality of every set in \ $\Lr$ \ exists.  

We must distinguish between the notion of an ordinal \ $\lambda$ \ being an \ $M$--{\it
cardinal} \ and the notion of \ $\lambda$ \ being a {\it standard\/} cardinal in \ $M$.

\begin{definition} Let \  $M$ \ be an inner model containing \ $\R$.
\be
\item The \ {\it cardinality\/} in \ $M$ \ of a set \ $a\in M$, \
denoted by \ $\abs{a}^M$, \ is the least ordinal \
$\lambda\in{\OR}$ \ such that \ $f\colon \lambda \maps{onto} a$ \ for some \ $f\in M$, \ if such an
\ $f$ \ exists. 
\item An ordinal \ $\lambda$ \ is called a {\it cardinal\/} in \ $M$ \ if \
$\lambda=\abs{\lambda}^M$.
\item For an ordinal \ $\lambda$, \ the least {\it cardinal\/} in \ $M$ \ greater than \ $\lambda$
\ is denoted by \ $\left(\lambda^+\right)^M$. 
\ee 
\end{definition}

\begin{rmk}
With respect to the above definition, we shall denote \ $\abs{\lambda}^M$ \ and \
$\left(\lambda^+\right)^M$ \ by \ $\abs{\lambda}$ \ and \ $\lambda^+$, \ respectively, unless stated otherwise.
\end{rmk}

Let \ $M$ \ be an inner model containing \ $\R$. \ We will show that ``almost all'' \ $M$--cardinals are standard
cardinals in \ $M$. \ First we prove the following lemma. Recall that the ordinal \
$\Theta$ \ is the supremum of the ordinals which are the surjective image of \ $\R$. 

\begin{lemma}\label{lem2}  Let \  $M$ \ be an inner model containing \ $\R$.  \ For
all ordinals \ $\xi\ge\lambda\ge\Theta^{M}$, \ if there is an \ $h\in M$ \ such that \
$h\colon\lambda\times\R\maps{onto}\xi$, \ then there is a function \ $g\in M$ \ such that \ 
$g\colon\lambda\maps{onto}\xi$.
\end{lemma}

\begin{proof}
We work in \ $M$. \ Let \ $\xi\ge\lambda\ge\Theta$. \ Suppose that \ $h$ \ is such
that \ $h\colon\lambda\times\R\maps{onto}\xi$. \ Since \ $\lambda\ge\Theta$, \
$\abs{\lambda}\ge\Theta$. \  For each \ $\alpha\in\lambda$, \ let \
$h_\alpha\colon\R\to \xi$ \ be defined by \ $h_\alpha(x)=h(\alpha,x)$. \ So the range of \
$h_\alpha$ \ has order type strictly less than \ $\Theta$. \ Let \
$\sigma_\alpha\colon\ran(h_\alpha)\maps{onto}\delta_\alpha$
\ be the ``collapse map'' of \ $\ran(h_\alpha)$ \ onto \ $\delta_\alpha<\Theta$. \ 
Define \ $q$ \ where \ $q\colon\lambda\times\Theta\to \xi$ \ by 
\[q(\alpha,\gamma)=
\begin{cases} \beta, &\text{if $\gamma<\delta_\alpha$, $\sigma_\alpha(h_\alpha(x))=\gamma$ and $\beta=h(\alpha,x)$ for}\\
&\text{ \ some $x\in\R$;}\\
 0,&\text{otherwise.}\end{cases}\]
The function \ $q$ \ is well-defined and onto. Since \ $\Theta\le\lambda$, \ it follows that \ 
$\abs{\lambda\times\Theta}=\abs{\lambda}$. \ Thus,  there exists a function \ 
$g\colon\lambda\maps{onto}\xi$.\end{proof}

\begin{corollary} \label{corM} Let \ $M$ \ be an inner model containing \ $\R$. \ Then the
following is true in \ $M$: \ for an ordinal \ $\xi$ \ we have
\setcounter{equation}{0}
\begin{equation}\label{rcard}
\abs{\xi}_{M} =
\begin{cases}
0,  &\textup{if \ $\xi<\Theta$;}\\
\abs{\xi},  &\textup{if \ $\xi\ge\Theta$,}
\end{cases}
\end{equation}
and
\begin{equation}\label{rcardplus}
\xi^+_{M} =
\begin{cases}
\Theta,  &\textup{if \ $\xi<\Theta$;}\\
\xi^+,  &\textup{if \ $\xi\ge\Theta$.}
\end{cases}
\end{equation}
In particular, \ the equations \ $\abs{\Theta}_{M}=\abs{\Theta}^{M}=\Theta$ \ and \
$\Theta^+_{M}=\Theta^+$ \ are true in \ $M$.
\end{corollary}

\begin{proof} We work in \ $M$. \ Since equation (\ref{rcardplus}) follows from equation
(\ref{rcard}), we shall just prove equation (\ref{rcard}). Let \ $\xi$ \ be an ordinal. If \
$\xi<\Theta$, \ then \
$\abs{\xi}_M=0$. \ If \ $\xi\ge\Theta$, \ then clearly \ $\abs{\xi}_M\le\abs{\xi}\le\xi$. \
Suppose for a contradiction, that \ $\abs{\xi}_M=\lambda<\abs{\xi}$. \ Hence, there is an
\ $h$ \ such that \ $h\colon\lambda\times\R\maps{onto}\xi$. \ Lemma~\ref{lem2} implies there is a
function \ $g$ \ such that \ $g\colon\lambda\maps{onto}\xi$. \ Therefore, \ $\abs{\xi}$ \ is not a
cardinal. This contradiction completes the proof. 
\end{proof}

\begin{remark}\label{recall} Let \ $M$ \ be any inner model containing the reals. Corollary~\ref{corM} implies
that
\begin{itemize}
\item $\abs{0}_M=0$ \ and \ $0^+_M=\Theta^M$ \ are the first two \ $M$--cardinals, and
\item  for any ordinal \ $\xi\ge\Theta^{M}$, \ $\xi$ \ is a (standard) cardinal in \ $M$ \ if and only if
\ $\xi$ \ is an \ $M$--cardinal.
\end{itemize}
\end{remark}

Let \ $M$ \ be an inner model of \ $\ZF$ \ and let \ $\P$ \ be a partial order in \ $M$. \ Given that \ $G$ \ is \
$\P$--generic over \ $M$, \ we let \ $M[G]$ \ be the resulting generic extension of \ $M$. \ If \ $M[G]\models\ZFC$ \ then
we shall call \ $M[G]$ \ a {\it $\ZFC$--generic extension.\/}  For any forcing concepts and notation which we do not
define, we refer the reader to Kunen \cite{Kunen}. 
We now define a partial order \ $\Q$ \
which will force the set of reals in the ground model to be a countable set in any generic extension.
\begin{definition}\label{col_to_omega}
We shall let \ $\Q=(Q,\le)$ \ denote the standard partial order that produces  a generic enumeration of the reals in
length \ $\omega$; \ that is, let  \ $Q=\{s\in {^n}\R : n\in\omega\}$ \  and for \ $s,t\in Q$, \ define \ $s\le t$ \ if
and only if \ $\dom(s)\ge\dom(t)$ and  $t=s\restriction\dom(t)$.
\end{definition} 

\begin{rmk} Let \ $M$ \ be an inner model containing \ $\R$. Then \ $\Q\in M$.
\end{rmk}

Our next result is the key lemma which allows us to establish certain cardinality calculations without the
axiom of choice. 
\begin{lemma} \label{lem1} Let \ $M$ \ be an inner model containing \ $\R$. \
Suppose that \ $G$ \ is \ $\Q$--generic over \ $M$.  \ Then for all ordinals \
$\lambda\ge\omega$ \ and all sets \ $X\in M$,
\[(\exists f\in M[G]) \left[f\colon\lambda\maps{onto}X\right] \Longleftrightarrow (\exists
g\in M)
\left[g\colon\lambda\times\R\maps{onto}X\right].\]
\end{lemma}

\begin{proof} Let \ $X\in M$ \ and let \ $\lambda$ \ be an ordinal where \ $\lambda\ge\omega$. \

\medskip
($\Rightarrow$). \ Assume that \ $f\in M[G]$ \ is such that \ $f\colon\lambda\maps{onto}X$.
\  Let \ $\dot{f}$ \  be a \ $\Q$--name  for \ $f\in M[G]$, \ and let \ $\check a$ \ be the canonical
\ $\Q$--name  for any \ $a\in M$. \ Let \ $p\in G$ \ be such that \[M\models\text{``}p\Vdash
\dot{f}\colon\check \lambda\maps{onto}\check X\text{''}.\] 
Define the map \ $h\colon\lambda\times\Q\maps{onto}X$ \ in \ $M$ \ by 
\[h(\alpha,q)=\begin{cases} a, &\text{if \ $q\in\Q$, $a\in X$ and
$p{^\frown}q\Vdash\dot{f}(\check\alpha)=\check a$};\\
 0,&\text{otherwise.}\end{cases}\]
Here \ $p{^\frown}q$ \ is the concatenation of \ $q$ \ to \ $p$. \ Since \ $h\in M$ \ and 
there is a map in \ $M$ \ from \ $\R$ \ onto \ $\Q$, \ it follows that there is a map \ $g\in M$ \
such that \ $g\colon\lambda\times\R\maps{onto}X$.

\medskip
($\Leftarrow$). \ Assume that \ $g\in M$ \ is such that \
$g\colon\lambda\times\R\maps{onto}X$. \ Because  there is a bijection in \ $M[G]$ \ between \
$\R^{M}$ \ and \ $\omega$, \ it follows that that there is a map \ $h\in M[G]$ \ such that \
$h\colon\lambda\times\omega\maps{onto}X$. \ Since
\ $M[G]\models\abs{\lambda\times\omega}=\abs{\lambda}$ \ (because \ $\lambda\ge\omega$), \ there is
a map \ $f\in M[G]$ \ such that \ $f\colon\lambda\maps{onto}X$.
\end{proof}

Remark \ref{recall} implies that for ordinals above or equal to \ $\Theta$ \ there is no difference between
\ $M$-cardinals and the ``standard'' cardinals in \ $M$. \ Therefore, \ $\Q$--forcing preserves all the standard
cardinals \
$\ge\,\Theta$.

\begin{corollary} \label{cardpreserve} Let \ $M$ \ be an inner model containing \
$\R$. \ Suppose that \ $G$ \ is \ $\Q$--generic over \ $M$.  \ Then for all ordinals \
$\kappa\ge\Theta^{M}$, \ the following are equivalent:
\be
\item[(1)] $M\models\text{``$\kappa$ is a standard cardinal''}$,
\item[(2)] $M\models\text{``$\kappa$ is an $M$--cardinal''}$,
\item[(3)] $M[G]\models\text{``$\kappa$ is a standard cardinal''}$.
\ee
Thus, for any ordinal \ $\kappa\ge\Theta^{M}$ \ in \ $M$, \ if \ $\kappa$ \ satisfies one of the above
three conditions, then \ $(\kappa^+)^M=\abs{\kappa}_M^+=(\kappa^+)^{M[G]}$.
\end{corollary}

\begin{proof} Let \ $\kappa\in M$ \ be such that \ $\kappa\ge\Theta^{M}$. \ Corollary~\ref{corM} shows that (1) and (2) are
equivalent.  \ Lemma \ref{lem1} implies that (2) and (3) are equivalent.
\end{proof}

We note that \ $\Q$--forcing collapses \ $\Theta$ \ to \ $\omega_1$.
\begin{corollary} \label{cor4} Let \ $M$ \ be an inner model containing \ $\R$. \
Suppose that \ $G$ \ is \ $\Q$--generic over \ $M$.  \ Then \ $\Theta^{M}=\omega_1^{M[G]}$.
\end{corollary}

\begin{proof} To prove that \ $\Theta^{M}=\omega_1^{M[G]}$, \ we note that for any ordinal \
$\xi>\omega$, \ Lemma~\ref{lem1} implies that \ $\xi$ \ is a cardinal in \ $M[G]$ \ if and only if \
$\xi$ \ is an \ $M$--cardinal. Hence, Corollary \ref{corM} implies that \ $\Theta^M$ \ is 
collapsed to \ $\omega_1$. \end{proof}

\begin{theorem}\label{cardM} Let \ $M$ \ be an inner model containing \
$\R$. \ Let \ $G$ \ be \ $\Q$--generic over \ $M$ \ and suppose that \ $M[G]$ \ is a \ $\ZFC$--generic extension.
Then the following hold:
\be
\item [(i)] For all sets \ $X\in M$, \  $M\models\text{``$|X|_M$  exists''}$.
\item[(ii)] For all \ $M$--cardinals \ $\kappa\ge\Theta^M$, \ $\kappa^+_M$ \ is an \ $M$--regular cardinal. 
\item[(iii)] Suppose that \ $\kappa$ \ is an \ $M$--regular cardinal. Let \ $\langle X_{\langle\alpha, x\rangle} :
\langle\alpha, x\rangle\in \lambda\times\R\rangle$ \ be any sequence of sets in \ $M$, \ where \ $\lambda<\kappa$. \ If \
$\abs{X_{\langle\alpha, x\rangle}}_M<\kappa$ \ for all \ $\langle\alpha, x\rangle\in \lambda\times\R$, \ then \
$\abs{\bigcup\limits_{\langle\alpha, x\rangle\in \lambda\times\R}X_{\langle\alpha, x\rangle}}_M<\kappa$.
\ee
\end{theorem}

\begin{proof} Assume \ $M[G]\models\ZFC$.

\medskip
\noindent (i) Let \ $X\in M$ \ and so, \ $X\in M[G]$. \ Since \ $M[G]\models\AC$, \ there is an ordinal \ $\lambda\in M$
\ and \ $f\in M[G]$ \ such that \ $f\colon\lambda\maps{onto}X$. \ Lemma \ref{lem1} now implies that \
$M\models\text{``$|X|_M$  exists''}$.

\medskip
\noindent (ii)  Let \ $\kappa\ge\Theta^M$ \ be an \ $M$--cardinal.  We shall prove that 
\[M\models \text{``$\kappa^+_M$ \ is an \ $M$--regular cardinal''.}\] 
Suppose otherwise, for a contradiction. It then follows that there is a function \ $g\in M$ \ such
that \ $g\colon \kappa \times {\R}\maps{cofinal} \kappa^+_M$. \ 
By Corollary \ref{cardpreserve} \ both \ $\kappa$ \ and \  $\kappa^+_V$ \ are cardinals in \ $M[G]$ \ and, in addition,
\ $\kappa^+_M=(\kappa^+)^{M[G]}$. \  By Lemma \ref{lem1}, the existence of the function \ $g$ \ implies that there is a
function \ $f\in M[G]$ \ such that \ $f\colon \kappa\maps{cofinal}(\kappa^+)^{M[G]}$. \ Since \ $M[G]\models\ZFC$, \ this is
impossible.

\medskip
\noindent (iii) Suppose that \ $\kappa$ \ is an \ $M$--regular cardinal. Let \ $\langle X_{\langle\alpha, x\rangle} :
\langle\alpha, x\rangle\in \lambda\times\R\rangle$ \ be a sequence of sets in \ $M$ \ such that  \
$\abs{X_{\langle\alpha, x\rangle}}_M<\kappa$ \ for all \ $\langle\alpha, x\rangle\in
\lambda\times\R$. \ By (i) we have that \ $\abs{\bigcup\limits_{\langle\alpha, x\rangle\in
\lambda\times\R}X_{\langle\alpha, x\rangle}}_M$
\ exists. We shall prove that \ $M\models\abs{\bigcup\limits_{\langle\alpha, x\rangle\in
\lambda\times\R}X_{\langle\alpha, x\rangle}}_M<\kappa$.  \ By Corollary \ref{cardpreserve}, \ $\kappa$ \  is a
cardinal in \ $M[G]$ \ and, in addition, Lemma
\ref{lem1} implies that \ $\kappa$ \ is regular in \ $M[G]$. \ Since \
$\abs{X_{\langle\alpha, x\rangle}}_M<\kappa$ \ for each \ $\langle\alpha, x\rangle\in \lambda\times\R$ \ holds in \ $M$,
\ Lemma
\ref{lem1} implies that \ $\abs{X_{\langle\alpha, x\rangle}}^{M[G]}<\kappa$ \ for each \ $\alpha<\lambda$ \ and \
$x\in\R^M$. \ Because, \ $\R^M$ \ is countable in \ $M[G]$, \ we have that \ $M[G]\models
\abs{\R^M\times\lambda}<\kappa$.
\ Therefore, since \ $M[G]\models\ZFC$, \ it follows that \ 
$M[G]\models \abs{\bigcup\limits_{\langle\alpha, x\rangle\in \lambda\times\R^M}X_{\langle\alpha,
x\rangle}}<\kappa$. \ Lemma \ref{lem1} now, again, implies that
\ $M\models\abs{\bigcup\limits_{\langle\alpha, x\rangle\in \lambda\times\R}X_{\langle\alpha, x\rangle}}_M<\kappa$.
\end{proof}

We will now concentrate on inner models of the form \ $L(\R)[A]$ \ where \ $A\subseteq \R\times(\OR)^{<\omega}$. 

\begin{lemma}\label{choice} Let \ $A\subseteq \R\times(\OR)^{<\omega}$. \ Assume \ $\ZF+ V\eq L(\R)[A]$. \ Suppose that \ $G$
\ is \ $\Q$--generic over \ $V$. \ Then \ $V[G]$ \ is a \ $\ZFC$--generic extension.
\end{lemma}

\begin{proof} Note that \ $V[G]=L(\R)[A][G]=L(G)[A]$ \ where \ $A\subseteq \R^V\times(\OR)^{<\omega}$ and \ $G$ \ is
(essentially) a single real which codes all of the reals in \ $V$. \ Hence, \ $V[G]$ \ is a model of the axiom of
choice. Thus, \
$V[G]\models\ZFC$.\end{proof}

In subsection \ref{realmice} we state and prove Theorem \ref{GenDJ}, which is a generalization of Lemma 5.21 of
Dodd-Jensen \cite{DJ}. The Dodd-Jensen proof of their lemma uses the following consequence of the axiom of choice.
\begin{proposition}\label{proper} Assume \ $\ZF+ \AC$. \ Let \ $\kappa\ge\omega$ \ be a cardinal. Then
\be
\item $\kappa^+$ \ is a regular cardinal.
\item Suppose that \ $\kappa$ \ is a regular cardinal. Let \ $\langle X_\alpha : \alpha<\lambda\rangle$ \ be any
sequence of sets, where \ $\lambda<\kappa$. \ If \ $\abs{X_\alpha}<\kappa$ \ for all \ $\alpha<\lambda$, \ then \
$\abs{\bigcup\limits_{\alpha<\lambda}X_\alpha
}<\kappa$.
\ee
\end{proposition}

The two conclusions stated in the above proposition need \ $\AC$ \ for their proof. In fact, there are models of \ $\ZF$ \ in
which these conclusions are false. Thus, to guarantee that the proof of our generalization (of Dodd-Jensen's Lemma
5.21) does not implicitly appeal to the axiom of choice, we must prove Theorem \ref{card} below.
The proof of this theorem is established by means of a forcing argument. To show that this forcing argument does indeed
provide a proof, we make the following observations:  Let
\ $\psi$ \ be a sentence of set theory. Barwise \cite[Theorem 8.10]{Barwise}  proves in \ $\ZF$ \ that  if \ $\psi$ \ has
a transitive model \ (e.g., \ $V_\alpha$), then \ $\psi$ \ has a transitive model in \ $L$. \
Therefore, a version of the L\"owenheim--Skolem theorem is provable in \ $\ZF$ \ without the axiom of choice, namely, 
\[
\text{``$\psi$ has a transitive model''} \Rightarrow \text{``$\psi$ has a countable
transitive model''}.\]
Hence, if a sentence \ $\varphi$ \ is true in every countable transitive model of a sufficiently large finite fragment of \
$\ZF$, \ then it follows  that \ $\ZF \vdash \varphi$ \ from the reflection principle. In particular, suppose that \
$\gamma$ \ is a sentence. Then  if the sentence \ $\gamma\rightarrow\varphi$ \ is true in every countable
transitive model of a sufficiently large finite fragment of \ $\ZF$, \ then it follows  that \ 
$\ZF+\gamma \vdash \varphi$. \ This completes our discussion of the observations that are used implicitly
to show that \ $\ZF+ V\eq L(\R)[A]$ \ is strong enough to prove the conclusions of following theorem.

\begin{theorem}\label{card} Let \ $A\subseteq \R\times(\OR)^{<\omega}$. \ Assume \ $\ZF+ V\eq L(\R)[A]$. \ Then
\be
\item [(i)] For all sets \ $X$, \ the $V\!\!$--cardinal \  $|X|_V$  exists.
\item[(ii)] Let \ $\kappa\ge\Theta$ \ be a \ $V\!\!$--cardinal. Then \ $\kappa^+_V$ \ is a \ $V\!\!$--regular cardinal.
\item[(iii)] Suppose that \ $\kappa$ \ is a \ $V\!\!$--regular cardinal. Let \ $\langle X_{\langle\alpha, x\rangle} :
\langle\alpha, x\rangle\in
\lambda\times\R\rangle$
\ be any sequence of sets, where \ $\lambda<\kappa$. \ If \ $\abs{X_{\langle\alpha, x\rangle}}_V<\kappa$ \ for all \
$\langle\alpha, x\rangle\in
\lambda\times\R$,
\ then \
$\abs{\bigcup\limits_{\langle\alpha, x\rangle\in \lambda\times\R}X_{\langle\alpha, x\rangle}}_V<\kappa$.
\ee
\end{theorem}

\begin{proof} This follows immediately from Lemma \ref{choice} and Theorem \ref{cardM}.
\end{proof}

One observes the apparent paradox: 
{\sl We used the axiom of choice to prove that there is no need to use this axiom.}

\section{The fine structure of real mice}\label{finestructure}
This paper is the first in a series of three papers devoted to determining the minimal complexity of scales in the inner model
\ $\Kr$. \ In an effort to make this series self-contained we now give an overview of the fundamental
notions presented in \cite{Crcm} and \cite{Cfsrm} which we will assume here and in Parts II \& III. More
specifically, in subsections \ref{abovereals}--\ref{realmice} we will cover the relevant definitions,
concepts and theorems concerning (respectively): 
\be
\item fine structure above the reals,
\item iterable real premice, 
\item real 1--mice and the definition of $\Kr$, and
\item real mice.
\ee
In addition, we shall establish some additional results that will also be used in Parts II \& III. 
\subsection{Fine structure above the reals}\label{abovereals}

Let \ $m\in\omega$. \ For \ $N\in\omega$, \ the language 
\[\Lng_N=\lbrace\en,\underline{\R},c_1,\dots,c_m,A_1,\dots,A_N\rbrace\]
consists of the constant symbols \ $\underline{\R}$ \ and \
$c_1,\dots,c_m$ \ together with the membership relation \ $\in$ \ and the predicate symbols \ $A_1,\dots,A_N$. \ The theory \
$\text{R}_N$ \ is the deductive closure of the following weak set theory above the reals:
\newcounter{plist}
\begin{list}
{\upshape (\arabic{plist})}
{\usecounter{plist}}
\item $\forall x\forall y(x=y \leftrightarrow \forall z(z\in x
\leftrightarrow  z\in y))$ \hfill (extensionality)
\item $\exists y\forall x(x\notin y)$ \hfill ($\emptyset$ exists)
\item $\forall x(x \ne \emptyset \Rightarrow \exists y(y\in x \land 
x\cap y=\emptyset))$ \hfill (foundation)
\item $\forall x\forall y\exists z\forall t(t\in z \leftrightarrow 
(t=x \lor  t=y))$ \hfill (pairing)
\item $\forall x\exists y\forall t(t\in y \leftrightarrow \exists z(z\in x
\land t\in z))$ \hfill (union)
\item $\exists w(\emptyset \in w \land \text{ord}(w) \land \text{lim}(w) \land \forall
\alpha \in w \lnot \text{lim}(\alpha))$ \hfill ($\omega$ exists)
\item\label{sep} $\forall u\forall \vec x \ \exists
 z \forall s (s\in z \leftrightarrow s \in u
\land \psi (s,\vec x))$ \hfill ($\Sigma_0$ separation)
\item\label{clo} $\forall u\forall \vec x \ \forall
w\exists y\forall z(z\in y \leftrightarrow \exists t\in w(z=\{ s \in u 
 : \varphi (s,t,\vec x)\}))$ \hfill ($\Sigma_0$ closure)
\item $\forall x(x\in \underline{\R} \leftrightarrow \forall y\in
x\exists n\in \omega \exists z\exists f
(\text{Tr}(z) \land y \subseteq z \land f\colon n \maps{onto} z))$
\hfill ($\underline{\R} = V_{\omega+1}$)  
\end{list}
where, in (\ref{sep}) and (\ref{clo}), \ $\psi$ \ and \ $\varphi$ \ range over \ $\Sigma_0$ \
formulae. The above predicates \ $\text{ord}(w)$, \ $\text{lim}(w)$, and \ $\text{Tr}(w)$ 
abbreviate ``$w$ is an ordinal'', ``$w$ is limit'', and ``$w$ is transitive'', respectively.

The set of reals \ \R \ is a proper subset of \ $V_{\omega+1}$ \ and is
easily ``separated'' from \ $V_{\omega+1}$. \ It is more convenient, however, to start
constructing new sets from the transitive set \ $V_{\omega+1}$ \ rather than from
\ $\R$. \ Since \ $V_{\omega+1}$ \ can be ``constructed'' from \ $\R$, \ we shall
consider \ $V_{\omega+1}$ \ as given and we will tacitly identify the two.

Recall the basis functions \ $F_1,\dots,F_{17},F_{18},\dots,F_{17+N}$ \ of Dodd
\cite[Definition 1.3]{Dodd} where each \ $F_j$ \ is a function of two variables. 
In particular, \ $F_{17+i}(a,b) = a\cap A_i$.  
\begin{definition} A function \ $\mathcal{F}$
\ is {\it rudimentary\/} in  \ $A_1,\dots,A_N$ \ provided that \ $\mathcal{F}$ \ is a composition of
the basis functions. 
\end{definition}

\begin{definition} A relation \ $\mathcal{R}$ \ is {\it rudimentary\/} in 
$A_1,\dots,A_N$ provided that for some rudimentary function \ $\mathcal{F}$, \  
$\mathcal{R}(x_1,\dots,x_n) \Longleftrightarrow \mathcal{F}(x_1,\dots,x_n) \ne \emptyset$, for all
\ 
$x_1,\dots,x_n$.
\end{definition}

We shall often just say that the function \ $\mathcal{F}$ \ (or the relation \ $\mathcal{R}$) \ is 
rudimentary, when the predicates \ $A_1,\dots,A_N$ \  are clear from the context.

The theory \ $\Rtheory^+_N$ \ is the theory \ $\text{R}_N$ \ together with the \
$\Pi_2$ \ sentence \ $V=L[A_1,\dots,A_N](\underline{\R})$. \ When there is no ambiguity we shall write \ $\Rplus$ \ for \
$\Rtheory^+_N$. \ Also, the theory \ $\Rplus$ \ can be axiomatized  by a single \ $\Pi_2$ \ sentence of the language \
$\Lng_N$. Given models \ $\mouse$ \ and \ $\nouse$ \ for the language \ 
$\Lng_N$,  \ a map \ $\pi\colon \mouse \mapsigma{0} \nouse$ \ 
is said to be  {\it cofinal\/} if, for all \ $a\in N$, \ there is a \ $b\in M$ \ such that
\ $\nouse \models a\in \pi(b).$  \ In \cite{Crcm} we show that if \ $\pi\colon \mouse
\mapsigma{0} \nouse$ \ is cofinal where \ $\mouse$ \ is a model of \ $\Rplus$, \ then \ 
$\pi\colon \mouse \mapsigma{1} \nouse$ \ and \ $\nouse \models \Rplus$. \ In some cases, the
embedding \ $\pi$ \ also preserves \ $\Pi_2$ \ sentences. Specifically, let \ $\theta$ \ be a \
$\Pi_2$ \ formula, say \ $\forall u \exists
v\psi(u,v)$, \ where \ $\psi$ \ is in \ $\Sigma_0$. \ Let \ $s$ \ be an assignment of elements
in \ $M$ \ to the free variables of \ $\theta$. \ We say that \ $\theta[s]$ \ is \
$\mouse$--{\it collectible\/}  if \ $\mouse \models \theta[s]
\text{ \ if and only if \ } \mouse \models  (\forall a) (\exists b)(\forall u \in a)(\exists v \in
b)\psi(u,v)[s]$. \ Consequently, if \ $\pi\colon \mouse
\mapsigma{0} \nouse$ \ is cofinal and \ $\mouse\models \theta[s]$, \ then \ $\nouse\models
\theta[\pi(s)]$ \ (see \cite[Lemma 1.19]{Crcm}).

We are interested in transitive models \ $\mouse
= (M,\en,\underline{\R}^{\mouse},A_1,\dots,A_N)$ \ of \ $\Rtheory^+_N$. \ We shall write \
${\R}^{\mouse} = \underline{\R}^{\mouse}$ for \ \mouse's \ version of the reals. For any \
$\alpha \in \OR^{\mouse}$ \ we let \  $S_\alpha^{\mouse}(\underline{\R})$ \ denote the unique set
in \ $\mouse$ \  satisfying  \ $\mathcal{M} \models \exists f(\varphi(f) \land \alpha\in \dom(f)
\land  S_\alpha^{\mouse}(\underline{\R}) = f(\alpha)),$ \  
where \ $\varphi$ \ is the \ $\Sigma_0$ \ sentence used to define the sequence \ 
$\langle S_\gamma^{\mouse}(\underline{\R})\ : \ \gamma <
\OR^{\mouse} \rangle$ \ (see Definition 1.5 of \cite{Crcm}).
For \ $\lambda = \OR^{\mouse}$, \ let \ 
$S_\lambda^{\mouse}(\underline{\R})  = \bigcup\limits_{\alpha < \lambda} S_\alpha^{\mouse}(\underline{\R}).$
\   Let \  $\widehat{\OR}$ \ denote the class \ $\lbrace \gamma : \text{ the ordinal }
\omega \gamma \text{ exists} \rbrace$ \ and let \ 
$\Jgamma = S_{\omega\gamma}^{\mouse}(\underline{\R})$,\ for \ $\gamma \leq
\widehat{\OR}^{\mouse}$. \ Since \ $\mouse \models \text{R}^+$, \ it follows that
\  $M=J_\alpha^{\mouse}(\underline{\R})$ \ where \ $\alpha=\widehat{\OR}^{\mouse}$. \ 
Let \ $\mouse^\gamma$ \ be the
substructure of \ $\mouse$ \ defined by \ $\mouse^\gamma =
(\Jgamma,\en,\underline{\R}^{\mouse},\Jgamma \cap A_1,\dots,\Jgamma \cap A_N)$ \  
for \ $1<\gamma \leq \widehat{\OR}^{\mouse},$ \ and let \ 
$M^\gamma = \Jgamma$. \ We can write \
${\mouse}^\gamma = (M^\gamma,\en,\underline{\R}^{\mouse},A_1,\dots,A_N),$ \ as this will
cause no confusion. In particular, \ $\mouse^\gamma$ \ is {\it amenable}, that is, \ 
$a\cap A_i \in M^\gamma$, \ for all \ $a\in M^\gamma$ \  where \ $1\leq i\leq N$.

For \ $F,G \in [\OR]^{<\omega}$ \ let
\ $F<_{\BK} G \iff \exists \alpha \in G(G=F-\alpha) \lor max(G \bigtriangleup F)
\in G$. \ Here, \ $\bigtriangleup$ \ is the symmetric difference operation.
The order \ $<_{\BK}$ \ is the Brouwer-Kleene order on finite sets of ordinals and is a \
$\Sigma_0$ \ well-order.

\begin{lemma}\label{satisfaction}
For each \ $n,k\in\omega$ \ let \ $\langle\varphi_i \ : \ i\in
\omega \rangle$ \ be an effective enumeration of the \ $\Sigma_n$ \ formulae in the language \ $\Lng_N$ \ containing \
$k+1$ \ many free variables. \ Now let \ $\mouse$ \ be a transitive model of \ $\Rplus$. \ Define the
satisfaction relation \ $\textup{Sat}_n^k$ \ by
\[\textup{Sat}_n^k(x,a_1,\dots,a_k) \iff \mouse \models
\varphi_{x(0)}(\lambda i.x(i+1), a_1,\dots,a_k)\]
for each \ $x\in\R^\mouse$ \ and \ $a_1,\dots,a_k\in M$.
\ Then the relation \ $\textup{Sat}_n^k$ \ is \ $\Sigma_n(\mouse)$.
\end{lemma}

\begin{proof}
Because \ $\mouse$ \ is transitive and rudimentarily closed, Corollary 1.13 of \cite{Jensen} implies that 
the satisfaction relation \ $\textup{Sat}_n^k$ \ is \ $\Sigma_n(\mouse)$.
\end{proof}

\begin{definition}\label{skolemdefn} Suppose \ $\mouse$ \ is a transitive model of \ $\Rplus$. \
Let \ $h\colon {\R}^{\mouse}\times M \rightarrow M$ \ be a partial \ $\boldface{\Sigma}{n}(\mouse)$ \ map.
Then \ $h$ \ is a \ $\Sigma_n$ \ {\it Skolem function\/} if and only if whenever \ $S$ \ is
\ $\Sigma_n(\mouse,\{a\})$ \ for some \ $a\in M$, \ and \ $S\ne\emptyset$, \ then \ 
$\exists x\in {\R}^{\mouse} \ (h(x,a) \in S)$.
\end{definition}

\begin{definition} Suppose \ $\mouse$ \ is a transitive model of \ $\Rplus$. \
Then\ $\mouse$ \ satisfies \ $\Sigma_n$ \ {\it selection} if and only if whenever
$R \subseteq M\times [{\OR}^{\mouse}]^{<\omega}$ \ is \ $\boldface{\Sigma}{n}(\mouse)$, \ 
there is a \ $\boldface{\Sigma}{n}(\mouse)$ \ relation \ $S \subseteq R$ \ such that
\[\forall a(\exists F R(a,F)\Rightarrow \exists ! F S(a,F)).\]
\end{definition}

\begin{lemma} Let \ $\mouse$ \ be a transitive model of \ $\Rplus$ \ and suppose that 
 \ $\mouse$ \  satisfies \ $\Sigma_n$ \ selection. Then there is a \ $\Sigma_n$ \ Skolem
function for \ $\mouse$. 
\end{lemma}

\begin{lemma}\label{skolem1} Let \ $\mouse$ \ be a transitive model of \ $\Rplus$. \ Then
there is a canonical \ $\Sigma_1$ \ Skolem function for \ $\mouse$. 
\end{lemma}

\begin{definition}\label{skolemdef} Let \ $\mouse$ \ be a transitive model of \ $\Rplus$. \ 
Then \ $h_{\mouse}$ \ is the canonical \ $\Sigma_1$ \ Skolem function with the parameter
free \ $\Sigma_1$ \ definition established by (the proof of) Lemma \ref{skolem1}. For \ $p\in M$, \
$h_{\mouse}^p$ \ is the function given by \ $h_{\mouse}^p(x,s) = h_{\mouse}(x,\langle s,p\rangle).$ 
\end{definition}

\begin{definition}\label{hulldef} Let \ $\mouse =
(M,\en,\underline{\R},c_1,\dots,c_m,A_1,\dots,A_N)$ \ be a transitive \ ${\mathcal L}_N$ \ model
and let \ ${\R}^{\mouse}\subseteq X\subseteq M$. \ Let \ $H = \{a\in M : \{a\}
\text{ \ is \ } \Sigma_n(\mouse,X)\}$. \ The \ $\Sigma_n$ \ {\it hull} \ of \ $X$ \ is the 
substructure \[{\Hull}_n^{\mouse}(X) = 
(H,\en,{\R}^{\mouse},c_1,\dots,c_m,H\cap A_1,\dots,H\cap A_N).\] We shall write \
${\Hull}_\omega^{\mouse}(X) =\bigcup\limits_{n\in\omega}{\Hull}_n^{\mouse}(X)$.  
\end{definition}

\begin{lemma}\label{skolemimage} Let \ $\mouse$ \ be a transitive model of \ $\Rplus$ \ and let \ 
$\mathcal{H} = {\Hull}_1^{\mouse}({\R}^{\mouse}\cup X\cup \{p\})$, \ where \ $X\cup
\{p\}\subseteq M$ \ and \ $X\ne\emptyset$. \ Then 
\be
\item $\mathcal{H} \prec_1 \mouse \text{\ and \ }\mathcal{H} \models \Rplus$ 
\item $H = {h_{\mouse}^p}^{\prime\prime}({\R}^{\mouse}\times(X)^{<\omega})$. 
\ee
\end{lemma}

\begin{definition} Let \ $\mouse$ \ be a transitive model of \ $\Rtheory^+$. 
 An ordinal \ $\lambda\le{\OR}^\mouse$ \ is a \ $\boldface{\Sigma}{1}(\mouse)$-{\it cardinal}
\ if for no \ $\gamma<\lambda$ \ does there exist a partial \ $\boldface{\Sigma}{1}(\mouse)$ \ function
\  $f\colon \gamma \times {\R}^\mouse\maps{onto} \lambda$.
\end{definition}

\begin{definition}\label{rhoref} Let \ $\mouse$ \ be a transitive model of \ $\Rplus$. \ 
The {\it projectum} \ $\rho^{}_\mouse$ \ is the least ordinal \ $\rho \le
\widehat{\OR}^{\mouse}$ \ such that \  $\pow(\R^{\mouse}\times
\omega\rho)\cap\boldface{\Sigma}{1}(\mouse)\not\subseteq M$, \ and \ $p^{}_\mouse$ \ is the \ 
$\le_{BK}$--least \ $p\in[\OR^{\mouse}]^{<\omega}$ \ such that \  $\pow(\R^{\mouse}\times
\omega\rho^{}_\mouse)\cap\Sigma_1(\mouse,\{p\})\not\subseteq M$. 
\end{definition}

\begin{lemma}\label{sigmacard} Let \ $\mouse$ \ be a transitive model of \ $\Rtheory^+$. 
\be
\item $\omega\rho^{}_\mouse$ \ is a \ $\boldface{\Sigma}{1}(\mouse)$-cardinal
\item $\omega\rho^{}_\mouse = \rho^{}_\mouse$
\item $\omega\rho^{}_\mouse$ \ is closed under the G\"odel pairing function (see \cite[p. 943]{Cfsrm}).
\ee
\end{lemma}

\begin{definition}\label{master} Let \ $\mouse$ \  be a transitive model of \ $\Rplus$. \ The \ $\Sigma_1$--{\it master
code\/} \ $A_\mouse$ \ of \ $\mouse$ \  is the set 
\[A_\mouse = \{(x,s)\in \R^\mouse\times (\omega\rho^{}_\mouse)^{<\omega} : \mouse\models
\varphi_{x(0)}[\lambda n.x(n+1),s,p^{}_\mouse]\}\]
where \ $\langle\varphi_i : i\in\omega\rangle$ \ is a fixed recursive listing of all the \
$\Sigma_1$ \ formulae of three variables in the language \ $\Lng_N$.
\end{definition}

\subsubsection{Acceptability above the reals}
\begin{definition}\label{acceptable} Suppose that \ $\mouse$ \  is a transitive model of \
$\Rtheory^+_N$.
\ We say that \ $\mouse$ \  is {\it acceptable above the reals\/} provided that for all \ $\gamma<{\widehat{\OR}}^\mouse$ \
and \ $\delta<{\OR}^{\mouse^\gamma}$, \
if \ $\mathcal{P}(\delta\times{\R}^{\mouse})\cap M^{\gamma+1}\nsubseteq M^\gamma$, \ then for each \ $u\in M^{\gamma+1}$ \ there is
an \ $f\in M^{\gamma+1}$ \ such that
\ $f=\langle f_{\xi} : \delta\le \xi<{\OR}^{\mouse^\gamma}\rangle$ \ 
and
\begin{equation*}f_{\xi}\colon \xi\times\R^\mouse \maps{onto}
\{\xi\}\,\cup\,(\mathcal{P}(\xi\times\R^\mouse)\cap u).
\end{equation*} 
\end{definition}

\begin{rmk}When the context is clear we will say that \ $\mouse$ \ is ``acceptable'' rather than say that  \
$\mouse$ \ is ``acceptable above the reals''.
\end{rmk}

\begin{rmk} Definition \ref{acceptable} is essentially the same as our definition of acceptability in \cite[Definition
1.1]{Cfsrm}. The definition in \cite{Cfsrm} involves functions \ $f_{\xi,x}$ \ with an additional index \
$x\in\R$. \ These two definitions are easily shown to be equivalent;
however, we now feel that the above Definition \ref{acceptable} is the more pertinent version. Thus from now on, we will
consider Definition \ref{acceptable} as the ``official definition'' of acceptability above the reals.
\end{rmk}

Dodd-Jensen \cite{DJ} first defined when a premouse without the reals is {\it strongly
acceptable.\/} Their definition assumes that such a premouse satisfies the axiom of choice. In our case,
however, a transitive model \ $\mouse$ \ of \ $\Rtheory^+_N$ \ does not necessarily satisfy the axiom
of choice. So before we can define when a real premouse is {\it strongly
acceptable above the reals,\/} we must modify the standard definition
of cardinality in \ $\mouse$.

\begin{definition}\label{realcard} Let \ $\mouse$ \ be a transitive model of \ $\Rtheory^+_N$. 
\be
\item For \ $a\in M$
\ the \ $\mouse$--{\it cardinality} \ of \ $a$, \ denoted by \ $|a|_\mouse$, \ is the least ordinal \
$\lambda\in{\OR}^\mouse$ \ such that \ $f\colon \lambda \times
{\R}^\mouse\maps{onto} a$ \ for some \ $f\in M$.
\item An ordinal \ $\lambda\le{\OR}^\mouse$ \ is an \ {\it $\mouse$--cardinal} \ if \ 
$\lambda = |\lambda|_\mouse$ \ or \ $\lambda = {\OR}^\mouse$.
\item For an ordinal \ $\lambda<{\OR}^\mouse$, \ $\lambda^+_\mouse$ \ is the least \
$\mouse$--cardinal greater than \ $\lambda$.
\ee
\end{definition}

\begin{rmk} We have decided to use the cross
product \ $\lambda\times{\R}^\mouse$, \ in the above definition of $\mouse$--{cardinality}. In the special case where \
$\lambda=0$, the   cross product \ $0\times{\R}^\mouse=\emptyset$ \ 
does not involve any reals. So, when applying the above definition, we shall abuse 
cross product notation slightly and define \ $0\times{\R}^\mouse={\R}^\mouse$.
\end{rmk} 

When \ $\mouse$ \  is a transitive model of \
$\Rtheory^+_N$, \ Lemma 1.7 of \cite{Crcm} implies that for any \ $a\in\mouse$ \ there is a function \ $f\in\mouse$
\ such that \ $f\colon [\alpha]^{<\omega}\times{\R}^\mouse \maps{onto} a$, \ for some \
$\alpha<{\OR}^\mouse$. \ Thus, Lemma 1.4 of \cite{Cfsrm} shows
that the \ $\mouse$--cardinality of a set in \ $\mouse$ \ always exists.

\begin{definition}\label{SA} Suppose that \ $\mouse$ \ is a transitive model of \
$\Rtheory^+_N$. \ We say that \ $\mouse$ \ is {\it strongly acceptable above the reals\/} if, whenever \
$\gamma\in{\widehat{\OR}}^\mouse$, \ $\delta<{\OR}^{\mouse^\gamma}$ \ and \ 
$\pow(\delta\times{\R}^{\mouse})\cap M^{\gamma+1}\nsubseteq M^\gamma$, \ then \
$\abs{\omega\gamma}_{\mouse^{\gamma+1}}\le\delta$. 
\end{definition}

\begin{rmk}When the context is clear, we may say that \ $\mouse$ \ is ``strongly acceptable'' rather than say that \ $\mouse$
\ is ``strongly acceptable above the reals''.
\end{rmk}

\begin{definition}\label{trans}
Let \ $\mouse$ \ be a transitive model of \ $\Rtheory^+_N$. \ Let \ $\xi\le{\OR}^{\mouse}$ \ and define \
$H^{\mouse}_{\xi}=\{a\in M :  \abs{T_c(a)}_\mouse<\xi\}$,
\ where \ $T_c(a)$ \ denotes the transitive closure of \ $a$.
\end{definition}

Let \ $\mouse$ \ be a transitive model of \ $\Rtheory^+_N$. \ Lemma 1.9 of
\cite{Crcm} states that \ $\mouse\models \forall a \exists y (y=T_c(a))$ \ and thus, 
 \ $\mouse\models \forall a \exists\lambda(\lambda=\abs{T_c(a)}_\mouse)$. \ Thus, the definition of \ $H^{\mouse}_{\xi}$ \ is well-defined.

We shall write \ $\langle J_\xi^{(\mouse, A_\mouse)}(\underline{\R}) :
\xi\in{\OR}\rangle$ \ for the Jensen hierarchy of sets which are relatively constructible above \ $\R^\mouse$ \ from the predicates \ $A_1\cap M,\dots,A_N\cap M,A_\mouse\cap M$. 
\begin{lemma}\label{kappa_small} Suppose that \ $\mouse$ \  is acceptable. Then \ $H^{\mouse}_{\omega\rho^{}_\mouse} = 
J_{\rho^{}_\mouse}^{(\mouse, A_\mouse)}(\underline{\R})$.
\end{lemma}
\begin{definition}\label{master2} Given that \ $\mouse$ \  is acceptable (above the reals), let \
$M^*=J_{\rho^{}_\mouse}^{(\mouse, A_\mouse)}(\underline{\R})$. \ The \ $\Sigma_1$-code of \ $\mouse$ \  is the
structure \ $\mouse^* = (M^*, \in,\underline{\R},c_1,\dots,c_m, A_1\cap M^*,\dots,A_N\cap M^*,A_\mouse\cap M^*)$ \
where we are assuming that the constants have the same interpretation in \ $\mouse^*$ \ as in \ $\mouse$. 
\end{definition}
\begin{lemma} $\mouse^*$ \ is strongly acceptable.
\end{lemma}
\begin{lemma} If \ $\mouse$ \ is
strongly acceptable, then \ $\mouse$ \ is
acceptable.
\end{lemma}
\begin{definition}\label{master3} Suppose that \ $\mouse$ \  is acceptable. Inductively define on \
$n\in\omega$ \ the \ $\Sigma_n$-code \ of \ \mouse, \ denoted by \ $\mouse^n$, \ as follows:
\be 
\item $\mouse^0 = \mouse$, \ $\rho_\mouse^0 = {\OR}^\mouse$, \ $p_\mouse^0=\emptyset$, \
and \ $\mathcal{A}_\mouse^0 = \emptyset$.
\item Assume that \ $\mouse^n$ \ has been defined and that \ $\rho^{}_{{\mouse^n}} > 1$. \ Define \ 
$\mouse^{n+1} = (\mouse^n)^*$, \ $\rho_\mouse^{n+1} = \rho^{}_{\mouse^n}$, \
$p_\mouse^{n+1}=p^{}_{\mouse^n}$, \ and \ $\mathcal{A}_\mouse^{n+1} = A_{\mouse^n}$. 
\ee
\end{definition}

One can show that \ $\rho_\mouse^i=\omega\rho_\mouse^i$ \ when \ $i\ge 1$.

\begin{rmk} The above notation is slightly inconsistent with previous notation. Namely, \ 
$\mouse^\gamma$ \ denotes a substructure of \ $\mouse$, \ while \ $\mouse^n$ \ denotes the \
$\Sigma_n$-code \ of \ $\mouse$. \ Nevertheless, we shall use integers and integer variables, 
for example \ $n$, \ exclusively for denoting \ $\mouse^n$, \ the \ $\Sigma_n$-code \ of \ $\mouse$,
\ and thereby resolve any confusion in notation.
\end{rmk}

\begin{definition}\label{sounddef} Let \ $\mouse$ \ be acceptable. We say that
\be
\item $\mouse$ \ is \ {\it sound} \ if
\ $h_{\mouse}^{p^{}_\mouse}\colon \R^{\mouse}\times\omega\rho^{}_\mouse \maps{onto} M$
\item $\mouse$ \ is \ {\it $n$--sound} \ if \ $\mouse^i$ \ is sound for all \ $i<n$.
\ee
\end{definition}

\begin{lemma}\label{soundlemma} Suppose that \ $\mouse$ \ is sound. Then
\be
\item $\rho^{}_\mouse$ \ is the least \ $\gamma\le{\widehat{\OR}}^\mouse$ \ so that there
is a \ $p\in[\OR^\mouse]^{<\omega}$ \ such that \ $h_{\mouse}^{p}\colon
\R^{\mouse}\times\omega\gamma \to M$  \ is onto
\item ${p^{}_\mouse}$ \ is the \ $\le_{\BK}$--least \ $p\in[\OR^\mouse]^{<\omega}$ \ such that 
\ $h_{\mouse}^{p}\colon\R^{\mouse}\times\omega\rho^{}_\mouse \to M$ \ is onto.
\ee
\end{lemma}

We now consider another way of iterating a ``projectum.''
\begin{definition}\label{defproj} Suppose that \ $\mouse$ \ is acceptable. Let
\ $\gamma_\mouse^0 = {\OR}^\mouse$ \ and, for \ $n\ge 1$, \ define \ $\gamma_\mouse^n$ \ to
be the least ordinal \ $\gamma\le\widehat{\OR}^{\mouse}$ \ such that \  $\pow(\R^{\mouse}\times
\omega\gamma)\cap\boldface{\Sigma}{n}(\mouse)\not\subseteq M$. 
\end{definition}

One can also show that \ $\gamma_\mouse^i=\omega\gamma_\mouse^i$ \ when \ $i\ge 1$.
\ For an arbitrary acceptable \ $\mouse$ \ the connection between \ $\gamma_\mouse^n$ \ and \ 
$\rho_\mouse^n$ \ is not clear; however, if \ $\mouse$ \ is $n$-sound, then   for all \ $i\le n$, \
$\gamma_\mouse^i=\rho_\mouse^i$ \ and \ $\mouse$ \ satisfies \ $\Sigma_{i+1}$ \ selection (see below). 

The statement and proof of Lemma 4.19 of Dodd \cite{Dodd} carry over to give our next result.
Recall that if \ $\mouse$ \ is a model of the language \ $\Lng_N$, \ then \ $\mouse^1$ \ is 
a model of the language \ $\Lng_{N+1}$.

\begin{lemma}\label{secgoose} Suppose \ $\mouse$ \ is sound and \ $\rho^{}_\mouse>1$. \ For \ 
$B\subseteq \R^{\mouse}\times(\omega\rho^{}_\mouse)^{<\omega}$,
\[B \text{ is } \boldface{\Sigma}{k}(\mouse^1) \iff  
B \text{ is }
\Sigma_{k+1}(\mouse,\R^{\mouse}\times(\omega\rho^{}_\mouse)^{<\omega}
\cup\{p^{}_\mouse\})\]
for all \ $k\ge 1$.
\end{lemma}

\begin{proof} See the proof of Lemma 4.19 of \cite{Dodd}.
\end{proof}

Iterating Lemma \ref{secgoose} gives

\begin{lemma}\label{thdgoose} Suppose \ $\mouse$ \ is  \ $n$--sound and \ $\rho_\mouse^n>1$. \ For \ 
$B\subseteq \R^{\mouse}\times(\omega\rho_\mouse^n)^{<\omega}$,
\[B \text{ is } \boldface{\Sigma}{k}(\mouse^n) \iff  
B \text{ is }
\Sigma_{k+n}(\mouse,\R^{\mouse}\times(\omega\rho_\mouse^n)^{<\omega}
\cup\{p_\mouse^1,\dots p_\mouse^n\})\]
for all \ $k\ge 1$.
\end{lemma}

\begin{corollary}\label{terms} Suppose \ $\mouse$ \ is \ $n$--sound. Then for \ $1\le i\le n$ \ there
is a  \ $\Sigma_i(\mouse,\R^{\mouse}\times(\omega\rho_\mouse^i)^{<\omega}
\cup\{p_\mouse^1,\dots p_\mouse^i\})$ \ partial function \ 
$f\colon \R^{\mouse}\times\omega\rho_\mouse^i \maps{onto} M$.
\end{corollary}

\begin{corollary}\label{gammaproj} If \ $\mouse$ \ is \ $n$--sound, then \
$\gamma_\mouse^i=\rho_\mouse^i$ \  for all \ $i\le n$.
\end{corollary}

\begin{corollary} Suppose \ $\mouse$ \ is \ $n$--sound and \ $\rho_\mouse^n > 1$. \ 
Then
\be
\item $M^n\cap \pow(\lambda \times \R^\mouse) = M\cap \pow(\lambda \times
\R^\mouse)$, \ for all \ $\lambda < \omega\rho_\mouse^n$
\item $\boldface{\Sigma}{k}(\mouse^n)\cap \pow(\omega\rho_\mouse^n \times \R^\mouse) =
\boldface{\Sigma}{k+n}(\mouse)\cap\pow(\omega\rho_\mouse^n \times\R^\mouse)$, \ 
for all \ $k\ge 1$.
\ee
\end{corollary}

\begin{theorem} Suppose \ $\mouse$ \ is \ $n$--sound and \ $\rho_\mouse^n > 1$. \ 
Then \ $\mouse$ \ satisfies \ $\Sigma_{i+1}$ \ selection, \ for all \ $i\le n$.
\end{theorem}

\begin{corollary} If \ $\mouse$ \ is \ $(n+k)$--sound and \ $\rho_\mouse^{n+k} > 1$,
\  then \ $\mouse^n$ \  satisfies \ $\Sigma_{k+1}$ \ selection.
\end{corollary}

The arguments which establish that \ $\gamma_\mouse^i=\rho_\mouse^i$ \ and that \ $\mouse$ \ 
satisfies  \ $\Sigma_{i+1}$ \ selection, \ for all \ $i\le n$, \ use the \ $n$--soundness of
\ $\mouse$. \ It should be noted that there are models \ $\mouse$ \ which are $n$--sound but not 
$(n+1)$--sound. In these cases the equality of \ $\gamma_\mouse^{n+1}$ \ and \
$\rho_\mouse^{n+1}$ \ is questionable, although the equality does hold for {\it real mice\/} (see Theorem \ref{gamma}).

\subsubsection{Projected types}\label{projtypes} We will show that a \ $\Sigma_{n+1}$--type \ $\Sigma$ \ (in the
language \ $\Lng_0$) can be ``translated'' to a \ $\Sigma_{1}$--type \ ${\Sigma}^*$ \ (in the language \ $\Lng_n$) such
that an \ $n$--sound \ $\mouse$ \  realizes \ $\Sigma$  \ if and only if \  $\mouse^n$ \ realizes \
${\Sigma}^*$.
\ We do this by means of two lemmas. The
first lemma makes the observation that the proof of the direction ($\Longleftarrow$) in Lemma \ref{secgoose} is uniform in \
$\mouse$, \ a transitive model of \ $R_N^+$. \ The second lemma will be used (in conjunction with Lemma
\ref{iteratedtypes} below) in the proof of Theorem
\ref{newthm} in Part II.

\begin{definition} Let \ $n,N\in\omega$. \ Then 
\begin{align*}
\Sigma^{N}_{n} &= \{ \theta\in \Sigma_n : \text{$\theta$ is an $\Lng_N$ formula of one free variable}\}\\
\Pi^{N}_{n} &= \{ \theta\in \Pi_n : \text{$\theta$ is an $\Lng_N$ formula of one free variable}\}\\
\Upsilon^{N}_{n} &= \Sigma^{N}_{n}\,\cup\,\Pi^{N}_{n}.
\end{align*}
\end{definition}

In the remainder of this subsection, we will be presuming that \ $\mouse$ \ is an \ $\Lng_N$--model. 
\begin{lemma}\label{wow} For every \ $k\ge 1$ \ and \ $N\ge 0$, \ there is a map \
$\overline{\phantom{\theta}}\colon \Sigma^{N}_{k+1}\to\Sigma^{N+1}_{k}$
\ (i.e., $\theta(v)\mapsto\overline{\theta}(v)$)  \ where for any sound \ $\mouse$ \ with \ $\rho^{}_\mouse>1$ \ the following
holds: for each \ $a\in M$ \ there exists a \ $q\in M^1$ \ such that
\begin{equation}
\mouse\models\theta(a) \iff  \mouse^1\models\overline{\theta}(q), \text{ \ for all \ $\theta(v)\in\Sigma^{N}_{k+1}$.} 
\label{typechange}\end{equation}
In addition, for each \ $q\in M^1$ \ there exists an \ $a\in M$ \ such that {\rm (\ref{typechange})} holds.
\end{lemma}

\begin{proof}[Sketch of Proof] We shall assume that the reader is familiar with the proof of Lemma 4.19 in
\cite{Dodd}. Now, given \ $\theta(v)\in\Sigma_{k+1}$ \ in the language \ $\Lng_N$, \ the following procedure for
obtaining \ $\overline{\theta}(v)$ \ is uniform in \ $\mouse$. \ First define the relation \ $B=\{ q\in
\R^\mouse\times(\omega\rho^{}_\mouse)^{<\omega} : \mouse\models \theta(h_{\mouse}^{p^{}_\mouse}(q))\}$. \ The proof of
the direction ($\Longleftarrow$) of Lemma \ref{secgoose} gives a \ $\Sigma_k$ \ formula \ $\overline{\theta}(v)$ \ in the
language \ $\Lng_{N+1}$ \ with no parameters such that
\[q\in B\iff \mouse^1\models \overline{\theta}(q).\] Recall that for arbitrary \ $a\in M$ \
there is a \ $q\in M^1$ \ such that \ $a=h_{\mouse}^{p^{}_\mouse}(q)$, \ because \ $\mouse$ \ is sound. Now, since \
$\overline{\theta}(v)$ \ depends only on 
\ $\theta(v)$ \ it follows that the map \ $\theta(v)\mapsto\overline{\theta}(v)$  \ is as required, that is, \ 
$(\forall a\in M)(\exists q\in M^1)(\forall \theta\in\Sigma^{N}_{k+1})\left[\,\mouse\models \theta(a) \iff
\mouse^1\models
\overline{\theta}(q)\,\right]$. \ In addition, \ $(\forall q\in M^1)(\exists a\in M)(\forall
\theta\in\Sigma^{N}_{k+1})\left[\,\mouse\models \theta(a) \iff \mouse^1\models
\overline{\theta}(q)\,\right]$.
\end{proof}

\begin{lemma}\label{wowwee} For all \ $k\ge 1$ \ and \ $N,n\ge 0$, \ there is a map \
$\overline{\phantom{\theta}}\colon \Sigma^{N}_{k+n}\to\Sigma^{N+n}_{k}$
\ (i.e.,  $\theta(v)\mapsto\overline{\theta}(v)$) \ where for any $n$--sound \ $\mouse$ \ with \ $\rho_\mouse^n>1$ \ the
following holds: for each \ $a\in M$ \ there exists a \ $q\in M^1$ \ such that
\begin{equation}
\mouse\models\theta(a)\iff\mouse^n\models\overline{\theta}(q), \text{ \ for all \
$\theta(v)\in\Sigma^{N}_{k+n}$.}\label{typen}
\end{equation}
In addition, for each \ $q\in M^n$ \ there exists an \ $a\in M$ \ such that {\rm (\ref{typen})} holds.
\end{lemma}

\begin{proof} The proof follows from Lemma \ref{wow} by induction on \ $n$.
\end{proof}

\begin{corollary}\label{wowweecor} For all \ $k,n\ge 1$ \ there is a map \
$\overline{\phantom{\theta}}\colon \Upsilon^{N}_{k+n}\to\Upsilon^{N+n}_{k}$
\ (i.e.,  $\theta(v)\mapsto\overline{\theta}(v)$) \ where for any $n$--sound \ $\mouse$ \ with \ $\rho_\mouse^n>1$ \ the
following holds: for each \ $a\in M$ \ there exists a \ $q\in M^1$ \ such that
\begin{equation}
\mouse\models\theta(a)\iff\mouse^n\models\overline{\theta}(q), \text{ \ for all \
$\theta(v)\in\Upsilon^{N}_{k+n}$.}\label{typenn}
\end{equation}
In addition, for each \ $q\in M^n$ \ there exists an \ $a\in M$ \ such that {\rm (\ref{typenn})} holds.
\end{corollary}

\begin{definition}\label{deftype} An \ $\Upsilon^{N}_{n}$--type is any nonempty subset of \ $\Upsilon^{N}_{n}$ \ where \
$N,n\in\omega$.
\end{definition}

\begin{definition}\label{relztype} Let \ $\Upsilon$ \ be an \ $\Upsilon^{N}_{n}$--type. \ Suppose that \ $\mouse$ \ is
a transitive model of \ $R_N^+$. \ We say that \ $\mouse$ \ {\it realizes\/} \ $\Upsilon$ \ if there is an \ $a\in M$ \
such that \ $\mouse\models\theta(a)$ \ for all \ $\theta\in\Upsilon$. 
\end{definition}

\begin{corollary}\label{cortype} Suppose that \ $\mouse$ \ is $n$--sound and \ $\rho_\mouse^n>1$. \ Then any 
\ $\Upsilon^{N}_{k+n}$--type \ $\Upsilon$ \ can be translated to an \ $\Upsilon^{N+n}_{k}$--type \ $\overline{\Upsilon}$
\ such that
\ $\mouse$ \ realizes \ $\Upsilon$ \ if and only if \ $\mouse^n$ \ realizes \ $\overline{\Upsilon}$.
\end{corollary}

\subsection{Iterable real premice}\label{realpremice}

In this subsection we bring together fine structure above the reals and the theory
of iterated ultrapowers. The mixture of these two techniques produces iterable
real premice and allows us in \cite{Crcm} to construct scales beyond those in \ $\Lr$.
A real premouse \ $\mouse$ \ is a premouse in the usual sense (see -- \cite{DJ}) but with two
additional  conditions: (i) $\mouse$  contains the set of reals \ $\R$ \ as an
element and (ii) $\mouse$  believes that its measure is  ``$\R$--complete.''

\begin{definition} Let \ $\mu$ \ be a normal measure on \
$\kappa$. \ We say that \ $\mu$ \ is an {\it $\R$--complete measure} on \ $\kappa$ \ 
if the following holds: if \ $\langle A_x : x\in \R\rangle$ \ is any sequence such that \ 
$A_x\in \mu$ \ for all \ $x\in \R$, \ then \ $\bigcap_{x\in \R}A_x \in \mu$.
\end{definition}

We now focus our attention on
transitive models \ $\mouse$ \ of \ $\Rplus$ \ such that  \ $\mouse$ \  believes that one of its
predicates is an \ ${\R}^{\mouse}$--complete measure on \ ${\pow(\kappa)} \cap
M$. \ For this reason we modify our official language by letting \[\Lng_n =
\{\,\en,\underline{\R},\underline{\kappa},\mu, A_1,\dots,A_n\,\},\] where \ $\mu$ \
is a new predicate symbol and  \ $\underline{\kappa}$ \ is a constant symbol. 
Models of the language \ $\Lng_0 = \{\,\en,\underline{\R},\underline{\kappa},\mu\,\}$ \ will 
be our main interest. Finally, we let \ $\Lng_n^{\underline{p}} = \Lng_n\cup\{\underline{p}\}$ \ when 
we need to add \ $\underline{p}$ \ as a new constant symbol.

\begin{definition}[Premice]\label{premice} A model \ $\mouse = 
(M,\en,\underline{\R}^{\mouse},\underline{\kappa}^{\mouse},\mu,A_1,\dots,A_N)$ \ is
a \ {\it premouse} (above the reals) \ if 
\be
\item $\mouse$ \ is a transitive model of \ $\Rplus$
\item $\mouse \models \text{``$\mu$ is an \ $\underline{\R}$--complete \ measure on 
\ $\underline{\kappa}$''}$. 
\ee
$\mouse$ \ is a \ {\it pure premouse} \ if \ $\mouse = 
(M,\underline{\R}^{\mouse},\underline{\kappa}^{\mouse},\mu)$. \ Finally,  \ $\mouse$ \  is a \ {\it
real premouse} \  if it is pure and \ ${\R}^{\mouse}= \R$.
\end{definition}

Note that ``$\mu$ is an \ $\underline{\R}$--complete \ measure on \
$\underline{\kappa}$'' \ is a \ $\Pi_1$ \ assertion.
\begin{definition}\label{PMdefined} The theory  \ $\PM$ \  is the theory  \ $\Rplus$ \  together
with the sentence \ ``$\mu$ is an \ $\underline{\R}$--complete \ measure on \
$\underline{\kappa}$''.
\end{definition}

The theory \ $\PM$ \ can be
axiomatized by a single \ $\Pi_2$ \ sentence. For a premouse \ $\mouse$, \ we shall
write \ $\kappa$ or $\kappa^{\mouse}$, \ for \
$\underline{\kappa}^{\mouse}$ \  when the context is clear.  
We may refer to the \ ``pointclass\footnote{See Remark \ref{Rparameter}.} \ $\Sigma_n(\mouse,X)$'',
\ or assert that \ ``$\Sigma_n(\mouse,X)$ \ has the scale property.''  Both cases actually refer to
\ $\Sigma_n(\mouse,X) \cap \mathcal{P}(\R)$, \ but the context should make this clear.  We may also
say that \ $a\in \mouse$ \ when we mean \ $a\in M$. \ Finally, to distinguish our definition of a
premouse from the premice of Dodd-Jensen, we may sometimes refer to our version as ``premice above
the reals.''

\subsubsection{Premouse iteration}
Given a premouse  \ $\mouse$ \  we now define its ultrapower, denoted by \ $\mouse_1$. \ 
Let \ ${^\kappa M}=\{\,f\in M : f\colon \kappa \rightarrow M \,\}$.
For $f,g\in {^\kappa M}$ \ define 
\[f\sim g \iff \mouse \models \{\, \xi \in \kappa : f(\xi)=g(\xi) \,\} \in\mu.\] 
Since  \ $\mouse$ \  satisfies \ $\Sigma_0$ \ separation, the above
set is in \mouse, and \ $\sim$ \ is an equivalence relation on \
${^\kappa M}$. For \ $f\in {^\kappa M}$, \ we denote the equivalence class of \
$f$ \ by \ $[f]$.  Let \ $M_1 = \{\,[f] : f\in {^\kappa M}\,\}$ \ and define
\begin{alignat*}{1}
\text{$[f]$} \ E \ [g] &\iff \mouse \models\{\, \xi \in \kappa : f(\xi)\in g(\xi) \,\} \in \mu\\
\text{$[f]$} \in \mu_1 &\iff \mouse \models\{\, \xi \in \kappa : f(\xi)\in \mu \,\} \in \mu\\
\text{$[f]$} \in A_i^1 &\iff \mouse \models \{\, \xi \in \kappa : f(\xi)\in A_i \,\} \in
\mu\text{, \ where } 1\leq i \leq N.
\end{alignat*}
By amenability, the sets on the right hand side are in  \ $\mouse$ \  and therefore
can be measured by \ $\mu$. \ For \ $a\in M$, \ let \ $c_a\ \in {^\kappa M}$ \ be the
constant function defined by \ $c_a(\xi)=a$ \ for all $\xi \in \kappa$. Now
define 
\[{\mouse}_1 = (M_1,E,[c_{_{{\R}^{\mouse}}}],[c_{{\kappa}^{\mouse}}],\mu_1,A_1^1,\dots,A_N^1).\]
Since the meaning will always be clear, we usually write
\[{\mouse}_1 = (M_1,\en,\underline{\R}^{{\mouse}_1},\underline{\kappa}^{{\mouse}_1},\mu,A_1,\dots,A_N).\]
A version of \L o\'s' Theorem holds for this ultrapower without the axiom of choice.
\begin{theorem}\label{Los} Let  \ $\mouse$ \  be a premouse. Then
\[{\mouse}_1 \models \varphi([f_1],\dots,[f_n]) \iff 
\mouse \models \{\,\xi\in\kappa : \varphi(f_1(\xi),\dots,f_n(\xi))\,\}\in \mu,\]
for every \ $\Sigma_0$ \ formula \ $\varphi$ \ and for all $f_1,\dots,f_n \in
{^\kappa M}$. 
\end{theorem}

For a premouse  \ $\mouse$ \  define \ $\pi^{\mouse}\colon \mouse \rightarrow 
{\mouse}_1$ \ by \ $\pi^{\mouse}(a)=[c_a]$ \ for \ $a\in M$. \ When the context
is clear we shall omit the superscript and write \ $\pi$ \ for \ $\pi^{\mouse}$.
\begin{lemma}
Let  \ $\mouse$ \  be a premouse. Then \ ${\mouse}_1 \models \PM$.
\end{lemma}

In general, we can iterate this ultrapower operation and get a
commutative system of models by taking direct limits at limit ordinals.
\begin{definition}\label{pmiteration}  Let  \ $\mouse$ \  be a premouse. Then 
\begin{equation}\premousesystem\label{lablepmiteration}\end{equation}
is the commutative system satisfying the inductive definition:
\be
\item $\mouse_0 = \mouse$
\item $\pi_{\gamma\gamma} = \text{identity map, \ and } \pi_{\beta\gamma} \circ
\pi_{\alpha\beta} = \pi_{\alpha\gamma} \text{ \ for all } 
\alpha\leq\beta\leq\gamma\leq\lambda$
\item If \ $\lambda=\lambda^\prime + 1$, \ then \ $\mouse_\lambda =
\text{ultrapower of } \mouse_{\lambda^\prime}$, \ and \ 
$\pi_{\alpha\lambda}=\pi^{\mouse_{\lambda^\prime}} \circ
\pi_{\alpha{\lambda^\prime}}$ \ for all \  $\alpha\leq\lambda^\prime$
\item If \ $\lambda$ \ is a limit ordinal, \ then \ $\langle\mouse_\lambda,
\langle \pi_{\alpha\lambda}\colon \mouse_\alpha \rightarrow \mouse_\lambda
\rangle_{\alpha<\lambda}\rangle$ \ is the direct limit of \ 
\[\la {\la \mouse_\alpha \ra}_{\alpha < \lambda} , 
\la \pi_{\alpha\beta}\colon \mouse_\alpha \rightarrow \mouse_\beta
\ra_{\alpha\leq\beta<\lambda}\ra.\]
\ee
\end{definition}

The commutative system in the above (\ref{lablepmiteration}) is called the  {\it premouse iteration} of \
$\mouse$. \  We note that the maps in the above commutative system are cofinal and 
are \ $\Sigma_1$ \ embeddings, that is,
\[\pi_{\alpha\beta}\colon \mouse_\alpha \updownmap{cofinal}{1} \mouse_\beta\]
for all \ $\alpha\leq\beta\in \OR$. \ We shall call \  $\pi_{0\beta}\colon \mouse\mapsigma{1}
\mouse_\beta$ \ the {\it premouse embedding\/} of \ $\mouse$ \ into its \ $\beta^{\,\ul{\text{th}}}$ \
{\it premouse iterate\/} \ $\mouse_\beta$.
\begin{definition} A premouse  \ $\mouse$ \  is an {\it iterable
premouse} if  \ $\mouse_\lambda$ \ is well-founded for all \ $\lambda \in \OR$.
\end{definition}

For an iterable premouse  \ $\mouse$ \  and \ $\alpha\in \OR$, \ we identify \
$\mouse_\alpha$ \ with its transitive collapse. Hence, 
\[\mouse_\alpha = (M_\alpha,\en,{\underline\R}^{\mouse_\alpha},
\underline\kappa^{\mouse_\alpha},\mu,A_1,\dots,A_N)\] 
is a premouse for all ordinals \ $\alpha$. \ In this case, we write
\ $\kappa_\alpha = \pi_{0\alpha}(\underline{\kappa}^{\mouse}) = 
\underline{\kappa}^{\mouse_\alpha}$ \ for \ $\alpha\in \OR$.

\begin{rmk} Note that \ $(\mouse_\gamma)^\alpha$ \ and \ $(\mouse^\alpha)_\gamma$ \ denote 
different orders of operations, and typically \ $(\mouse_\gamma)^\alpha\ne (\mouse^\alpha)_\gamma$.
\ In this paper, the notation \ $\mouse_\gamma^\alpha$ \ is to be interpreted as
\ $(\mouse_\gamma)^\alpha$.
\end{rmk}

\subsubsection{A minimal criterion for premouse iterability}\label{mouseitercrit} Theorem \ref{thmttho} (see below)
offers a ``minimal'' relative criterion which will assure that a model of \ $\PM$ \ is an iterable premouse.
This criterion (see \cite[Theorem 2.31]{Crcm}) was used in \cite{Crcm} to produce scales \ $\Sigma_1$--definable over
an iterable real premouse.  We shall now review this criterion.
Let \ $\mouse$ \ be a model of \ $\PM$ \ and define   
\[F^{\mouse}=\{\,f\in M :  \exists n\in\omega \ \mouse \models f\colon {^n}\underline\kappa \rightarrow \OR\,\}.\] 
For \ $f\in F^{\mouse}$, \ write \ $d(f) = n$ \ if and only if \ $n\in\omega$  \ and \ $\mouse \models f\colon
{^n}\underline\kappa \rightarrow \OR$. \ We shall assume the convention that \ $f\in F^{\mouse}$ \ and  $d(f) = 0$ \
whenever \ $f\in {\OR}^{\mouse}$. \ Finally, for \ $n\in\omega$, \ define
\  $F^{\mouse}_n = \{\,f\in F^{\mouse} : d(f) = n \,\}.$

For a model \ $\nouse = (N,\en,{\underline\R}^{\nouse},
\underline\kappa^{\nouse},\mu,A_1,\dots,A_K)$ \ of  \ $\PM$, \ we now define a
predicate \ $\mu_n$ \ on \ $\{\,a\in N: \nouse \models a\subseteq 
{^n\underline\kappa}\,\}$ \ by induction on \ $n$. \ For \ $n=1$ \ let \
$\mu_1=\mu$.  Given \ $a,\xi_0 \in N$, \ where \ $\nouse\models
a\subseteq{^{n+1}\underline\kappa} \land \xi_0\in\underline\kappa$, \ let \
$a_{\xi_0}\in N$ \ be such that 
\[\nouse \models a_{\xi_0}=\{\,\langle\xi_1,\xi_2,\dots\xi_n\rangle\in {^n\underline\kappa} : 
\langle \xi_0,\xi_1,\dots\xi_n\rangle \in a \,\}.\]  
Now, assuming that \ $\mu_n$ \ is defined, 
let \ $\mu_{n+1}$ \ be defined by
\[a\in \mu_{n+1} \iff \nouse \models a\subseteq{^{n+1}\underline\kappa} \land
\{\,\xi_0 \in \underline\kappa : a_{\xi_0}\in\mu_n\,\}\in\mu\]
for all \ $a\in N$.  \ Clearly, for each
\ $n\in\omega$, \ $\mu_n$ \ is ``rudimentary over \ $\nouse$,'' \ that is, there is a
rudimentary function \ $\mathcal{F}$ \ such that 
\[a\in\mu_n \iff \nouse \models \mathcal{F}(a)\ne\emptyset\]
for all \ $a\in N$.

Let \ $R \subseteq (\OR\times \OR)^{\mouse}$ \ be any rudimentary relation. Given $f, g \in
F^{\mouse}$, \ let \ $n=d(f)$ \ and \ $m=d(g)$.  For any  \ $s\in{^n(n+m)}\!\uparrow$ \ and \ for
any \  $t\in{^m(n+m)}\!\uparrow,$ \ we shall write
\ $\mouse\models \,f\,R^{s,t}\,g$ \
if and only if
\[
\mouse\models\{\langle \xi_0,\dots,\xi_{n+m-1}\rangle \in{^{n+m}\underline\kappa}: 
f(\xi_{s(0)},\dots,\xi_{s(n-1)}) R g(\xi_{t(0)},\dots,\xi_{t(m-1)})\}\in\mu_{n+m}\]

\begin{definition}\label{defttho} Let \ $\mouse$ \ and
\ $\mathcal{A}$ \ be models of \ $\PM$. \ A map \ $\sigma\colon F^{\mathcal{A}} \rightarrow
F^{\mouse}$ \ is said to be \ $\leq$--{\it extendible} \  if, for all \ $f,g \in F^{\mathcal{A}}$, 
\be
\item $d(f)=d(\sigma(f))$
\item for all \ $s\in{^{d(f)}(d(f)+d(g))}\!\uparrow$ \ and \ for all \ 
$t\in{^{d(g)}(d(f)+d(g))}\!\uparrow$, 
\[\mathcal{A} \models f\le ^{s,t}g \iff \mouse \models \sigma(f)\le ^{s,t}\sigma(g).\] 
\ee
\end{definition}

As noted earlier, the next theorem presents a criterion that was used in \cite{Crcm} to produce scales
which are $\Sigma_1$--definable over an iterable real premouse. 

\begin{theorem}\label{thmttho} Let \ $\mouse$ \ be an iterable premouse and let \ $\mathcal{A}$ \
be a model of \ $\PM$. \ Suppose that \ $\sigma\colon F^{\mathcal{A}} \rightarrow
 F^{\mouse}$ \ is \ $\leq$--extendible.  Then \ $\mathcal{A}$ \ is
(isomorphic to) an iterable premouse. 
\end{theorem}

The following proposition will be used implicitly here, and in \cite{Part2} and \cite{Part3}.
\begin{proposition} Let \ $\mouse=(M,\R,\kappa,\mu)$ \ be an iterable real premouse. 
Then \ $\mouse^\gamma$ \ is also an iterable real premouse for
each ordinal \ $\gamma$ \ where \  $\kappa<\gamma<\widehat{\OR}^{\mouse}$.
\end{proposition}

\begin{proof} Clearly \ $\mouse^\gamma$ \ is a transitive model of \ $\PM$. \ Hence, \
$\mouse^\gamma$ \ is a real premouse. To see that \ $\mouse^\gamma$ \ is iterable, let \
$\sigma\colon F^{\mouse^\gamma}\to F^{\mouse}$ \ be the identity map. Since \ $\sigma$ \ is \
$\le$--extendible, Theorem \ref{thmttho} implies that \ $\mouse^\gamma$ \ is iterable.
\end{proof}

\subsubsection{Iterated types}
Our next result is Corollary 2.41 of \cite{Cfsrm} and states that a \ 
$\Sigma_{n+1}$--type \ $\Sigma$
\ can be ``translated'' to a \ $\Sigma_{n+1}$--type \ ${\Sigma}^*$ \ such that
a premouse iterate \ $\mouse_\alpha$ \  realizes \ $\Sigma$ \ if and only if \  $\mouse$ \ realizes \
${\Sigma}^*$, \ whenever \
$\alpha$ \ is a multiple of \ $\omega^\omega$. \ This result 
will be used in \cite{Part2} to prove Theorem \ref{newthm}. Recall the notation presented in subsection \ref{projtypes}.

\begin{lemma}\label{iteratedtypes} Suppose \ $\mouse$ \ is an
iterable premouse in the language \ $\Lng_N$ \ and \ $\alpha\in{\OR}$ \ is a multiple of \
$\omega^\omega$. \ Let \ $n\in\omega$. \ There is a map (independent of
\ $\mouse$, \ $\alpha$) \ 
$^*\, \colon \Sigma_{n+1}^N \to \Sigma_{n+1}^N$ \ mapping each \ $\Lng_N$ \ formula of
one free variable \ $\vartheta(v)\in\Sigma_{n+1}^N$ \ to \ $\vartheta^*(t,v)\in\Sigma_{n+1}^N$ \
($t\in\omega^{<\omega}$) with the following property: \  $(\forall a\in M_\alpha)(\exists
f\in M)(\exists t\in{^{d(f)}}(n+1))$ \ such that
\begin{equation}\mouse_\alpha\models \vartheta(a) \iff \mouse\models \vartheta^*(t,f), \ 
\text{ \ for all \ } \vartheta\in\Sigma_{n+1}^N.\label{typeiterate}
\end{equation} 
In addition, \ $(\forall f\in M)(\forall t\in{^{d(f)}}(n+1))(\exists a\in M_\alpha)$ \
such that {\rm (\ref{typeiterate})} holds.
\end{lemma}

\begin{corollary} Suppose \ $\mouse$ \ is an
iterable premouse in the language \ $\Lng_N$ \ and \ $\alpha\in{\OR}$ \ is a multiple of \
$\omega^\omega$. \ Let \ $n\in\omega$. \ There is a map (independent of
\ $\mouse$, \ $\alpha$) \ 
$^*\, \colon \Upsilon_{n+1}^N \to \Upsilon_{n+1}^N$ \ mapping each \ $\Lng_N$ \ formula of
one free variable \ $\vartheta(v)\in\Upsilon_{n+1}^N$ \ to \ $\vartheta^*(t,v)\in\Upsilon_{n+1}^N$ \
($t$ is a definable term) with the following property: \  $(\forall a\in M_\alpha)(\exists
f\in M)(\exists t\in{^{d(f)}}(n+1))$ \ such that
\begin{equation}\mouse_\alpha\models \vartheta(a) \iff \mouse\models \vartheta^*(t,f), \ 
\text{ \ for all \ } \vartheta\in\Upsilon_{n+1}^N.\label{typeiterate2}
\end{equation} 
In addition, \ $(\forall f\in M)(\forall t\in{^{d(f)}}(n+1))(\exists a\in M_\alpha)$ \
such that {\rm (\ref{typeiterate2})} holds.
\end{corollary}

\begin{corollary}\label{cortwotype} Suppose \ $\mouse$ \ is an iterable premouse in the language \ $\Lng_N$ \ and \
$\alpha\in{\OR}$ \ is a multiple of \ $\omega^\omega$. \ Then any  \ $\Upsilon_{n+1}^N$--type \ $\Upsilon$ \ can be translated
to an \ $\Upsilon_{n+1}^N$--type \ $\Upsilon^*$ \ such that the premouse iterate \ $\mouse_\alpha$ \ realizes \
$\Upsilon$ \ if and only if \ $\mouse$ \ realizes \ $\Upsilon^*$.
\end{corollary}

\subsection{Real 1--mice and the definition of $\Kr$}\label{realonemice}

Recall that \ $\Lr$ \ is the smallest inner model of \ $\ZF$ \ containing the reals.
An extensive theory of the structure of \ $\Lr$ \ has been
developed under the hypothesis that \ $\Lr$ \ is a model of \
$\AD$. \ Assuming determinacy for sets of reals in \ $\Lr$, \ researchers have
essentially settled all the important problems of descriptive
set theory in \ $\Lr$. \ In particular, Steel \cite{steel} determines the
complexity of scales in \ $\Lr$ \ under the hypothesis that \ $\Lr$ \ is a model of \ $\AD$. \ One concludes that this 
hypothesis is sufficient to develop the structure and descriptive set theory of \  $\Lr$. \ This success inspires one to look
for inner models of \ $\AD$ \ larger than \ $\Lr$. \ We now briefly describe how to construct one such inner model, namely \
$\Kr$.
\begin{definition} Let \ $\mouse$ \  be a transitive model of \ $\Rplus$. \
The {\it projectum} \ $\rho^{}_\mouse$ \ is the least ordinal \ $\rho \le
\widehat{\OR}^{\mouse}$ \ such that \ $\mathcal{P}(\R^{\mouse}\times
\omega\rho)\cap\boldface{\Sigma}{1}(\mouse)\not\subseteq M$, \ and \ $p^{}_\mouse$ \ is the \ 
$\le_{\BK}$--least \ $p\in[\OR^{\mouse}]^{<\omega}$ \ such that \ $\mathcal{P}(\R^{\mouse}\times
\omega\rho^{}_\mouse)\cap\Sigma_1(\mouse,\{p\})\not\subseteq M$. 
\end{definition}

\begin{definition}\label{onemouse}  An iterable pure premouse \ $\mouse$ \ is a \ {\it 1--mouse} \ if \
$\omega\rho^{}_{\mouse}\le \kappa^{\mouse}$. \ In addition, if \ $\R^{\mouse} = \R$, \ then \ $\mouse$ \ is
said to be a {\it real 1--mouse}.
\end{definition}

Using real 1--mice there a natural way to construct an inner model of \ $\AD$ \ larger than \ $\Lr$.
\begin{definition} The real core model is the class 
\ $Kr = \{\,x : \exists\, \nouse (\nouse \text{ \ is a real 1--mouse} \land x\in N)\,\}$.
\end{definition}

One can prove that \ $\Kr$ \ is an inner model of \ $\ZF$ \ and
contains a ``constructible'' set of reals not in \ $\Lr$ \ (see \cite{Crcm}). It turns out that the structure of \ $\Kr$ \
can also be developed under the hypothesis that \ $\Kr$ \ is a model of \ $\AD$. \ For example, using a mixture of descriptive
set theory, fine structure and the theory of iterated ultrapowers,  one can produce definable scales in \ $\Kr$ \ beyond those
in \ $\Lr$ \ and prove that \ $\Kr\models\DC$.

\begin{rmk} There exists an iterable real premouse if and only if \ $\R^\shrp$ \ exists.
Hence \ $\Kr$, \ as defined above,  is nonempty if and only if \ $\R^\shrp$ \ exists. Therefore, 
in the case where there are no iterable real premice, one should assume the
convention that \ $\Kr = \Lr$.
\end{rmk}

\subsection{Real mice}\label{realmice} 
Real 1--mice are, in a sense, the basic building blocks of \ $\Kr$,\ however, to gain a better understanding of the structure of \ $\Kr$, \ we need to introduce {\it real mice.}

\begin{definition}\label{master4} Let \ $\mouse$ \ be an acceptable pure premouse. We say that  \
$\mouse$ \  is {\it critical\/} if \ $\rho_\mouse^{n+1} \le \kappa^\mouse < \rho_\mouse^n$, \ for some \
$n\in\omega$. \ This integer \ $n$ \ will be denoted by \ $n(\mouse)$ \ and we shall write
\ $\overline{\mouse} = \mouse^{n(\mouse)}$ \ and \ $\overline{M} = M^{n(\mouse)}$. 
\end{definition}

In \cite{Crcm} we established that \ $\Sigma_1(\mouse)$ \  has the scale property when \ 
$\mouse = (M,\R,\kappa,\mu)$ \ is an iterable premouse satisfying the axiom of determinacy. The key
fact used to prove this theorem is that \ $\mouse$ \  is an iterate of its core (see subsection \ref{coremice} below)
when \ $\rho^{}_{\mouse} = 1$. \ This strategy fails when \ $\rho^{}_{\mouse} > 1$ \ and \ $\rho_{\mouse}^m = 1$ \ for
some \ $m>1$ \ because our iteration maps are only \ $\Sigma_1$--elementary and not necessarily
\ $\Sigma_m$--elementary. \ In this case, however, there is an \ $n$ \ such that \ $\rho_{\mouse}^{n+1} \le \kappa <
\rho_{\mouse}^n$ \ and we shall be able to define an iteration procedure which is \ $\Sigma_{n+1}$--elementary. Using
these iterations we can show that \ $\boldface{\Sigma}{m}(\mouse)$ \  has the scale property when  \ $\mouse$ \  is a
weak mouse and \ $m=m(\mouse)$.
\  Before we begin, we give an overview of this iteration procedure.
The remainder of this section is devoted to the review of such iteration, called
{\it mouse iteration.\/}  

Given an acceptable premouse  \ $\mouse$ \  with \ 
$\rho_{\mouse}^{n+1} \le \kappa^\mouse < \rho_{\mouse}^n$, \ its \ $\Sigma_n$-code \ $\mouse^n$
\ is a premouse satisfying a set of axioms \ $\mathcal{T}^n$ \ (see \cite[pp. 951--954]{Cfsrm}). First we shall take the
ultrapower of \ $\mouse^n$, \ thereby obtaining \ $\pi^n\colon \mouse^n \mapsigma{1}{\nouse^n}$. \
The embedding \ $\pi^n$ \ also preserves certain \ $\Pi_2$ \ sentences, namely those that are \
$\mouse^n$--collectible (see subsection \ref{abovereals}). It turns out that each axiom in \ $\mathcal{T}^n$ \
is  equivalent (in \ $\mouse^n$) to an \ $\mouse^n$--collectible sentence and therefore, \ $\nouse^n$ \
satisfies the axioms in \ $\mathcal{T}^n$. \  Since \ $\nouse^n$ \ is a model of \ $\mathcal{T}^n$, \ it
believes that it is the \ $\Sigma_n$--code of a structure \ $\nouse$. \ Thus we can extend \ $\pi^n$,
\ by decoding master codes, to a map \
$\pi\colon
\mouse
\mapsigma{{n+1}}{\nouse}$. \ The premouse \ $\nouse$ \ is the first mouse iterate of \ $\mouse$, \ and we
shall be able to iterate this procedure through the ordinals.

Inspired by Dodd \cite[Definitition 3.25]{Dodd}, we inductively defined in \cite{Cfsrm} the set of
axioms \ $\mathcal{T}^n$ \ (in the language \ $\Lng_n$) \ which are true in the \ $\Sigma_n$--code, \ $\mouse^n$, \
of any acceptable pure premouse \ $\mouse$ \ with \ $\rho_{\mouse}^n>\kappa^\mouse$.  \ We shall not repeat the inductive
definition of the axioms \ $\mathcal{T}^n$ \ here; however, we will review some preliminary technical notions that were used
in the definition of \ $\mathcal{T}^n$. \ Given a real \ $x$ \ and an \ $i\in\omega$, \ we shall let \ $i{^\frown}x$ \ denote
the real \ $y$ \ such that \ $y(0)=i$ \ and \ $y(m)=x(m-1)$ \ for \ $1\le m \in\omega$. \ Recall that \  $(x)_i$, \ or \ $x_i$
\ when the context is clear, denotes the \ $i^{\ul{\text{th}}}$ \ real coded by \ $x$. \ Given a finite
sequence \
$x_0,\dots,x_i$ \ of reals we shall write \ $\langle x_0,\dots,x_i \rangle$ \ to denote an effective coding of these reals by
a real \ $y$ \ where \ $y=\langle x_0,\dots,x_i \rangle$ \ such that \ $(y)_j= x_j$ \ for all \ $j\le i$. \ Given a sequence \
$s$ \ of ordinals, we shall let \ $s_i$  \ denote the \ $i^{\ul{\text{th}}}$ \ element of the sequence \
$s$. \ Let \
$\Lng_n^{\underline{p}} = \Lng_n\cup\{\underline{p}\}$, \ where \ $\underline{p}$ \ is a new constant symbol. We
are assuming an effective G\"odel numbering of all \ $\Sigma_1$ \ formulae  \ $\varphi$ \ in \ $\Lng_n^{\underline{p}}$
\ where the natural number  \ $\ulcorner\varphi\urcorner$ \ denotes the G\"odel number of \ $\varphi$. \ Finally, given any
\ $\Sigma_1$ \ formula \ $\varphi(v_{j_0},v_{j_1},\dots,v_{j_m})$ \ of \ $\Lng_n^{\underline{p}}$, \ with free
variables as displayed, \ let \ $\varphi(v_{j_0},v_{j_1},\dots,v_{j_m})^+$ \ denote the \ $\Sigma_1$ \ formula \ $\psi(x,s)$ \
given by
\[\begin{aligned}
 x\in\underline{\R}\,&\land\, s\in{^{j_m+1}{\OR}}\\
&\land\,
\exists v_{j_0}\exists v_{j_1}\cdots\exists
v_{j_m}\left[\varphi(v_{j_0},v_{j_1},\dots,v_{j_m})\land \bigwedge_{i=0}^{m}h(x_{j_i},\langle
s_{j_i},\underline{p}\rangle)=v_{j_i}\right]\\
\end{aligned}\]
where \ $h$ \ is the canonical \ $\Sigma_1$ \ Skolem function (see Definition \ref{skolemdef}).
We assume some recursive map \
$\ulcorner\varphi\urcorner\mapsto\ulcorner\varphi^+\urcorner$.
\ We shall use the abbreviation \ $D(a,x,s)$ \ for the \ $\Lng_{n+1}$-formula  
\begin{equation}\label{domain}
\begin{aligned}
a\in [\omega]^{<\omega}&\land s\in(\OR)^{<\omega} \land x\in \ul{\R}\, \land\\
&\forall i\in a [i\in\text{dom}(s) \land  A_{n+1}(\,\ulcorner(\exists v(v=v_i))^+\urcorner{^\frown}
x,s\,)].
\end{aligned}
\end{equation}

\begin{remark}\label{code} The expression \ $D(a,x,s)$ \ asserts that for each \ $i\in a$, \ ``$(x_i,s_i)$ \ is a code,''
that is, \ ``$h(x_i,\langle s_i,\ul{p}\rangle)$ \ exists.''
\end{remark}

In order to characterize those models of \ $\mathcal{T}^n$ \ with well-founded
extensions, we also defined an \ $\Lng_n$--formula \ $u\,\ul{E}_n\, u^\prime$
\ inductively on \ $n$ \ (see \cite[Definition 2.10]{Cfsrm}). 
For an acceptable pure premouse \ $\mouse$ \ with \ $\rho^n_\mouse >
\kappa^\mouse$, \ one can check that \ $\mouse^n\models\mathcal{T}^n$ \ and that \ $E_n^{\mouse^n}$ \ is
well-founded (see the proof of \cite[Corollary 2.13]{Cfsrm}). The next theorem provides a converse  (see
\cite[Theorem 2.11]{Cfsrm}).

\begin{theorem}[Model Extension Theorem]\label{MET} Suppose that \ $\mathcal{A}$ \ is an \ $\Lng_n$--model of \
$\mathcal{T}^n$
\ with \
$E_n^{\mathcal{A}}$ \ well-founded. Then there is a premouse \ $\mouse$ \ isomorphic to \ $\mathcal{A}$ \ and
if \ $n>0$, \ then there is also a premouse \ $\nouse$ \ which is an $\Lng^{\ul{p}}_{n-1}$--model such
that
\be
\item $\ul{p}^\nouse\in [\OR]^{<\omega}$ \ and
\ 
$h_\nouse^{\ul{p}^\nouse}\colon \R^\nouse\times \omega\rho^{}_\nouse \maps{onto} N$.
\item $\omega\rho^{}_\nouse = {\OR}^\mouse$.
\item $\R^\mouse = \R^\nouse$ \ and \ $\kappa^\mouse = \kappa^\nouse$.
\item $A_k^\mouse\cap M = A_k^\nouse\cap M$ \ for \ $1\le k \le n$.
\item $M = H_{\omega\rho^{}_\nouse}^\nouse$ \ and \ $\mu^\mouse\cap M = \mu^\nouse\cap N$.
\item $\nouse\models \mathcal{T}^{n-1}$ \ and \ $E_{n-1}^\nouse$ \ is well-founded.
\item  For every \ $\Sigma_1$ \ formula \ $\phi(v_0,\dots,v_m)$ \ in the language \ $\Lng^{\ul{p}}_{n-1}$ 
\[\nouse\models \phi(h(x_0,\langle s_0,\ul{p}^\nouse\rangle),\dots,
h(x_m,\langle s_m,\ul{p}^\nouse\rangle))\]
if and only if
\[\mouse\models A_n(\ulcorner\phi(v_0,\dots,v_m)^+\urcorner{^\frown}x,s)
\]
for every \ $x\in\R^\nouse$ \ and \ $s\in(\omega\rho^{}_\nouse)^{<\omega}$ \ where
\ $m+1\subseteq\text{dom}(s)$.
\item $p^{}_\nouse = \ul{p}^\nouse$.
\item $\nouse$ \ is sound.
\item $A_n^\mouse$ \ is the \ $\Sigma_1$ \ master code of \ $\nouse$. \ Hence, \ $\mouse$ \ is the 
\ $\Sigma_1$--code of \ $\nouse$.
\ee
\end{theorem}

\begin{theorem}\label{kappasound} Suppose \ $\mouse$ \ is an acceptable pure premouse and let \
$n\in\omega$. \ If \
$\omega\rho_\mouse^{n+1}\ge\kappa^\mouse$, \ then \ $\mouse$ \ is \ $(n+1)$--sound.
\end{theorem}

\begin{corollary}\label{kappagamma} Suppose \ $\mouse$ \ is an acceptable pure premouse and let \
$n\in\omega$. \ If \
$\omega\rho_\mouse^{n+1}\ge\kappa^\mouse$, \ then
$\gamma_\mouse^i=\rho_\mouse^i$ \  for all \ $i\le n+1$.
\end{corollary}

\begin{theorem}\label{skolemthm}
Suppose \ $\mouse$ \ is an acceptable pure premouse and let \
$k\in\omega$. \ If \
$\omega\rho_\mouse^{k}\ge\kappa^\mouse$, \ then there is a \ $\Sigma_{k+1}$ \ Skolem
function for \ $\mouse$.
\end{theorem}

\begin{theorem}\label{SCL}
Let \ $\mouse$ \ be an acceptable pure premouse and let \
$k\in\omega$. \ Suppose that  \ $\omega\rho_\mouse^{k}\ge\kappa^\mouse$ \ and let \ $\lambda=\textup{max}\{\kappa^\mouse,
\omega\rho_\mouse^{k+1}\}$. \ Then there is a \ $\boldface{\Sigma}{k+1}(\mouse)$ \ function \
$F\colon \lambda\times\R \maps{onto} M$.
\end{theorem}

\begin{proof} By Theorem \ref{kappasound} and Corollary \ref{kappagamma},
\ $\omega\rho_\mouse^{k+1}=\omega\gamma_\mouse^{k+1}$. \ So let \ $p\in M$ \ be such that there is a \
$\Sigma_{k+1}(\mouse,\{p\})$ \ set \ $D\subseteq\omega\rho_\mouse^{k+1}$ \ such that \
$D\notin M$ \ and there is a \ $\Sigma_{k+1}$ \ Skolem function \ $F_p$ \ for \ $\mouse$ \ which is
\ $\Sigma_{k+1}(\mouse,\{p\})$. \ Let \ $\mathcal{H} =\Hull_{k+1}^{\mouse}({\R\cup\lambda} \cup \{p\})$ \ and let \ $H$ \
denote the domain of this structure. Because \ $F_p$ \ is a \ $\Sigma_{k+1}$ \ Skolem function for \ $\mouse$, \ it follows
that \ $\mathcal{H}\prec_{k+1}\mouse$ \ and  \ $F_p\colon \lambda\times\R \maps{onto} H$. \ Let \ $\nouse$ \ be the transitive
collapse of \ $\Hull_{k+1}^{\mouse}({\R\cup \lambda} \cup \{p\})$ \ and let \ $\tau\colon \nouse \mapsigma{k+1} \mouse$ \ be
the inverse of the collapse map. Since \ $\kappa^{\nouse} = \kappa^{\mouse}$, \ it follows (by a condensation argument) that \
$\nouse$ \ is an initial segment of \ $\mouse$ \ and \ $\nouse = \Hull_{k+1}^{\nouse}({\R\cup \lambda} \cup \{\tau^{-1}(p)\})$.
\ Let \ $F_{\tau^{-1}(p)}$ \ be the natural interpretation of \ $F_p$ \ in \ $\nouse$. \ It follows that  \
$F_{\tau^{-1}(p)}$ \ is \ $\Sigma_{k+1}(\nouse,\{\tau^{-1}(p)\})$ \ and \ $F_{\tau^{-1}(p)}\colon \lambda\times\R \maps{onto}
N$. \ Now, since \  $D\in \pow(\R\times \omega\rho^{}_{\mouse})\cap\Sigma_{k+1}(\mouse,\{p\})$, \ it
follows that \ $D\in \pow(\R\times \omega\rho_\mouse^{k+1})\cap\Sigma_{k+1}(\nouse,\{\tau^{-1}(p)\})$.
Therefore \ $\nouse = \mouse$, \ otherwise \ $D\in M$. \ Thus, \
$F_{\tau^{-1}(p)}\colon \lambda\times\R \maps{onto} M$ \ and is \ $\Sigma_{k+1}(\mouse,\{\tau^{-1}(p)\})$.
\end{proof}

\subsubsection{Mouse iteration} The following extension of embeddings lemma (see \cite[Lemma 2.14]{Cfsrm}) is the key
result which allows us to define mouse iteration. Since a proof of this fundamental result was
not presented in \cite{Cfsrm}, we shall now provide a proof of this important technical lemma. 
\begin{lemma}[Extension of Embeddings Lemma]\label{EOEL} Let \ $\mathcal{A}$ \ and \ $\mathcal{B}$ \ be
premice in the language \ $\Lng_n$. \ Suppose that \ $\mathcal{A}$ \ and  \ $\mathcal{B}$ \ are models of \
$\mathcal{T}^n$ \ where \ $E_n^{\mathcal{A}}$ \ and \ $E_n^{\mathcal{B}}$ \ are well-founded. If \ $\pi^n\colon
\mathcal{A}\mapsigma{k} \mathcal{B}$ \ for \ $k\ge 1$, \ then there are acceptable pure premice \ $\mouse$ \ and \ $\nouse$ \
together with a map \ $\pi\colon \mouse \mapsigma{{n+k}} \nouse$ \  such that
\be
\item[(1)] $\mathcal{A}=\mouse^n$ \ and \ $\mathcal{B}=\nouse^n$
\item[(2)] $\pi^n \subseteq \pi$ \ and for \ $i\le n$ \ $\pi(p_\mouse^i) = p_\nouse^i$.
\ee
\end{lemma}

\begin{proof} By Theorem \ref{MET} there are pure premice \ $\mouse$ \ and \ $\nouse$ \ such that \
$\mouse^n=\A$ \ and \ $\nouse^n=\B$. \ So, \ $\pi^n\colon
\mouse^n\mapsigma{k} \nouse^n$. \ We show how to construct a map \ $\pi^{n-1}\colon
\mouse^{n-1}\mapsigma{k+1} \nouse^{n-1}$ \ such that \ $\pi^n \subseteq \pi^{n-1}$ \ and \
$\pi^{n-1}(p^{}_{\mouse^{n-1}}) = p^{}_{\nouse^{n-1}}$. \ If we then iterate this construction, then \
$\pi^0\colon\mouse\mapsigma{k+n}\nouse$ \ will be our desired map \ $\pi$. 

For notational convenience let \ $\fmouse=\mouse^{n-1}$ \ and \ $\fnouse=\nouse^{n-1}$. \ 
Define \ $\pi^{n-1}\colon
\fmouse\to\fnouse$ \ by 
\begin{equation}\pi^{n-1}\left(h^{}_\fmouse\,(x,\langle \gamma,p^{}_{\fmouse}\,\rangle)\right)= 
h^{}_\fnouse\,(x,\langle\pi^n(\gamma),p^{}_{\fnouse}\,\rangle)\label{pimap}\end{equation}
where \ $x\in\R^\fmouse$ \ and \ $\gamma\in\OR^\fmouse$. \ 
Since \ $\pi^n\colon \mouse^n\mapsigma{k} \nouse^n$, \ Theorem \ref{MET} and its proof (see \cite[Theorem 2.11]{Cfsrm})
imply that 
\be
\item[(a)] $\pi^{n-1}$ \ is well-defined
\item[(b)] $\pi^{n-1}\colon \fmouse\mapsigma{1}\fnouse$
\item[(c)] $\pi^{n-1}(x)=\pi^n(x)=x$ \ for all \ $x\in\R^\fmouse$
\item[(d)] $\pi^{n-1}(\gamma)=\pi^{n}(\gamma)$ \ for all \ $\gamma\in\omega\rho^{}_\fmouse$
\item[(d)] $\pi^{n-1}(p^{}_{\fmouse})=p^{}_{\fnouse}\,$.
\ee
\renewcommand{\theclam}{\arabic{clam}}
\begin{clam} $\pi^{n-1}\colon \fmouse\mapsigma{k+1}\fnouse$.
\end{clam}

\begin{proof}[Proof of Claim 1] In the interest of simplifying notation we will assume that \ $k=2$ \ without any loss
of generality. Thus, we are assuming that \ $\pi^n\colon \mouse^n\mapsigma{2} \nouse^n$ \ and we want to prove that \
$\pi^{n-1}\colon\fmouse\mapsigma{3}\fnouse$.  \ Let \ $\psi(v_0,\dots,v_i,v_{i+1},v_{i+2})$ \ be a \ $\Sigma_1$ \
formula in the language \ $\Lng_{n-1}$. \ Let \ $\varphi(v_0,\dots,v_i)$ \ be the \ $\Sigma_3$ \ formula \ $(\exists
v_{i+1})(\forall v_{i+2})\psi(v_0,\dots,v_i,v_{i+1},v_{i+2})$. \ Let \ $m_0,\dots,m_i$ \ be arbitrary
elements in \ $\fmouse$, \ that is, in the domain of \ $\fmouse$. \ We show that 
\[\fmouse\models\varphi[m_0,\dots,m_i]\iff\fnouse\models\varphi[\pi^{n-1}(m_0),\dots,\pi^{n-1}(m_i)].\]
By (1) of Lemma \ref{MET}, there exist \ $x_0,\dots,x_i\in\R^\fmouse$ \ and \ 
$\gamma_0,\dots,\gamma_i\in\omega\rho^{}_\fmouse$ \ such that \ $m_j=h^{}_\fmouse\,(x_j,\langle
\gamma_j,p^{}_{\fmouse}\,\rangle)$ \ for \ $0\le j\le i$. \  Let \ $x=\langle x_0,\dots,x_i\rangle\in\R^\fmouse$ \ and \
$s=\langle\gamma_0,\dots,\gamma_i\rangle\in (\omega\rho^{}_\fmouse)^{<\omega}$. \ Recalling (\ref{domain}) and Remark
\ref{code},  let \ $\chi(x,s)$ \ be the formula 
\[\begin{aligned}
(\exists x_{i+1})&(\exists\gamma_{i+1})(\forall x_{i+2})(\forall\gamma_{i+2}) 
(D(i+1,x,s)\land D(1,\langle x_{i+1}\rangle,\langle\gamma_{i+1}\rangle)\\
&\land[ D(1,\langle x_{i+2}\rangle,\langle\gamma_{i+2}\rangle)\rightarrow
A_n(\ulcorner\varphi^+\urcorner{^\frown}(x{^\frown}\langle x_{i+1},x_{i+2}\rangle),s{^\frown}\langle
\gamma_{i+1},\gamma_{i+2}\rangle) ])\\
\end{aligned}\]
where \ $x{^\frown}\langle x_{i+1},x_{i+2}\rangle=\langle x_0,\dots,x_i,x_{i+1},x_{i+2}\rangle\in\R^\fmouse$ \ and \
$\langle y\rangle\in\R^\fmouse$ \ is such that \ $\langle y\rangle_0=y$. \ The
above formula \ $\chi(x,s)$ \ is \ $\Sigma_2$ \ in the language \ $\Lng_n$.\footnote{$\Sigma_2$ in the theory
$\Rplus_n$.}
\  Note that \ $\pi^n(x)=x$ \ and \
$\pi^n(s)=\langle\pi^n(\gamma_0),\dots,\pi^n(\gamma_{i})\rangle=\langle\gamma_0',\dots,\gamma_i'\rangle
\in(\omega\rho^{}_\fnouse)^{<\omega}$.
 \ The following holds
\begin{alignat*}{2}
&\fmouse\models\varphi[m_0,\dots,m_i]&&\\ 
&\iff \fmouse\models\varphi[h^{}_\fmouse\,(x_0,\langle
\gamma_0,p^{}_{\fmouse}\,\rangle),\dots,h^{}_\fmouse\,(x_i,\langle \gamma_i,p^{}_{\fmouse}\,\rangle)] &&\quad\text{by Lemma
\ref{MET}(1)}\\ 
&\iff \mouse^n\models\chi(x,s)&&\quad\text{by Lemma \ref{MET}(1,2)}\\
&\iff \nouse^n\models\chi(x,\pi^n(s))&&\quad\text{by $\pi^n\colon \mouse^n\mapsigma{2}\nouse^n$ \& (c)}\\
&\iff \fnouse\models\varphi[h^{}_\fnouse\,(x_0,\langle
\gamma_0',p^{}_{\fnouse}\,\rangle),\dots,h^{}_\fnouse\,(x_i,\langle \gamma_i',p^{}_{\fnouse}\,\rangle)] &&\quad\text{by Lemma
\ref{MET}(1,2)}\\
&\iff \fnouse\models\varphi[\pi^{n-1}(m_0),\dots,\pi^{n-1}(m_i)]&&\quad\text{by equation (\ref{pimap}).}
\end{alignat*}
This completes the proof of Claim 1.\end{proof}

\begin{clam} $\pi^{n}(B)=\pi^{n-1}(B)$, \ for all \ $B\in M^n$ \ such that \ $B\subseteq\OR^{\mouse^n}\times\R^{\mouse^n}$.
\end{clam}

\begin{proof} Let \ $h$ \ be the canonical \ $\Sigma_1$--Skolem
function (see Definition \ref{skolemdef}) and let \ $i\in\omega$ \ be such that \ 
`$h(i{^\frown}x,\langle a,b\rangle)=\langle x,a\rangle$' \ is true in any transitive model of \ $\Rplus_n$. \ Now,
since
\ $B\in M^n$, \ it follows that \ $B\in \fmouse$ \ (i.e., the domain of \ $\fmouse$). \ Let \ $x\in\R^\fmouse$ \ and
\ $\gamma\in\omega\rho^{}_\fmouse$ \ be such that \ $h^{}_\fmouse\,(x,\langle
\gamma,p^{}_{\fmouse}\,\rangle)=B$. \ Let \ $\psi(v,x,\gamma)$ \ be the \ $\Pi_1$ \ formula
\[(\forall \xi)(\forall y) [\, \langle\xi,y\rangle\in v \leftrightarrow A_n( \ulcorner(v_0\in
v_1)^+\urcorner{^\frown}\langle i{^\frown}y,x\rangle, \langle\xi,\gamma\rangle)\,].\]
The following holds
\begin{alignat*}{2}
&\fmouse\models B= h(x,\langle\gamma,p^{}_{\fmouse}\,\rangle)&&\\ 
&\iff \mouse^n\models\psi(B,x,\gamma) &&\quad\text{by Lemma \ref{MET}}\\ 
&\iff \nouse^n\models\psi(\pi^n(B),x,\pi^n(\gamma))&&\quad\text{$\pi^n\colon \mouse^n\mapsigma{k}\nouse^n$ \& (c)}\\
&\iff \fnouse\models \pi^n(B) = h(x,\langle\pi^n(\gamma),p^{}_{\fnouse}\,\rangle) &&\quad\text{by Lemma \ref{MET}}\\ 
&\iff \fnouse\models \pi^n(B) = h(x,\langle\pi^{n-1}(\gamma),p^{}_{\fnouse}\,\rangle) &&\quad\text{by (d)}.
\end{alignat*}
Since \ $\pi^{n-1}(B)= h^{}_\fnouse(x,\langle\pi^{n-1}(\gamma),p^{}_{\fnouse}\,\rangle)$ \ by the definition of \
$\pi^{n-1}$, \ it now follows that \ $\pi^{n}(B)=\pi^{n-1}(B)$. \ This completes the proof of Claim 2.
\end{proof}

\begin{clam} $\pi^{n}\subseteq\pi^{n-1}$.
\end{clam}

\begin{proof} Let \ $a\in M^n$ \ be arbitrary. We will show that \ $\pi^{n-1}(a)=\pi^{n}(a)$. \ Since \ $a\in M^n$,
\ there is a function \ $g\colon\lambda\times\R^{\mouse^n}\maps{onto}T_c(\{a\})$ \ in \ $\mouse^n$ \ where \
$\lambda<\OR^{\mouse^n}$. \ Let \ $B$ \ be the relation on \ $(\lambda\times\R^{\mouse^n})^2$ \
defined by
\ $\langle\xi,x \rangle B \langle\xi',x'\rangle$ \iff $g(\langle\xi,x \rangle) \in g(\langle\xi',x' \rangle)$. 
Because \ $B\in \mouse^n$ \ can be coded in \ $\mouse^n$ \ by a subset of \ $\OR^{\mouse^n}\times\R^{\mouse^n}$, \ we
shall assume that \ $B\subseteq\OR^{\mouse^n}\times\R^{\mouse^n}$. \ Since \ $\fmouse$ \ is acceptable above the
reals and because \ ${\OR}^{\mouse^n}$ \ is a \ $\boldface{\Sigma}{1}(\fmouse)$ \ cardinal, it follows that 
\ $J_{\lambda+\gamma}(\R^\mouse)[B]\in M^n$ \ where \ $\gamma=\text{rank}(\{a\})<{\OR}^{\mouse^n}$ \ (see \cite[Lemma
1.15]{Cfsrm} and \cite[Lemma 3.18]{Dodd}). \ Thus, the transitive collapse of \ $B$ \ is in \ $\mouse^n$ \ and hence,
also in \ $\fmouse$. \ Therefore, there is a uniformly \ $\Sigma_1$ \ function \ $f$ \ which collapses the relation
\ $B$ \ and is absolute between \ $\fmouse$ \ and \ $\mouse^n$. \ Note that there exists \ $\eta\in{\OR}^{\mouse^n}$ \
 and \ $y\in\R^{\mouse^n}$ \ such that \ $\mouse^n\models a=f(B,\eta,y)$. \ Since this collapse is absolute, we also
have that \ $\fmouse\models a=f(B,\eta,y)$. \ Because the maps \ $\pi^n$ \ and \ $\pi^{n-1}$ \ are at least \
$\Sigma_1$ \ elementary, we have by (c) above that
\bi
\item $\nouse^n\models \pi^n(a)=f(\pi^n(B),\pi^n(\eta),y)$
\item $\fnouse\models \pi^{n-1}(a)=f(\pi^{n-1}(B),\pi^{n-1}(\eta),y)$.
\ei
However, \ $\pi^n(B)=\pi^{n-1}(B)$ \ by Claim 2 and \ $\pi^n(\eta)=\pi^{n-1}(\eta)$ \ by (d) above. Therefore,
\ $\pi^n(a)=\pi^{n-1}(a)$ \ by absoluteness. The proof of Claim 3 is complete.
\end{proof}

The proof of the extension of embeddings lemma is now complete.\end{proof}

We now review the formal definition of mouse iteration. Let \ $\mouse$ \ be an acceptable pure premouse
and let \ $n\in\omega$ \ be such that \ $\rho_\mouse^n > \kappa^\mouse$. \ Let \
$\nouse=\mouse^n$. \ Since \ $\nouse$ \ is a premouse, let
\begin{equation}\nousensystem\tag{$*$}\end{equation}
be the {\it premouse iteration} of \ $\nouse$.
Since \ $\nouse$ \ is a transitive model of \ $\mathcal{T}^n$ \ and
the theory \ $\mathcal{T}^n$ \ is preserved by cofinal \ $\Sigma_1$ \ embeddings of \ $\nouse$, \ 
it follows that \
$\nouse_\alpha\models \mathcal{T}^n$ \ for all ordinals \ $\alpha$. \ If \ $E_n^{\nouse_\alpha}$ \ is
well-founded for each \ $\alpha$, \ then the extension of embeddings Lemma \ref{EOEL} yields a
commutative system
\begin{equation}\mousesystem\tag{$\#$}\end{equation}
such that, for all ordinals \ $\alpha\le\beta$
\be
\item $\nouse_\alpha = (\mouse_\alpha)^n$
\item $\pi_{\alpha\beta}(p^i_{\mouse_\alpha}) = p^i_{\mouse_\beta}$ \ for all \ $i\le n$.
\ee

\begin{definition}\label{niteratable} Let  \ $\mouse$ \  be an acceptable pure premouse. Given \ $n\in\omega$ \
such that \ $\rho_\mouse^n > \kappa^\mouse$, \ let \ $\nouse=\mouse^n$. \ Suppose that
\ $E_n^{\nouse_\alpha}$ \ is well-founded for all ordinals \ $\alpha$, \ where \ $\nouse_\alpha$ \ is
as defined in the above commutative system ($*$). \ Then we say that
 \ $\mouse$ \  is \ {\it $n$--iterable} \ and we call the above commutative system \ ($\#$) \ the \ {\it
$n$--iteration} \ of \ $\mouse$.
\end{definition} 

\begin{rmk} Note that \ $(\mouse_\alpha)^n$ \ and \ $(\mouse^n)_\alpha$ \ denote 
different orders of operations, and typically \ $(\mouse_\alpha)^n\ne (\mouse^n)_\alpha$.
\ In this paper, when we use the notation \ $\mouse^n_\alpha$ \ our intended order of operations shall be made clear
either explicitly or from the context.
\end{rmk}

\begin{definition}\label{mouseiteratable} Suppose that  \ $\mouse$ \  is a critical pure premouse which
is \ $n(\mouse)$--iterable. Then
 \ $\mouse$ \  is called a \ {\it mouse} \ and the \ $n(\mouse)$--iteration of
 \ $\mouse$ \  is called the {\it mouse iteration} of \ $\mouse$. In addition, if  \ $\mouse$ \ 
contains all the reals, that is, if \ $\R^\mouse = \R$, \ then  \ $\mouse$ \  is said to be a {\it real mouse.}
\end{definition} 

For a mouse  \ $\mouse$ \  with mouse iteration \[\mousesystem\] where \ $n=n(\mouse)$, \ we have that each one of the
 mouse iterates \ $\mouse_\alpha$ \ is critical and \ $n(\mouse_\alpha) = n(\mouse)$, \ by
applying Corollary 2.14(2) of \cite{Crcm} to the premouse iteration of \ $\overline{\mouse}=\mouse^n$.
\ We note that for a mouse  \ $\mouse$ \  and \ $\alpha\in \OR$, \ we always identify \
$\overline{\mouse}_\alpha$ \ and \  $\mouse_\alpha$ \ with their respective transitive collapses.

{\sc To summarize:} A {\it real mouse\/} \ $\mouse$ \  contains all the reals and has the form \ 
$\mouse=(M,\R,\kappa,\mu)$. \ In addition, \ $\mouse$ \ is acceptable and critical. Let \ $n=n(\mouse)$ \ be the unique
integer such that  \ $\rho_{\mouse}^{n+1}\le \kappa^\mouse < \rho_{\mouse}^n$.   \ Now  let \ 
$\overline{\mouse}=\mouse^n$. \ Since \ $\overline{\mouse}$ \  is an iterable real premouse, let
\begin{equation}\oversystem\end{equation}
be the  {\it premouse iteration\/} of \ $\overline{\mouse}$ \ as in Definition \ref{pmiteration}.
We can extend this system of transitive models via the extension of embeddings lemma and
obtain the commutative system of transitive structures
\begin{equation}\mousesystem.\label{mouse.iter}\end{equation}
The system (\ref{mouse.iter}) is called the {\it mouse iteration\/} of \ $\mouse$. \ We shall call  \
$\pi_{0\beta}\colon \mouse\mapsigma{n+1} \mouse_\beta$ \ the {\it mouse embedding\/} of \ $\mouse$ \ into its \
$\beta^{\,\ul{\text{th}}}$ \ {\it mouse iterate\/} \ $\mouse_\beta$.

\begin{rmk} A real 1--mouse \ $\mouse=(M,\R,\kappa,\mu)$ \ is the simplest of real mice; because \ $\mouse$ \ is
iterable and \ $\mathcal{P}(\R\times\kappa)\cap\boldface{\Sigma}{1}(M)\not\subseteq M$.
\end{rmk}

\begin{theorem}\label{gamma} If \ $\mouse$ \  is a real mouse, then \ 
$\gamma_\mouse^n=\rho_\mouse^n$ \ whenever \ $\rho_\mouse^n$ \ is defined.    
\end{theorem}

\subsubsection{Core mice}\label{coremice}
Let \ $\mouse$ \ be a mouse, let \ 
$\mathcal{H} = \Hull_1^{\overline{\mouse}}({\R^{\overline{\mouse}}\cup \omega\rho^{}_{\,\overline{\mouse}}}
\cup \{p^{}_{\,\overline{\mouse}}\})$, \ and let  
\ $\mathcal{A}$ \ be the transitive collapse of \ $\mathcal{H}$. \ Let \ 
$\overline{\sigma}\colon \mathcal{A} \mapsigma{1} \overline{\mouse}$ \ be the inverse of the collapse map. It follows
that \ $\mathcal{A}$ \  is a transitive model of \ $\mathcal{T}^n$, \ and \ $\R^{\mathcal{A}}=\R^\mouse$.
\ By the extension of embeddings lemma there is  an acceptable pure premouse \ $\core$ \
and a map \ $\sigma \supseteq \overline{\sigma}$ \ such that 
\be
\item $\core^n = \mathcal{A}$, \ where \ $n=n(\mouse)$
\item $\sigma\colon \core \mapsigma{{n+1}}\mouse$.
\ee
 We denote this acceptable pure premouse \ $\core$ \ by \ $\core(\mouse)$,
\ and denote \ $\sigma$, $\overline{\sigma}$ \ by \ $\sigma_{\mouse}$,
$\overline{\sigma}_{\mouse}$, \ respectively. We call \ $\core$ \ {\it the core of $\mouse$.\/} We note that \
$\core$ \ is also a mouse. 
\begin{lemma} Suppose that \ $\mouse$ \ is a mouse and let \
$\core=\core(\mouse)$. \ Then \ $\core$ \ is a mouse with \ $n(\core) = n(\mouse) = n$ \ and \
$\rho_{\core}^{n+1}\le\rho_{\mouse}^{n+1}$.
\end{lemma}

\begin{lemma}\label{relation1} Let \ $\mouse$ \ be a mouse and let \ $\pi_{0\alpha}\colon
\mouse \mapsigma{{n+1}} \mouse_\alpha$, \ where \ $n=n(\mouse)$, \ be the mouse embedding of \
$\mouse$
\ into its \ $\alpha^{\ul{\text{th}}}$ \ mouse iterate \ $\mouse_\alpha$. \ Then
\be
\item $\mouse_\alpha$ \ is a mouse \ and \ $\rho_{\mouse_\alpha}^{n+1} = \rho_{\mouse}^{n+1}$
\item $p_{\mouse_\alpha}^{n+1} = \pi_{0\alpha}(p_\mouse^{n+1})$
\item $\core(\mouse_\alpha) = \core(\mouse)$.
\ee
\end{lemma}

For a proof of the following theorem see Theorem 2.33 of \cite{Cfsrm}.
\begin{theorem}\label{relationtwo} Let \ $\mouse$ \ 
be a mouse with core \ $\core=\core(\mouse)$ \ and let \ $n=n(\mouse)$. \ Then the following
hold: 
\be
\item There is a premouse iterate \ $\overline{\core}_\theta$ \ for some ordinal \ $\theta$, \ such that \
$\overline{\core}_\theta = \overline{\mouse}$; \ and so, \ $\core_\theta = \mouse$
\item $\omega\rho_{\core}^{n+1} = \omega\rho_{\mouse}^{n+1}$
\item $p_{\core}^{n+1}=\overline{\sigma}^{\,-1}_{\mouse}(p_{\mouse}^{n+1})$.
\ee
\end{theorem}

Our next lemma implies the existence of definable Skolem functions for
core mice (see \cite[Lemma 2.34, Corollaries 2.35, 2.36]{Cfsrm}).

\begin{lemma}\label{soundnplusone} Let \ $\core$ \ be a core mouse and let \ $n\in\omega$.
If \ $\rho_\core^n > 1$, \ then \ $\core$ \ is \ $(n+1)$--sound.
\end{lemma}

\begin{corollary}\label{skolem} Let \ $\core$ \ be a core mouse and let \ $k\in\omega$.
If \ $\rho_{\overcore}^k > 1$, \ then \ $\overcore$ \ satisfies \ $\Sigma_{k+1}$ \
selection and thus, there is a \ $\Sigma_{k+1}$ \ Skolem function for \ $\overcore$.
\end{corollary}

\subsubsection{Indiscernibles}
We shall now review how Dodd's analysis of indiscernibles in \cite{Dodd} generalizes to ``premice above
the reals.'' 
\begin{definition} Let \ $\mouse$ \ be a premouse and let \ $X\subseteq
M$. \ A set \ $I \subseteq  \OR^{\mouse}$ \ is a set of order \ $\Sigma_n(\mouse,X)$
\ {\it indiscernibles\/ } if, for any \ $\Sigma_m$ \ formula \
$\varphi(v_0,v_1,\dots,v_{k-1})$ \ with parameters allowed from \ $X$,
\[\mouse \models \varphi(\alpha_0,\alpha_1,\dots,\alpha_{k-1}) \iff 
\mouse \models \varphi(\beta_0,\beta_1,\dots,\beta_{k-1}),\]
for all \ $\alpha_0<\alpha_1<\dots<\alpha_{k-1} \text{ \ and all \ } 
\beta_0<\beta_1<\dots<\beta_{k-1}$ \ taken from \ $I$.
\end{definition}

\begin{definition} 
Every ordinal is said to be \ $0$--{\it good}. \ Suppose that \
$m\in\omega$ \ and that the notion of \ $m$--{\it good} \ has been defined. An ordinal \
$\alpha$ \ is said to be \ $(m+1)$--{\it good} \ if \ $\alpha$ \ is a limit of \ $m$--{\it good} \
ordinals.
\end{definition}

\begin{rmk} If \ $\alpha$ \ is  \ $(m+1)$--good, \ then \ $\alpha$ \ is  \
$m$--good. \ Note that \ $\alpha$ \ is  \ $m$--good \ if and only if \ $\alpha$ \ is a multiple
of \ $\omega^m$.
\end{rmk}

The key notion that allows us to obtain indiscernibles is that of a full sequence of
indiscernibles. For ordinals \ $\beta < \alpha$, \  we shall say that \ $\alpha$ \ is \ $m$--{\it better}
\ than \ $\beta$ \ if, \[\beta \ \text{\ is \ $m$--good } 
\implies \ \alpha \ \text{\ is \ $(m+1)$--good.}\]
\begin{definition} Every increasing sequence \ $\langle \alpha_0, \dots, \alpha_k
\rangle$ \ with \ $\alpha_0 = 0$ \ is said to be \ $0$--{\it full}. \ Suppose that \
$m\in\omega$ \ and that the notion of \ $m$--{\it full} \ has been defined. A sequence \ $\langle
\alpha_1, \dots, \alpha_k \rangle$ \  is said to be \ $(m+1)$--{\it full}, \ if 
\be
\item $\langle \alpha_1, \dots, \alpha_k \rangle$ \ is  \ $m$--{\it full}
\item $(\forall j<k)(\forall \beta)(\alpha_j < \beta < \alpha_{j+1} \Rightarrow \alpha_{j+1} \ \text{is  $m$--better
than} \ \beta)$.
\ee
\end{definition}

\begin{definition} $\text{ch}_m(\alpha) = \max\{\,i\le m : \alpha \text{ is \ $i$--good}\,\}.$ 
\end{definition}

\begin{lemma}\label{fullness} Suppose that \ $\langle\alpha_1,\dots,\alpha_k\rangle$ \ is an increasing sequence of
ordinals, where \ $k\in\omega$. \ For each \ $m\in\omega$ \ there is an $m$--full sequence \
$\langle\beta_0,\dots,\beta_h\rangle$ \ such that \ $\{\alpha_1,\dots,\alpha_k\}\subseteq\{\beta_1,\dots,\beta_h\}$
\ and \ $\alpha_k=\beta_h$.
\end{lemma}

\begin{proof} See Dodd's proof of Lemma 7.12 of \cite{Dodd}.
\end{proof}

\begin{definition} Let \ $\langle \alpha_1, \dots, \alpha_k
\rangle$ \ and \ $\langle \beta_1, \dots, \beta_k
\rangle$ \ be increasing sequences of ordinals. For \ $m\in\omega$, \ we say that 
$\langle \alpha_1, \dots, \alpha_k
\rangle\sim_m\langle \beta_1, \dots, \beta_k
\rangle$ \ if and only if
\be
\item  $\langle \alpha_1, \dots, \alpha_k \rangle$ \ and \ $\langle \beta_1, \dots, \beta_k
\rangle$ \ are \ $m$--full, and
\item $(\forall j\le k) [\text{ch}_m(\alpha_j) = \text{ch}_m(\beta_j)]$.
\ee
\end{definition}

We now quote a technical lemma of
Dodd \cite[Lemma 7.17]{Dodd}, whose proof easily generalizes to ``premice above the reals'', and a corollary on the
existence of indiscernibles. First, given an iterable premouse \ $\mouse$ \ let \ $\pi_{0\alpha}\colon \mouse\mapsigma{1}
\mouse_\alpha$ \ be the premouse embedding of \ $\mouse$ \ into its premouse iterate \ $\mouse_\alpha$. \ Similarly,
let  \ $\kappa_\alpha = \pi_{0\alpha}(\kappa^{\mouse})$.
\begin{lemma} Suppose \ $\varphi$ \ is a \ $\Sigma_{m+1}$ \
formula of two free variables. Then there is a \ $\Sigma_{m+1}$ \ formula \ $\varphi^*$ \ of three
free variables such that, for any iterable premouse \ $\mouse$ \ and for all \ $a\in M$,
\[\mouse_\theta \models \varphi(\langle \kappa_{\alpha_1},\dots,\kappa_{\alpha_k}\rangle,
\pi_{0\theta}(a)) \ \text{\,iff\,}\ \mouse \models \varphi^*(\langle
\text{ch}_m(\alpha_1),\dots,\text{ch}_m(\alpha_k)\rangle,\text{ch}_m(\theta), a)\]
whenever \ $\langle \alpha_1,\dots,\alpha_k, \theta\rangle$ \ is \ $m$--full. 
\end{lemma}

\begin{corollary} Let \ $\mouse$ \ be an iterable premouse and let \
$\varphi(v_1,\dots,v_k,v_{k+1})$ \ be a  \ $\Sigma_{m+1}$ \ formula. Then for all \  $\langle
\alpha_1, \dots, \alpha_k, \theta \rangle\sim_m\langle \beta_1, \dots, \beta_k, \theta
\rangle$, \ and for all \ $a\in M$,
\[\mouse_\theta \models \varphi(\kappa_{\alpha_1},\dots,\kappa_{\alpha_k},\pi_{0\theta}(a)) 
\leftrightarrow \varphi(\kappa_{\beta_1},\dots,\kappa_{\beta_k}, \pi_{0\theta}(a))\]
where \ $\pi_{0\theta}\colon\mouse\mapsigma{1}\mouse_\theta$ \ is the premouse embedding.
\end{corollary}

We now have indiscernibles for certain premouse and mouse iterates.
\begin{corollary}\label{helper} Let \ $\mouse$ \ be an iterable premouse. Suppose that \ $\theta$ \ is \
$m$--good, \ and let \ $I_m=\{\,\kappa_\beta : \beta \text{ is $m$--good} \land \beta <
\theta\,\}$. \ Then \ $I_m$ \ is a set of order \
$\Sigma_{m+1}(\mouse_\theta,\{\,\pi_{0\theta}(a) : a\in M\,\})$ \ indiscernibles.
\end{corollary}

\begin{corollary}\label{liftup} Let \ $\mouse$ \ be a mouse with \ $n=n(\mouse)$. \ Let \ $\pi_{0\theta}\colon
\mouse\to \mouse_\theta$ \ be the mouse embedding of \ $\mouse$ \ into its \ $\theta^{\,\ul{\text{th}}}$ \ mouse
iterate \ $\mouse_\theta$. \ Let \ $I_k=\{\,\kappa_\alpha : \alpha \text{ is $k$--good} \land \alpha < \theta\,\}$, \
where  \ $\kappa_\alpha = \pi_{0\alpha}(\kappa^{\mouse})$. \ Then \ $I_k$ \ is a set of order \
$\Sigma_{n+(k+1)}(\mouse_\theta,\{\,\pi_{0\theta}(a) : a\in M\,\})$ \ indiscernibles.
\end{corollary}

\begin{proof} Let  \ $\pi^n_{0\theta}\colon {\mouse^n}\mapsigma{1} \mouse^n_\theta$ \ be the premouse embedding
of \ $\mouse^n$ \ into its \ $\theta^{\,\ul{\text{th}}}$ \ premouse iterate \ $\mouse^n_\theta$. \ By Corollary
\ref{helper} \ $I_k$ \ is a set of \ $\Sigma_{k+1}(\mouse^n_\theta,\{\pi^n_{0\theta}(a) : a\in
M^n\})$ \ indiscernibles. By Lemma \ref{EOEL} and its proof, we can construct the map \
$\pi^{n-1}_{0\theta}\colon\mouse^{n-1}\to\mouse^{n-1}_\theta$. \ By applying the idea in the proof of Claim
1 of Lemma \ref{EOEL}, we will now show that \ $I_k$ \ is a set of \
$\Sigma_{1+(k+1)}(\mouse^{n-1}_\theta,\{\pi^{n-1}_{0\theta}(a) : a\in M^{n-1}\})$ \ indiscernibles. 
By repeating this result we will get that \ $\pi_{0\theta}=\pi^0_{0\theta}$ \ and that \ $I_k$ \ is
a set of \ $\Sigma_{n+(k+1)}(\mouse_\theta,\{\pi_{0\theta}(a) : a\in M\})$ \ indiscernibles, as desired. 

In the interest of simplifying notation, we will assume that \ $k=1$.  \ So, we have that
\begin{equation}\label{indisc}\text{$I_k$ \ is a set of \ $\Sigma_{2}(\mouse^n_\theta,\{\pi^n_{0\theta}(a) : a\in
M^n\})$ \ indiscernibles.}\end{equation} Using the notation in the proof of Claim
1 of Lemma \ref{EOEL}, we will show that \ $I_k$ \ is also a set of \
$\Sigma_{3}(\mouse^{n-1}_\theta,\{\pi^{n-1}_{0\theta}(a) : a\in M^{n-1}\})$ \ indiscernibles.  
For notational convenience let \ $\fmouse=\mouse^{n-1}_\theta$. \ Let \ $\psi(v_0,\dots,v_i,v_{i+1},v_{i+2})$ \ be a
\ $\Sigma_1$ \ formula in the language \ $\Lng_{n-1}$. \ Let \
$\varphi(v_0,\dots,v_{i-1},v_i)$ \ be the \ $\Sigma_3$ \ formula \ $(\exists v_{i+1})(\forall
v_{i+2})\psi(v_0,\dots,v_{i-1},v_i,v_{i+1},v_{i+2})$. \ Let \ $a$ \ be arbitrary element in \
$M^{n-1}$. \ Let \ $\kappa_{\alpha_0}<\dots<\kappa_{\alpha_{i-1}}$ \ and \
$\kappa_{\beta_0}<\dots<\kappa_{\beta_{i-1}}$ \ be taken from \ $I_k$. \ We show that 
\[\fmouse \models \varphi(\kappa_{\alpha_0},\dots,\kappa_{\alpha_{i-1}},\pi^{n-1}_{0\theta}(a)) 
\leftrightarrow \varphi(\kappa_{\beta_0},\dots,\kappa_{\beta_{i-1}}, \pi^{n-1}_{0\theta}(a)).\]
By (1) of Lemma \ref{MET} and the definition of \ $\pi^{n-1}_{0\theta}$, \ there exist \ $x_i\in\R^{\mouse^n_\theta}$
\ and \  $\gamma_i\in\omega\rho^{}_{\mouse^n}$ \ such that \ $\pi^{n-1}_{0\theta}(a)=h^{}_{\fmouse\,}(x_i,\langle
\pi^{n}_{0\theta}(\gamma_i),p^{}_{\fmouse}\,\rangle)$. \ Let \ $y\in\R^{\mouse^n_\theta}$ \ be such that
 \ $\gamma=h^{}_{\fmouse\,}(y,\langle \gamma,p^{}_{\fmouse}\,\rangle)$ \ for all \
$\gamma\in\omega\rho^{}_{\fmouse}$ \ (recall that \ $I_k\subseteq \omega\rho^{}_{\fmouse}$). \
 Let \ $x=\langle y,y,\dots,y,x_i\rangle\in\R^{\mouse^n_\theta}$ \ and \
$s=\langle\kappa_{\alpha_0},\dots,\kappa_{\alpha_{i-1}},\pi^{n}_{0\theta}(\gamma_i)\rangle\in
(\omega\rho^{}_{\fmouse})^{<\omega}$. \ Recalling (\ref{domain}) and Remark
\ref{code},  let \ $\chi(x,s)$ \ be the formula 
\[\begin{aligned}
(\exists x_{i+1})&(\exists\gamma_{i+1})(\forall x_{i+2})(\forall\gamma_{i+2}) 
(D(i+1,x,s)\land D(1,\langle x_{i+1}\rangle,\langle\gamma_{i+1}\rangle)\\
&\land[ D(1,\langle x_{i+2}\rangle,\langle\gamma_{i+2}\rangle)\rightarrow
A_n(\ulcorner\varphi^+\urcorner{^\frown}(x{^\frown}\langle x_{i+1},x_{i+2}\rangle),s{^\frown}\langle
\gamma_{i+1},\gamma_{i+2}\rangle) ]).\\
\end{aligned}\]
The above formula \ $\chi(x,s)$ \ is \ $\Sigma_2$ \ in the language \ $\Lng_n$ \ 
and the following holds
\begin{alignat*}{2}
&\fmouse\models\varphi(\kappa_{\alpha_0},\dots,\kappa_{\alpha_{i-1}},\pi^{n-1}_{0\theta}(a))&&\\ 
&\iff \fmouse\models\varphi[h^{}_{\fmouse}\,(y,\langle
\kappa_{\alpha_0},p^{}_{\fmouse}\,\rangle),\dots, h^{}_\fmouse\,(x_i,\langle
\pi^{n}_{0\theta}(\gamma_i),p^{}_{\fmouse}\,\rangle)] &&\quad\text{by Lemma \ref{MET}}\\ 
&\iff
\mouse^n_\theta\models\chi(x,\langle\kappa_{\alpha_0},\dots,\kappa_{\alpha_{i-1}},\pi^{n}_{0\theta}(\gamma_i)\rangle)&&\quad\text{by
Lemma \ref{MET}}\\ 
&\iff
\mouse^n_\theta\models\chi(x,\langle\kappa_{\beta_0},\dots,\kappa_{\beta_{i-1}},\pi^{n}_{0\theta}(\gamma_i)\rangle)&&
\quad\text{by (\ref{indisc})}\\ 
&\iff \fmouse\models\varphi[h^{}_{\fmouse}\,(y,\langle
\kappa_{\beta_0},p^{}_{\fmouse}\,\rangle),\dots,h^{}_\fmouse\,(x_i,\langle
\pi^{n}_{0\theta}(\gamma_i),p^{}_{\fmouse}\,\rangle)] &&\quad\text{by Lemma \ref{MET}}\\ 
&\iff\fmouse\models\varphi(\kappa_{\beta_0},\dots,\kappa_{\beta_{i-1}},\pi^{n-1}_{0\theta}(a))&&\quad\text{by Lemma
\ref{MET}}.
\end{alignat*}
This completes our proof of the Corollary.
\end{proof}

The next result is essentially Corollary 2.42 of \cite{Crcm}. We have just specified some relevant parameters. 
\begin{lemma}\label{defone} Let \ $\mouse$ \ be a mouse with \ $n=n(\mouse)$. \ Let \
$\pi_{0\theta}\colon \mouse\to \mouse_\theta$ \ be the mouse embedding of \ $\mouse$ \ into its \
$\theta^{\,\ul{\text{th}}}$ \ mouse iterate \ $\mouse_\theta$. \ Similarly, let 
\ $\overline{\pi}_{0\theta}\colon \overline{\mouse}\to \overline{\mouse}_\theta$ \ be the premouse embedding of
\ $\overline{\mouse}$ \ into its \
$\theta^{\,\ul{\text{th}}}$ \ premouse iterate \ $\overline{\mouse}_\theta$.
\ Suppose that \ $\theta$ \ is \ $m$--good  for each \ $m\in\omega$ \
and \ let \ $I_m=\{\,\kappa_\alpha : \alpha \text{ is $m$--good} \land \alpha < \theta\,\}$, \
where  \ $\kappa_\alpha =
\pi_{0\alpha}(\kappa^{\mouse})$. \  Then
\be
\item[(1)] $I_0$ \ is uniformly \ $\Pi_1(\overline{\mouse}_\theta,\{p^{}_{\,\overline{\mouse}_\theta},\kappa_0\})$
\item[(2)] $I_0$ \ is \ $\boldface{\Pi}{n+1}({\mouse}_\theta)$
\item[(3)] for \  $\beta<\theta$, \ $I_m\setminus \kappa_{\beta+1}$ \ is a set of order \
$\Sigma_{m+1}(\overline{\mouse}_\theta,\{\overline{\pi}_{\beta\theta}(q) : q\in \overline{M}_\beta\})$ \ indiscernibles
and a set of order \
$\Sigma_{n+m+1}(\mouse_\theta,\{\pi_{\beta\theta}(a) : a\in M_\beta\})$ \ indiscernibles
\item[(4)] $I_m$ \ is uniformly \ $\Sigma_\omega(\overline{\mouse}_\theta,\{\kappa_0\})$
\item[(5)] $I_m$ \ is \ $\boldface{\Sigma}{\omega}({\mouse}_\theta)$. 
\ee
\end{lemma}

\begin{proof} (1) is established in the proof of Corollary 2.42 of \cite{Crcm}. Lemma \ref{thdgoose} implies (2).
For (3), let \ $\beta<\theta$ \ and let \ $\theta'$ \ be the unique ordinal such that \ $\beta+\theta'=\theta$. \
Since \ $\theta$ \ is \ $m$--good  for each \ $m\in\omega$, \ it follows that \ $\theta'$ \ is also \ $m$--good  for
all \ $m\in\omega$. \ Corollaries \ref{helper} and \ref{liftup} imply  (3). Now, (1) implies (4) because
the parameter \ $p^{}_{\,\overline{\mouse}_\theta}$ \ is definable over \ $\overline{\mouse}_\theta$ \ and because the
notion of being ``$m$--good'' is definable.\footnote{Note: $\alpha$ is $m$--good iff $\kappa_\alpha$  is
$m$--good relative to the ordinals in $I_0$; that is, $I_i=\text{limit points of }I_{i-1}$.} Finally, (4) and Lemma
\ref{thdgoose} imply (5).
\end{proof}

When \ $\mouse$ \ is a mouse with core \ $\core=\core(\mouse)$, \ 
we know that there is an ordinal \ $\theta$ \ such that the mouse iterate \ $\core_\theta$ \ is such that \
$\core_\theta=\mouse$. \ If  \ $\mouse$ \ is a proper initial segment of an iterable premouse, then one can easily predict
this ordinal \ $\theta$; \ namely, either \ $\theta=0$ \ or \ $\theta=\kappa^\mouse$. \ Lemma \ref{deftwo} below  will
establish this result together with some other observations that will be used in Parts II \& III. 

First, we shall define \ $\Sigma_\omega^\mu$ \ formulae and prove some relevant propositions and a corollary.  Recall that the
premouse \ $\nouse=(N,\R^\nouse,\kappa,\mu)$  \ is a model of the language \
$\Lng_0=\{\in,\underline{\R},\underline{\kappa},\mu\}$ \ where \ $\mu$ \ is a predicate. For the remainder of this
subsection we shall let \ $\Lng=\Lng_0$. \ We will now add a quantifier to the language \ $\Lng$. \ Since the quantifier
extends the predicate \ $\mu$ \ in our intended structures, we shall use the same symbol \ $\mu$ \ for this quantifier.
We shall denote this expanded language by \ $\Lng^\mu$ \ and write \
$\Sigma_\omega^\mu$ \ for the formulae in this expanded language. For \ $\gamma$ \ such that \
$\kappa<\gamma<\widehat{\OR}^{\nouse}$, \ let \ $\mu^{\gamma+1}=N^{\gamma+1}\cap\mu$. \ Then \
$(\nouse^\gamma,\mu^{\gamma+1})$ \ is an \ $\Lng^\mu$ \ structure, where the new
quantifier symbol is to be interpreted by \ $\mu^{\gamma+1}$. \ That is, \ $(\nouse^\gamma,\mu^{\gamma+1})\models
(\mu\,\varkappa\in\kappa)\psi(\alpha)$ \ if and only if \ $\{\varkappa\in\kappa : (\nouse^\gamma,\mu^{\gamma+1})\models
\psi(\varkappa)\}\in\mu^{\gamma+1}$.

Letting \ $L=\bigcup\limits_{\gamma\in\OR}J_\gamma$, \ where  \ $J_\gamma$ \ is the \ $\gamma^{\,\ul{\text{th}}}$ \ level of the
Jensen hierarchy for the constructible universe \ $L$, \ recall that 
\[\forall x [x\in\pow(J_\gamma)\cap J_{\gamma+1}\implies x\in\boldface{\Sigma}{\omega}(J_\gamma)].\]
This fact and its proof easily generalize to give the following two propositions.
\begin{proposition}\label{quant} Let \ $\nouse=(N,\R^\nouse,\kappa,\mu)$ \ be a premouse and let \ $\gamma$ \ be such that \
$\kappa<\gamma<\widehat{\OR}^{\nouse}$. \ Then 
\[\forall x [x\in\pow(N^\gamma)\cap N^{\gamma+1}\implies
x\in\boldfacemu{\Sigma}{\omega}(\nouse^\gamma,\mu^{\gamma+1})].\]
\end{proposition}

\begin{proposition}\label{theendprop} Let \ $\phi(v)$ \ be a \ $\Sigma_1$ \ formula of \ $\Lng$. \ For each \ $k\in\omega$
\ there is a formula \ $\psi_k(v)$ \ in \ $\Sigma_\omega^\mu$ \ such that 
\[(S_{\gamma+k}^{\nouse}(\underline{\R}),\R^\nouse,\kappa,\mu)\models\phi(a)\iff
(\nouse^\gamma,\mu^{\gamma+1})\models\psi_k(a)\]
for all \ $a\in N^\gamma$, \ whenever \ $\nouse=(N,\R^\nouse,\kappa,\mu)$ \ is a premouse and \
$\kappa<\gamma<\widehat{\OR}^{\nouse}$.
\end{proposition}

\begin{definition}\label{preddef} Let \ $\nouse=(N,\R^\nouse,\kappa,\mu)$ \ be a premouse and 
let \ $\kappa<\gamma<\widehat{\OR}^{\nouse}$. \ We shall say that \ $\mu^{\gamma+1}$ \ is {\it predictable\/} if 
the following holds: For each \ $\Lng$ \ formula \ $\chi(v_0,v_1,\dots,v_k)$ \ there is another  \ $\Lng$ \
formula \ $\psi(v_0,v_1,\dots,v_k)$ \ and a \ $d\in N^\gamma$ \ such that for all \ $a_1,\dots,a_k\in N^\gamma$ \ 
\[B_{a_1,\dots,a_k}\in\mu^{\gamma+1}\iff \nouse^\gamma\models\psi(d,a_1,\dots,a_k)\]
where \ $B_{a_1,\dots,a_k}=\{\varkappa\in\kappa : \nouse^\gamma\models \chi(\varkappa,a_1,\dots,a_k)\}$. 
\end{definition}

\begin{lemma}\label{predlemma} Let \ $\nouse=(N,\R^\nouse,\kappa,\mu)$ \ be a premouse and \
$\kappa<\gamma<\widehat{\OR}^{\nouse}$. \ If \ $\mu^{\gamma+1}$ \ is predictable, then
\ $\pow(N^\gamma)\cap N^{\gamma+1}\subseteq \boldface{\Sigma}{\omega}(\nouse^\gamma)$.
\end{lemma}

\begin{proof} Let \ $A\subseteq N^\gamma$ \ be such that \ $A\in N^{\gamma+1}$. \ By Proposition \ref{quant} \
$A\in{\boldfacemu{\Sigma}{\omega}}(\nouse^\gamma,\mu^{\gamma+1})$. \ Let \ $\psi\in\Sigma_\omega^\mu$ \ and \
$a_1,\dots,a_k\in N^\gamma$ \ be such that for all \ $a\in N^\gamma$ \ 
\begin{equation}
a\in A\iff (\nouse^\gamma,\mu^{\gamma+1})\models\psi(a,a_1,\dots,a_k).\label{wendy}
\end{equation}
Since \ $\mu^{\gamma+1}$ \ is predictable, we will exhibit a \ $\Sigma_\omega$ \ formula \
$\varphi(v_0,\dots,v_k,v_{k+1})$ \ and a \ $d\in N^\gamma$ \ such that for all \ $a,a_1,\dots,a_k\in N^\gamma$
\begin{equation}
(\nouse^\gamma,\mu^{\gamma+1})\models\psi(a,a_1,\dots,a_k) \iff \nouse^\gamma\models
\varphi(a,a_1,\dots,a_k,d).\label{wendell}
\end{equation} The formula \ $\varphi$ \ is constructed by induction on the complexity of \ $\psi$ \ using the
assumption that \ $\mu^{\gamma+1}$ \ is predictable. The \ $\mu$--quantifier case is the only inductive case that
requires checking; that is, suppose that \ $\psi(v_0,\dots,v_k,v_{k+1})$ \ has the form \
$(\mu\,\varkappa\in\underline{\kappa})\chi(\varkappa,v_0,\dots,v_k,v_{k+1})$. \ By the induction hypothesis there is a
formula \ $\varphi'\in\Sigma_\omega$ \ and a \ $d'\in N^\gamma$ \ such that for all \ $\varkappa,a,a_1,\dots,a_k\in
N^\gamma$
\[(\nouse^\gamma,\mu^{\gamma+1})\models\chi(\varkappa,a,a_1,\dots,a_k)\iff \nouse^\gamma\models
\varphi'(\varkappa,a,a_1,\dots,a_k,d').\]
By predictability, we can now obtain the required formula \ $\varphi$ \ and the parameter \ $d\in N^\gamma$ \ that
verifies (\ref{wendell}). Therefore, (\ref{wendy}) and (\ref{wendell}) imply that \
$A\in\boldface{\Sigma}{\omega}(\nouse^\gamma)$.
\end{proof}

\begin{lemma}\label{deftwo}
Let  \ $\mouse$ \ be a mouse with core \ $\core$ \ and let \ $n=n(\mouse)$. \ Suppose that \ $\rho_{\mouse}^{n+1}
< \kappa^{\mouse}$ \ and that \ $\mouse=\nouse^\gamma$ \ for an iterable premouse \ $\nouse$ \ where \
$\kappa^{\mouse}=\kappa^{\nouse}=\kappa$
\ and \ 
$\kappa<\gamma<\widehat{\OR}^{\nouse}$. \ Let
\ $\theta$
\ be such that the mouse iterate \ $\core_\theta=\mouse$. \ Then 
\be
\item $I_m^{\overcore_\theta}\in\mu^\nouse$ \ for all \ $m\in\omega$
\item $\theta=\kappa$
\item $\theta$ \ is \ $m$--good  for all \ $m\in\omega$
\item $\theta$ \ is a multiple of \ $\omega^\omega$
\item $\pow(M)\cap N^{\gamma+1}\subseteq \boldface{\Sigma}{\omega}(\mouse)$.
\ee
\end{lemma}

\begin{proof} Since  \ $\mouse$ \ is a proper initial segment of \ $\nouse$, \ it follows that \
$\mouse\in\nouse$ \ and so, \ $\overmouse\in\nouse$. \ Also, because \ $\core_\theta=\mouse$ \  we have that \
$\overcore_\theta=\overmouse$. \ Let \ $\kappa_0=\kappa^\core$. \ Lemma \ref{defone} implies that \ 
$I_m^{\overcore_\theta}$ \ is uniformly definable (in the constant \ $\kappa_0$) \ over
\ $\overcore_\theta$ $(=\overmouse)$ \ and hence, \ $I_m^{\overcore_\theta}\in\nouse$ \ for each \ $m\in\omega$. \
Clearly, \ $\kappa_0\le\kappa^\mouse=\kappa^\nouse$. \ We first show that \ $\kappa_0<\kappa^\nouse$. \ Suppose, for a
contradiction, that \ $\kappa_0=\kappa^\mouse$. \ Thus, \ $\overcore=\overmouse$. \ It follows (see Definition
\ref{skolemdef} and subsection \ref{coremice}) that there is a function
\ $f$ \ in \ $\nouse$ \ such that \ $f\colon\rho_{\mouse}^{n+1}\times\R^\nouse\maps{onto}\kappa^\nouse$ \ where \ 
$\rho_{\mouse}^{n+1}< \kappa^{\nouse}$ \ by assumption. This contradiction shows that \ $\kappa_0<\kappa^\nouse$.

Now, let \ $\varkappa>\abs{\nouse}$ \ be a regular cardinal\footnote{This appeal to $\AC$ is removable. Let
$M=L(\nouse)$. One can prove the lemma in a $\ZFC$--generic extension $M[G]$  and thus, by absoluteness, the result
holds in $L(\nouse)$ and hence in $V$.} and  let \ $\mouse_\varkappa$ \ and $\core_\varkappa$ \ be the respective
mouse iterates of \ $\mouse$ \ and \ $\core$. \  Let \ $\nouse_\kappa$ \ be the premouse iterate of \ $\nouse$ \ with
the corresponding premouse embedding \ $\pi\colon\nouse\mapsigma{1} \nouse_\kappa$. \ We conclude that
\ $\underline{\kappa}^{\mouse_\varkappa}=\underline{\kappa}^{\core_\varkappa}=\underline{\kappa}^{\nouse_\varkappa}=\varkappa$,
\ $\overcore_\varkappa=\overmouse_\varkappa$ \ and \ $\mouse_\varkappa$ \ is a proper initial segment of \
$\nouse_\varkappa$. \ We can also deduce that \ $\overcore_\varkappa=\pi(\overcore_\theta)$.
\ Lemma \ref{defone} implies that \ $I_m^{\overcore_\varkappa}$ \ is uniformly definable (in the constant \ $\kappa_0$) \ over
\ $\overcore_\varkappa$ $(=\overmouse_\varkappa)$ \ and hence, \ $I_m^{\overcore_\varkappa}\in\nouse_\varkappa$ \ for all \
$m\in\omega$. \ Since  \ $I_m^{\overcore_\varkappa}$ \ is a club in \ $\varkappa$ \ for each \ $m\in\omega$, \ 
it follows that \ $\nouse_\varkappa\models I_m^{\overcore_\varkappa}\in\mu^{\nouse_\varkappa}$. \ Because \
$\pi\colon\nouse\mapsigma{1} \nouse_\kappa$, \
$\overcore_\varkappa=\pi(\overcore_\theta)$ \ and \ $\pi(\kappa_0)=\kappa_0$, \ it now follows that  \ $\nouse\models
I_m^{\overcore_\theta}\in\mu^{\nouse}$ \ for all \ $m\in\omega$. \ Thus,
\ $\theta=\kappa^{\mouse}=\kappa^{\nouse}$ \ and so, \ $\theta$ \ is \ $m$--good  for all \ $m\in\omega$. \ Hence, \ $\theta$ \
is a multiple of \ $\omega^m$ \ for all \ $m\in\omega$ \ and therefore, \ $\theta$ \ is a multiple of \ $\omega^\omega$. \
After we prove the following claim, we will show that \ $\pow(M)\cap N^{\gamma+1}\subseteq
\boldface{\Sigma}{\omega}(\mouse)$.
\begin{claiim} $\mu^{\gamma+1}$ \ is predictable.  
\end{claiim}

\begin{proof}[Proof of Claim]
Let \ $\chi(v_0,v_1,\dots,v_k)\in\Sigma_{n+m+1}$ \ be an \ $\Lng$ \ formula where \ $n=n(\mouse)$. \  We must find 
an \ $\Lng$ \ formula \ $\psi(v_0,v_1,\dots,v_k)$ \ and a \ $d\in M$ \ such that for all \ $a_1,\dots,a_k\in M$ \ 
\begin{equation}
B_{a_1,\dots,a_k}\in\mu^{\gamma+1}\iff \mouse\models\psi(d,a_1,\dots,a_k)\label{idiot1}
\end{equation}
where \ $B_{a_1,\dots,a_k}=\{\varkappa\in\kappa : \mouse\models \chi(\varkappa,a_1,\dots,a_k)\}$. \ 
Let \ $I_m=I_m^{\overcore_\theta}$. \ Lemma \ref{defone}(5) implies that there is a \ $\psi\in\Sigma_\omega$ \ and \
$d\in M$ \ such that 
\[\mouse\models(\exists\eta\in\kappa)(\forall\varkappa)\left[\varkappa\in I_m\setminus\eta\rightarrow
\chi(\varkappa,v_1,\dots,v_k)\right]\iff \mouse\models\psi(d,v_1,\dots,v_k)\]
for all \ $v_1,\dots,v_k\in M$.

Now let \ $a_1,\dots,a_k\in M$. \ We will show that (\ref{idiot1}) holds.  
Let
\begin{equation}\label{ocsystem}\overcoresystem\end{equation}
be the premouse iteration of \ $\overcore$ \ and let
\begin{equation}\label{csystem}\coresystem\end{equation}
be the mouse iteration of \ $\core$ \ via Lemma \ref{EOEL}, the extension of embeddings lemma. In addition, let \
$\kappa_\alpha = \pi_{0\alpha}(\underline{\kappa}^{\core})$ \ for \ $\alpha\in \OR$. \ Recall that \
$\core_\theta=\mouse$, \ $\theta=\kappa$, \ and \ $\theta$ \ is a multiple
of \ $\omega^\omega$. \  Because \ $\overcore_\theta$ \ is a direct limit in the system (\ref{ocsystem}), it follows
that \ $\core_\theta$ \ is a direct limit in the system (\ref{csystem}). Thus, there is an ordinal \ $\xi<\theta$ \
and \ $c_1,\dots,c_k\in C_\xi$ \ (the domain of \ $\core_\xi$) \ such that 
\ $a_i= \pi_{\xi\theta}(c_i)$ \ for \ $1\le i\le k$. \ Let \ $\eta=\xi+1$. \ Lemma
\ref{defone}(3) implies that \ $I_m\setminus\eta$ \ is a set of order \
$\Sigma_{n+(m+1)}(\core_\theta,\{\pi_{\xi\theta}(c) : c\in C_\xi\,\})$ \ indiscernibles.
Hence, for all \ $\varkappa, \varkappa'\in I_m\setminus\eta$,
\[\mouse\models\chi(\varkappa,a_1,\dots,a_k)\leftrightarrow\chi(\varkappa',a_1,\dots,a_k).\]
Since \ $I_m\in\mu^{\gamma+1}$, \ (\ref{idiot1}) now follows.
This completes the proof of the Claim.\end{proof} 
 
Lemma \ref{predlemma} and the above Claim now imply that $\pow(M)\cap N^{\gamma+1}\subseteq
\boldface{\Sigma}{\omega}(\mouse)$.\end{proof}

\subsubsection{A minimal criterion for mouse iterability}\label{minmouseiterable}

Assuming acceptability, the following theorem gives a ``coarse'' condition  for mouse iterability.
\begin{theorem}\label{intismi} Let \ $\mouse=(M,\R^\mouse,\kappa,\mu)$ \ be an iterable premouse. Let \ $\gamma$ \ be such
that \ $\kappa<\gamma<\widehat{\OR}^{\mouse}$. \ Suppose that \ $\mouse^\gamma$ \ is acceptable and
critical.  Then \ $\mouse^\gamma$ \ is a mouse.
\end{theorem}

\begin{proof} Since \ ${\mouse^\gamma}$ \ is acceptable and critical, we just need to show that \ ${\mouse^\gamma}$ \ 
is \ $n(\mouse^\gamma)$--iterable (see Definitions \ref{niteratable} and \ref{mouseiteratable}). Let \
\[\premouseiteration{\mouse}{\alpha}{\beta}\]
be the premouse iteration of \ $\mouse$. \ For each ordinal \ $\alpha$, \ we shall also write \ 
$\mouse^{\gamma_\alpha}=\mouse^{\pi_{0\alpha}(\gamma)}_\alpha$ \ where \ $\mouse^{\pi_{0\alpha}(\gamma)}_\alpha=
(\mouse_\alpha)^{\pi_{0\alpha}(\gamma)}$.

Let \ $n=n(\mouse^\gamma)$. \ In addition, \ let \
$\overline{\mouse^{\gamma_\alpha}}=\left({\mouse^{\gamma_\alpha}}\right)^{n}$ \ and let \ $\overline{M^{\gamma_\alpha}}$ \
denote the domain of this structure. For all ordinals \ $\alpha\le\beta$, \ the following hold:
\be
\item $\mouse^{\gamma_\alpha}\in M_\alpha$
\item $\overline{\mouse^{\gamma_\alpha}}\in M_\alpha$
\item $\pi_{\alpha\beta}\left(\mouse^{\gamma_\alpha}\right)=\mouse^{\gamma_\beta}$
\item $\pi_{\alpha\beta}\left(\,\overline{\mouse^{\gamma_\alpha}}\,\right)=\overline{\mouse^{\gamma_\beta}}$
\item $\pi_{\alpha\beta}\colon\mouse^{\gamma_\alpha}\mapsigma{\omega}\mouse^{\gamma_\beta}$
\item $\pi_{\alpha\beta}\colon\overline{\mouse^{\gamma_\alpha}}\mapsigma{\omega}\overline{\mouse^{\gamma_\beta}}$
\item $\mouse^{\gamma_\alpha}$ \ is acceptable, critical and \ $n=n(\mouse^{\gamma_\alpha})=n(\mouse^\gamma)$.
\ee
Let \ 
\[\premouseiterationover{\left(\,\overline{\mouse^{\gamma}}\,\right)}{\alpha}{\beta}\]
be the premouse iteration of \ $\overline{\mouse^{\gamma}}$.
\ It can be shown that for each ordinal \ $\alpha$ \ there exists an embedding \ $\sigma_\alpha\colon
\left(\,\overline{\mouse^{\gamma}}\,\right)_\alpha\mapsigma{0}
\overline{\mouse^{\gamma_\alpha}}$ \ (see, for example, \cite[Lemma 10.32]{Dodd}). Now, since \
$E_n^{\overline{\mouse^{\gamma_\alpha}}}$ \ is well-founded (because \ $\mouse^{\gamma_\alpha}$ is transitive), it
follows that \ $E_n^{\left(\,\overline{\mouse^{\gamma}}\,\right)_\alpha}$ \ is well-founded for all \ $\alpha$. \
Therefore, \ $\mouse^\gamma$ \ is a mouse.
\end{proof}

Our next theorem will establish a ``minimal'' relative criterion ensuring that a model of \
$\mathcal{T}^n$ \ is \ $n$--iterable. This criterion will be used to produce scales definable over a
weak real mouse.  First, we will give some definitions. 
\begin{definition}  Let \ $\mouse$ \ be a model of \ $\mathcal{T}^m$. \ The set of {\it
ordinal codes\/} in \ $\mouse$, \ denoted by \ $\cOR^\mouse$, \ is defined by 
\[\cOR^\mouse = \{ (x,\gamma)\in (\R\times {\OR})^\mouse : \mouse \models (x,\gamma)\,
\ul{E}_m\,(x,\gamma) \}.\] 
The equivalence relation \ $\equiv$ \ on the set \
$\cOR^\mouse$ \ is given by
\[(x,\gamma)\equiv (y,\beta) \iff \mouse \models (x,\gamma)\,\ul{E}_m\,(y,\beta) \land
(y,\beta)\,\ul{E}_m\,(x,\gamma),\]
and the equivalence class of an ordinal code \ $(x,\gamma)$ \ in \ $\mouse$ \ is denoted by \
$[(x,\gamma)]=[(x,\gamma)]_\equiv$. \ Define the set of equivalence classes as
\[\cOR^\mouse_{\equiv} = \{ [(x,\gamma)] : (x,\gamma)\in \cOR^\mouse\}\]
and define the relation \ $\mathcal{E}$ \ on \ $\cOR^\mouse_{\equiv}$ \ by
\[[(x,\gamma)]\,\mathcal{E}\,[(y,\beta)] \iff \mouse \models (x,\gamma)\,E_m\,(y,\beta).\]
\end{definition}

Let \ $\mouse$ \  be a model of \ $\mathcal{T}^m$ \ and define   
\[F^{\mouse}=\{\,f\in M :  \exists n\in\omega \
\mouse \models f\colon {^n}\ul{\kappa} \rightarrow\OR\,\}.\] 
For \ $f\in F^{\mouse}$, \ write \ $d(f) = n$ \ if and only if \ $n\in\omega$  \ and \ $\mouse \models
f\colon {^n}\ul{\kappa} \rightarrow\OR$. \ We shall assume the convention that \ $f\in F^{\mouse}$ \ and 
$d(f) = 0$ \ whenever \ $f\in {\OR}^{\mouse}$. \ Finally, for \
$n\in\omega$, \ define \  $F^{\mouse}_n = \{\,f\in F^{\mouse} : d(f) = n \,\}.$

\begin{definition}\label{FunnyF} Suppose that \ $\mouse$ \ is a model of \ $\mathcal{T}^m$. \
For \ $n\in\omega$ \ let 
\[\F_n^{\mouse} = \{ (x,f)\in(\R\times F_n)^{\mouse} :
\mouse \models \forall \xi_0,\dots,\xi_n\, (x,f(\xi_0,\dots,\xi_n))\in\cOR\}\] 
and let \ $\F^{\mouse} = \bigcup\limits_{n\in\omega}\F_n^{\mouse}$. \ 
For \ $\f = (x,f) \in \F^{\mouse}$ \ we let \ $\mathfrak{f}_0 = x$, \ $\mathfrak{f}_1 = f$, \ 
and we write \ $\f(\xi_0,\dots,\xi_{n-1}) = (x,f(\xi_0,\dots,\xi_{n-1}))$. \ In addition, we write \
$d(\f) = d(f)$.
\end{definition}

\begin{rmk} Given an \ $(x,f)\in(\R\times F_n)^{\mouse}$ \ such that \ $\mouse\models
\{\,\xi \in {^n\ul{\kappa}} : (x,f(\xi_0,\dots,\xi_{n-1}))
\in \cOR\,\}\in \mu_n,$ \ then \ $(x,f)$ \ is equivalent (modulo $\equiv_{\mu_n}$) to an \
$\f\in\F_n^{\mouse}$. \ In this case we shall implicitly assume that \ $(x,f)=\f$. \ In addition,
by \ $\R$--completeness in \ $\mouse$, \ {\it any\/}
function \ $h\colon {^{n}(\ul{\kappa}^{\mouse})}\rightarrow \cOR^\mouse$ \ in \ $\mouse$ \ is
equivalent, modulo \ $\mu_n$, \ to a function \ $\f\in \F^{\mouse}_n$.
\end{rmk}

Let \ $R$ \ be a rudimentary relation on \ $\cOR^{\mouse}$. \ Given
$\f, \g \in \F^{\mouse}$, \ let \ $n=d(\f)$ \ and \ $m=d(\g)$.  For any  \
$s\in{^n(n+m)}\!\uparrow$
\ and
\ for any \  $t\in{^m(n+m)}\!\uparrow,$ \ we shall write
\ $\mouse\models \,\f\,R^{s,t}\,\g$ \ 
if and only if
\[
\mouse \models\{\xi
\in{^{n+m}\ul{\kappa}}:\f(\xi_{s(0)},\dots,\xi_{s(n-1)})\,R\,\g(\xi_{t(0)},\dots,\xi_{t(m-1)})\}\in \mu_{n+m}.
\] 
Theorem \ref{thmttho} establishes a relative criterion assuring that a premouse \ $\mathcal{A}$ \ is iterable.
This criterion required the existence of an \ $\le$--extendible map \ $\sigma\colon F^{\mathcal{A}}\to F^\mouse$ \ (see
Definition \ref{defttho}) where \ $\mouse$ \ is ``premouse iterable.'' Our next definition generalizes this notion
and provides a sufficient condition for an  acceptable pure premouse to be \ $k$--iterable.
\begin{definition}\label{Extendible} Let \ $\mouse$ \ and \ $\mathcal{A}$ \ be models of \ $\mathcal{T}^k$. \ 
A map \ $\sigma\colon \F^{\mathcal{A}} \rightarrow \F^{\mouse}$ \ is said to be \
$\ul{E}_k$--{\it extendible}
\  if, for all \ $\f, \g\in \F^{\mathcal{A}}$
\be
\item $d(\f)=d(\sigma(\f))$
\item for all \ $s\in{^{d(\f)}(d(\f)+d(\g))}\!\uparrow$ \ and \ for all \ 
$t\in{^{d(\g)}(d(\f)+d(\g))}\!\uparrow$, 
\[\mathcal{A} \models \f\,\ul{E}_k ^{s,t}\,\g \iff \mouse \models
\sigma(\f)\,\ul{E}_k ^{s,t}\,\sigma(\g).\] 
\ee
\end{definition}

The following theorem (see \cite[Theorem 2.28]{Cfsrm}) establishes a relative criterion ensuring that a model of \
$\mathcal{T}^n$ \ is \ $n$--iterable and is used in Part II to produce scales definable over a
weak real mouse.

\begin{theorem}\label{criterion} Let \ $\nouse$ \ be a \ $k$--iterable premouse and 
let \ $\mouse=\nouse^k$. \ Suppose that \ $\mathcal{A}$ \ is a model of \ $\mathcal{T}^k$ \ and that \ $\sigma\colon
\F^{\mathcal{A}}
\rightarrow
\F^{\mouse}$ \ is \ $\ul{E}_k$--extendible.  Then \ $\mathcal{A}$ \ is
(isomorphic to) the \ $\Sigma_k$--code of a \ $k$--iterable pure premouse \ $\mathcal{B}$. 
\end{theorem}

\begin{remark}\label{clear} Let \ $\mouse$ \ be a mouse. When the
context is clear, we shall write \ $\Eq=\Eq_{n(\mouse)}$ \ and \ $E=E_{n(\mouse)}$.
\end{remark}

\begin{definition} Let \ $\mouse$ \ be a mouse. For \ $\f,\g\in \F^{\overmouse}_n$, \
where \ $n\in\omega$, \ let
\be
\item $\f\,\Eq_{\mu_n}\, \g \iff \overmouse \models
\{\langle \xi_0,\dots,\xi_{n-1}\rangle : \f(\xi_0,\dots,\xi_{n-1})\,\Eq\,
\g(\xi_0,\dots,\xi_{n-1})\}\in \mu_n$
\item $\f\,E_{\mu_n}\,\g \iff \overmouse \models
\{\langle \xi_0,\dots,\xi_{n-1}\rangle : \f(\xi_0,\dots,\xi_{n-1})\,E\,
\g(\xi_0,\dots,\xi_{n-1})\}\in \mu_n$
\item $\f\equiv_{\mu_n}\g \iff \f\,\Eq^{\overmouse}_{\mu_n}\,\g \text{ \ and \ }
\g\,\Eq^{\overmouse}_{\mu_n}\,\f$.
\ee
\end{definition}

For a mouse \ $\mouse$, \ $\equiv_{\mu_n}$ \ is an equivalence
relation on \ $\F^{\overmouse}_n$, and we let \ $[\f]_{\mu_n}$ \ represent the
equivalence class of \  $\f\in \F^{\overmouse}_n$. \ We assume the convention
that \ $\F^{\overmouse}_0={\cOR}^{\overmouse}$ \ and \ $\Eq^{\overmouse}_{\mu_0}\ =\ \Eq^{\overmouse}$ \
the ordering on \ ${\cOR}^{\overmouse}$.

The following result is Lemma 3.2 of \cite{Cfsrm} and will be used when we construct definable scales over a weak real mouse.
\begin{lemma} Suppose that \ $\mouse$ \ is a mouse with \ $m=n(\mouse)$
\ and let \ $n\in\omega$. \ Then \ $E^{\overmouse}_{\mu_n}$ \ is well-founded and hence,
\ $\Eq^{\overmouse}_{\mu_n}$ \ is a prewellordering on \ $\F^{\overmouse}_n$.
\end{lemma}

\begin{definition}\label{nicedef}
Let \ $\mouse$ \ be a mouse. For \ $\f\in \F_n^{\overmouse}$,
\ we shall write \ $\abs{[\f]_{\mu_n}}$ \ for the \ $E_{\mu_n}^{\overmouse}$--rank \ of \ $\f$. \ 
In addition, \ we say that \  $\tau\subseteq \F^{\overmouse}
\times \F^{\overmouse}$ \ is \ {\it nice} \ if \ $\tau$ \ is finite and \ $d(\h) = d(\f)$ \
for all \ $(\h,\f)\in\tau$. 
\end{definition}

This completes our review of the relative criterion for mouse iterability presented in Theorem \ref{criterion}. In
the next subsection we shall review two theorems concerning the definability of this criterion. 

\subsubsection{Defining $\ul{E}$--extendible maps}
We proved in \cite{Crcm} that \ $\Sigma_1(\mouse)$ \ has the scale property when \ $\mouse$ \ is an iterable real
premouse satisfying the axiom of determinacy. Our proof required us to show that the condition
``there exists a \ $\le$--extendible map \ $\sigma\colon F^\mouse\to F^\mouse$ \ such that \ $\sigma\supset\tau$'' \ is
definable over \ $\mouse$ \ for \ $\tau\in\mouse$. \ In Part II \cite{Part2} we shall prove that when  \ $\mouse$ \  is
a weak real mouse satisfying the axiom of determinacy, then \ $\boldface{\Sigma}{m}(\mouse)$ \ has the scale property
where \
$m=m(\mouse)$. \ Our method of constructing these scales requires us to show, for \ $\tau\in\overmouse$, \ that the
existence of an \ $\Eq$--extendible map \ $\sigma\colon\F^{\overmouse}\to\F^{\overmouse}$, \ such that \
$\sigma\supset\tau$, \  is definable over \ $\overmouse$. \ It will be necessary, however, to revise this iterability
condition slightly.\footnote{The structure $\overmouse$ does not necessarily satisfy $\AC$.}

\begin{definition}\label{quasi} Let \ $\mouse$ \ be a mouse. We say
that a relation \
$\Phi\subseteq \F^{\overmouse} \times  \F^{\overmouse}$ \ is an \ $\Eq$--{\it extendible
quasi-map}, denoted by \ $\Phi\colon \F^{\overmouse}  \leadsto \F^{\overmouse}$, \  if the
following conditions hold:   
\be
\item $\dom(\Phi)=\F^{\overmouse}$ 
\item $(\forall (\h,\f)\in \Phi)\,[d(\h)=d(\f)]$ 
\item $(\forall \h\in \F^{\overmouse})\,(\exists
\f\in \F^{\overmouse})\,[(\h,\f)\in\Phi]$  
\item $(\forall \h, \h^\prime\in \F^{\overmouse})\,(\forall \f,
\f^\prime\in  \F^{\overmouse})$ \[(\h,\f)\in \Phi\,\land\,
\h\equiv_{\mu_n}^{\overmouse}\h^\prime\,
\land
\,\f\equiv_{\mu_n}^{\overmouse}\f^\prime\Longrightarrow  (\h^\prime,\f^\prime)\in \Phi,\]
 where \ $n=d(\h)$
\item $(\forall (\h,\f)\in \Phi)\,(\forall (\h^\prime,\f^\prime)\in\Phi)\,(\forall
s\in{^{n}(n+k)}\!\uparrow)\,(\forall t\in{^{k}(n+k)}\!\uparrow)$  \[\overmouse \models \h
\,\Eq^{s,t}\,
\h^\prime \Longleftrightarrow \overmouse \models \f\,\Eq^{s,t}\,\f^\prime,\]
 where \
$n=d(\h)=d(\f)$ \ and \
$k=d(\h^\prime)=d(\f^\prime)$. 
\ee
In addition, given any \ $\h\in \F^{\overmouse}$, \ we shall write \ $\Phi(\h)$ \ to denote some function \ $\f\in
\F^{\overmouse}$ \ such that \ $(\h,\f)\in\Phi$.
\end{definition} 

\begin{rmk} If there exists a \ $\sigma\colon \F^{\overmouse}  \rightarrow \F^{\overmouse}$ \ which  is \ $\Eq$--extendible,
then one can easily define an \ $\Eq$--extendible quasi-map \ $\Phi\colon \F^{\overmouse}\leadsto \F^{\overmouse}$. \ The
converse needs some form of the axiom of choice. That is, given an \ $\Eq$--extendible quasi-map \ $\Phi$, \ one can ``thin'' 
\ $\Phi$ \ to a function \ $\sigma$ \ by {\it choosing} representatives from the appropriate equivalence classes. So 
our method of constructing the desired scales actually requires us to show that for \ $\tau\in\overmouse$, \ the existence of
an
\
$\Eq$--extendible quasi-map \ $\Phi\colon \F^{\overmouse}  \leadsto \F^{\overmouse}$, \ such that \ $\Phi\supset\tau$, \  is
definable over \ $\overmouse$.
\end{rmk}

The following theorem is in \cite[Theorem 3.26]{Cfsrm} and will be used in Part II as part of the proof of Theorem
\ref{newthm}. Given a structure \ $\calH$ \ we shall let \ $H$ \ denote the domain of this structure.
\begin{theorem}\label{defontau} Suppose \ $\mouse$ \ is a mouse, \ $\mathcal{H} \prec_1
\overmouse$, \ and \ $\R^{\overmouse}\subseteq H$. \ Let \ $\tau\subseteq \F^{\mathcal{H}}\times \F^{\mathcal{H}}$ \ be
nice. Then the condition on \ $\tau$, \ $(\exists \Phi\supseteq\tau)\left(\Phi\colon
\F^{\overmouse} 
\leadsto \F^{\overmouse}\right)$, \ is \ $\Sigma_\omega(\mathcal{H})$.
\end{theorem}

The proof of the above theorem requires that \ $\mathcal{H} \prec_1 \overmouse$. 
\ However, when \ $\mathcal{H} \prec_0\overmouse$ \ belongs to a certain good covering (see below) 
 of \ $\overmouse$, \ the conclusion of Theorem
\ref{defontau} holds, not necessarily for \ $\calH$, \  but for a slightly 
larger substructure \ $\calH^\prime$ \ of \ $\overmouse$, \ where \
$\calH\prec_0\calH^\prime\prec_0\overmouse$. 
\begin{definition}\label{goodcovdef} Let \ $\nouse =
(N,\in,\underline{\R},c_1,c_2,\dots,c_m,A_1,A_2,\dots,A_N)$ \ 
be a transitive model of \ $\Rplus$ \ containing all the reals, that is, \ $\ul{\R}^\nouse=\R$. \ 
Suppose that \ $\langle\mathcal{H}_i : i\in\omega\rangle$ \ is a sequence of substructures of \
$\nouse$ \ such that \ $\R\subseteq H_0\subseteq H_1\subseteq\cdots \subseteq\bigcup\limits_{i=0}^\omega H_i=N$.
\ We shall say that \ $\langle\mathcal{H}_i : i\in\omega\rangle$ \ is a {\it good covering\/} of \ $\nouse$ \ if
for each \
$i\in\omega$,
\be
\item $\mathcal{H}_i\prec_0\mathcal{H}_{i+1}$
\item $\mathcal{H}_i$ \ is in \ $\boldface{\Sigma}{\omega}(\mathcal{H}_{i+1})$
\item $\pow(\R)\cap\boldface{\Sigma}{\omega}(\mathcal{H}_i)\subseteq N$
\item there exists a \ $\boldface{\Sigma}{\omega}(\mathcal{H}_{e(i)})$ \ function \
$f\colon\R\maps{onto}H_i$
\ee
for some fixed \ $e\in\R$ \ where \ $e(j)\ge j$ \ for all \ $j\in\omega$.
\end{definition}

The next theorem is in \cite[Theorem 3.32]{Cfsrm} and will also be used in Part II to prove Theorem \ref{newthm}.
\begin{theorem}\label{defontautwo}  Let \ $\mouse$ \ be a real mouse satisfying \ $\AD$, \ and suppose that \ $n(\mouse)\ge
1$. \ Let \ $\langle\beta_i : i\in \omega\rangle$ \ be a sequence of ordinals such that \
$\kappa^{\overmouse}<\beta_i<\beta_{i+1}<\widehat{\OR}^{\overmouse}$ \ for all \ $i\in\omega$. \ Given \
$G\in\overmouse^{\beta_0}$, \ let \  $\langle\mathcal{H}_i : i\in\omega\rangle$ \ be defined by
\[\calH_i=\Hull^{\overmouse^{\beta_{i+1}}}_\omega(\{G,\beta_0,\dots,\beta_i\}\cup\R).\]
If \ $\overline{M}=\bigcup\limits_{i=0}^\omega H_i$, \ then \ 
\be
\item[(i)] $\langle\mathcal{H}_i : i\in\omega\rangle$ \ is a good covering of \ \overmouse, \ and
\item[(ii)] for \ $s\in\omega$ \ the condition on \ $\tau\in H_s$, \ $(\exists
\Phi\supseteq\tau)\left(\Phi\colon
\F^{\overmouse} \leadsto \F^{\overmouse}\right)$, \ is \ $\boldface{\Sigma}{\omega}(\mathcal{H}_{s+1})$.
\ee
\end{theorem}

This completes our overview, together with some new results, concerning the fundamental notions presented in \cite{Crcm} and
\cite{Cfsrm} which will be used here and in Parts II \& III.

\section{Iterable real premice are acceptable}\label{realaccept} In this section we shall prove the following
key theorem which will support our analysis of scales in \ $\Kr$.
\begin{theorem}\label{GenDJ} Suppose that \ $\mouse$ \ is an iterable real premouse.  Then \ $\mouse$ \ is
acceptable above the reals.
\end{theorem}

The above theorem is a generalization of Lemma 5.21 of Dodd-Jensen
\cite{DJ} and is essential for our work in Part III \cite{Part3}.  Theorem \ref{GenDJ} and its proof
\bi
\item justify our analysis of scales in \ $\Kr$ \ at the ``end of a gap''
(see \cite{Part3}), and
\item allow us to obtain scales of minimal complexity (e.g., see Lemma \ref{minicomp} below).
\ei

In Parts II and III we will be assuming the axiom of determinacy (AD). It is well-known that AD refutes the
axiom of choice. The Dodd-Jensen proof (of \cite[Lemma 5.21]{DJ}) uses the axiom of choice (see Remark
\ref{Justine}).   Our proof of Theorem \ref{GenDJ} is, in many respects, a generalization of the Dodd-Jensen
proof (see \cite[Corollary 11.27]{Dodd} and \cite[Lemma 5.21]{DJ}).  Nevertheless, we present our proof in the
interest of completeness\footnote{See Footnote \ref{justification}.} and, more importantly, to show explicitly that
our proof does not require the axiom of choice. In addition, a number of the lemmas used
in our proof will be applied in Part III.

First, we shall briefly indicate the basic strategy used in the proof of
Theorem \ref{GenDJ}. Let \ $\mouse$ \ be an iterable pure premouse above the reals and let \ $\kappa=\kappa^\mouse$. 
\ One proves that \ $\mouse^\gamma$ \ is acceptable by induction on \ $\gamma \le \widehat{\OR}^{\mouse}$. \ 
To prove that \ $\mouse^{\gamma+1}$ \ is acceptable for \ $\gamma<\widehat{\OR}^{\mouse}$, \ assuming that \
$\mouse^\gamma$ \ is acceptable, one must prove  Lemmas \ref{lemma1}--\ref{lemmatwo} and Lemmas
\ref{lemma3}--\ref{lemma6} below.  We shall now begin to prove these technical lemmas.

\begin{lemma}\label{lemma1} Let \ $\mouse=(M,\R,\kappa,\mu)$ \ be  an iterable real premouse. Let \ $\gamma <
\widehat{\OR}^{\mouse}$. \ Assume that \
${\mouse^\gamma}$ \ is acceptable and critical. Let
\ $n=n(\mouse^\gamma)$. \ If \ $\rho_{\mouse^\gamma}^{n+1} < \kappa$, \ then \ $\pow(M^\gamma)\cap
M^{\gamma+1}\subseteq \boldface{\Sigma}{\omega}(\mouse^{\gamma})$.
\end{lemma}

\begin{proof} Theorem \ref{intismi} implies that \ ${\mouse^\gamma}$ \ is a mouse. Now (5) of Lemma \ref{defone}
implies the desired conclusion. 
\end{proof}

The next lemma is a generalization of Lemma 4.9 of Dodd-Jensen \cite{DJ}. Our proof is modeled in part
after the  proofs presented in \cite[see pp. 88--91]{Dodd} and \cite[pp. 62--63]{DJ}. 
\begin{lemma}\label{lemmatwo} Let \ $\mouse=(M,\R,\kappa,\mu)$ \ be  an iterable real premouse. Let \ $\gamma <
\widehat{\OR}^{\mouse}$. \ Assume that \
${\mouse^\gamma}$ \ is acceptable and critical. Let \ $n=n(\mouse^\gamma)$. \
If \ $\rho_{\mouse^\gamma}^{n+1} = \kappa$, \ then \
$H^{\mouse^\gamma}_{\kappa}=H^{\mouse^{\gamma+1}}_{\kappa}$.
\end{lemma}

\begin{proof} We will use the premouse iteration of \ $\mouse^{\gamma+1}$. \ To help simplify the notation
in this proof we will violate our notational convention as follows: We shall write
\[\premouseiterationt{{\mouse^{\gamma_\alpha+1}}}{\alpha}{\beta}{\mouse^{\gamma_\beta+1}}\]
for the premouse iteration of \ $\mouse^{\gamma+1}$ \ where we let \
$\mouse^{\gamma_\alpha+1}=(M^{\gamma_\alpha+1},\R,\kappa_\alpha,\mu_\alpha)$ \ denote the \ $\alpha^{\ul{\text{th}}}$ \
iterate of \ $\mouse^{\gamma+1}=(M^{\gamma+1},\R,\kappa,\mu)$, \ for each ordinal \ $\alpha$; \ that is, 
\ $\mouse^{\gamma_\alpha+1}=\left(\mouse^{\gamma+1}\right)_\alpha$ \ (this is our main notational violation).  
For each such ordinal \ $\alpha$, \ we shall also write \ 
$\mouse^{\gamma_\alpha}=(M^{\gamma_\alpha},\R,\kappa_\alpha,\mu_\alpha)$ \ for the ``predecessor'' of \
$\mouse^{\gamma_\alpha+1}$, \ that is, \
$M^{\gamma_\alpha+1}$ \ is the ``rudimentary in $\mu_\alpha$'' closure of \
$M^{\gamma_\alpha}$.
\ In addition, \ let \ $\ov{\mouse^{\gamma_\alpha}}=\left({\mouse^{\gamma_\alpha}}\right)^{n+1}$ \
and let  \ $\ov{M^{\gamma_\alpha}}$ \ denote the
domain of this structure. We are using the notation \
$\ov{\mouse^{\gamma_\alpha}}=\left({\mouse^{\gamma_\alpha}}\right)^{n+1}$ \
to distinguish from \ $\overline{\mouse^{\gamma_\alpha}}=\left({\mouse^{\gamma_\alpha}}\right)^{n}$, \ where \
$n=n(\mouse^{\gamma_\alpha})$.

Note that (recall Definitions \ref{master}, \ref{master2}, \ref{master3}, \ref{master4}) for all ordinals \
$\alpha\le\beta$ \ we have the following:
\newcounter{plist_recall}
\newcounter{plist_recall2}
\begin{list}
{\upshape (\arabic{plist})}
{\setlength{\itemsep}{3pt}
\setlength{\labelwidth}{.5in}
\usecounter{plist}}
\item $\mouse^{\gamma_\alpha}\in M^{\gamma_\alpha+1}$
\item $\pi_{\alpha\beta}\left(\mouse^{\gamma_\alpha}\right)=\mouse^{\gamma_\beta}$
\item $\pi_{\alpha\beta}\colon\mouse^{\gamma_\alpha}\mapsigma{\omega}\mouse^{\gamma_\beta}$
and so,
\item $\mouse^{\gamma_\alpha}$ \ and \ $\mouse^{\gamma_\beta}$ \ are acceptable and critical
\item $\pi_{\alpha\beta}\left(\rho_{\mouse^{\gamma_\alpha}}^k\right)=\rho_{\mouse^{\gamma_\beta}}^k$, \ 
$\pi_{\alpha\beta}\left(p_{\mouse^{\gamma_\alpha}}^k\right)=p_{\mouse^{\gamma_\beta}}^k$, \ 
$\pi_{\alpha\beta}\left(\mathcal{A}_{\mouse^{\gamma_\alpha}}^{k}\right)=\mathcal{A}_{\mouse^{\gamma_\beta}}^{k}$ \ for all
$k\le n=n(\mouse^{\gamma})$
\item \setcounter{plist_recall2}{\value{plist}} $\rho_{\mouse^{\gamma_\alpha}}^{n+1} = \kappa_\alpha$ \  and hence, \
$n(\mouse^{\gamma_\alpha})=n=n(\mouse^{\gamma})$
\item $\ov{\mouse^{\gamma_\alpha}}\in M^{\gamma_\alpha+1}$
\item $\pi_{\alpha\beta}\left(\ov{\mouse^{\gamma_\alpha}}\right)=\ov{\mouse^{\gamma_\beta}}$
\item \setcounter{plist_recall}{\value{plist}} 
$\pi_{\alpha\beta}\colon\ov{\mouse^{\gamma_\alpha}}\mapsigma{\omega} \ov{\mouse^{\gamma_\beta}}$.
\end{list}
\renewcommand{\theclam}{\arabic{clam}}
\begin{clam}\label{claimzero} For ordinals \ $\alpha\le\beta$, \ 
\be
\item[(i)] $\pow(\R)\cap M^{\gamma_\alpha}=\pow(\R)\cap M^{\gamma_\beta}$,
\item[(ii)] $\pow(\R)\cap \ov{M^{\gamma_\alpha}}=\pow(\R)\cap \ov{M^{\gamma_\beta}}$.
\ee
\end{clam}

\begin{proof}[Proof of Claim \ref{claimzero}] For (i) we note that \ $\pow(\R)\cap M^{\gamma_\alpha+1}=\pow(\R)\cap
M^{\gamma_\beta+1}$ \ because the measures are $\R$--complete. Because \
$\pi_{\alpha\beta}\colon\mouse^{\gamma_\alpha+1}\mapsigma{1}\mouse^{\gamma_\beta+1}$, \ it follows that 
$\pow(\R)\cap M^{\gamma_\alpha}=\pow(\R)\cap M^{\gamma_\beta}$. \ To prove (ii), observe that since \
$\mouse^{\gamma_\alpha}$
\ and  \ $\mouse^{\gamma_\beta}$ \ are acceptable and critical, it follows that \ $\pow(\R)\cap
\ov{M^{\gamma_\alpha}}=\pow(\R)\cap
\ov{M^{\gamma_\beta}}$ \ by Corollaries 1.34 and 2.13 of \cite{Cfsrm}.
\end{proof}

\begin{clam}\label{claimone} For ordinals \ $\alpha\le\beta$, \ 
$\ov{\mouse^{\gamma_\alpha}} \prec_\omega \ov{\mouse^{\gamma_\beta}}$.
\end{clam}

\begin{proof}[Proof of Claim \ref{claimone}] Let \ $\alpha\le\beta$ \ be ordinals. Let \
$\ov{M^{\gamma_\alpha}}$ \ be the domain of the structure \ $\ov{\mouse^{\gamma_\alpha}}$. \
Since \
$\rho_{\mouse^{\gamma_\alpha}}^{n+1} =
\kappa_\alpha$, \  Lemma \ref{kappa_small} implies that \
$\ov{M^{\gamma_\alpha}}=H^{\mouse^{\gamma_\alpha}}_{\kappa_\alpha}$. \ Since every set in \
$H^{\mouse^{\gamma_\alpha}}_{\kappa_\alpha}$ \ has \ $\mouse^{\gamma_\alpha}$--cardinality less than \ $\kappa_\alpha$, \
it follows that \ $\pi_{\alpha\beta}(a)=a$ \ for all \ $a\in \ov{M^{\gamma_\alpha}}$. \ Therefore, 
the above (\arabic{plist_recall}) implies that \ $\ov{\mouse^{\gamma_\alpha}} \prec_\omega
\ov{\mouse^{\gamma_\beta}}$. \ This completes the proof of Claim 2. 
\end{proof}

\begin{clam}\label{claimtwo} For ordinals \ $\alpha<\beta$, \ 
$\ov{\mouse^{\gamma_\alpha}}\in \ov{M^{\gamma_\beta}}$, \ where \
$\ov{M^{\gamma_\beta}}$ \ is the domain of the structure \ $\ov{\mouse^{\gamma_\beta}}$.
\end{clam}

\begin{proof}[Proof of Claim \ref{claimtwo}] Let \ $\alpha<\beta$ \ be any ordinals. Definitions \ref{master},
\ref{master2}, \ref{master3}, and \ref{master4} imply that 
\[\ov{\mouse^{\gamma_\alpha}}=
\left(\ov{M^{\gamma_\alpha}}, \in, \R, \mathcal{A}_{\mouse^{\gamma_\alpha}}^{1}\cap
\ov{M^{\gamma_\alpha}},
\mathcal{A}_{\mouse^{\gamma_\alpha}}^{2}\cap\ov{M^{\gamma_\alpha}},
\dots,\mathcal{A}_{\mouse^{\gamma_\alpha}}^{n}\cap\ov{M^{\gamma_\alpha}}\right)\] 
and that \ 
$\ov{M^{\gamma_\alpha}}=J_{\kappa_\alpha}^{\ov{\mouse^{\gamma_\alpha}}}(\underline{\R})$, \
 the Jensen hierarchy of sets which are relatively constructible above \ $\R$ \ from
the predicates \ $\mathcal{A}_{\mouse^{\gamma_\alpha}}^{1}\cap \ov{M^{\gamma_\alpha}},
\dots,\mathcal{A}_{\mouse^{\gamma_\alpha}}^{n}\cap\ov{M^{\gamma_\alpha}}$. \ Note that \
$\mathcal{A}_{\mouse^{\gamma_\alpha}}^{i}\cap\ov{M^{\gamma_\alpha}}\subseteq
\R\times(\kappa_\alpha)^{<\omega}$  \ for all \ $i\le n$.
\ Now, since \ $\pi_{\alpha\beta}(a)=a$ \ for all \ $a\in \R\times(\kappa_\alpha)^{<\omega}$, \ it follows that 
\[\mathcal{A}_{\mouse^{\gamma_\alpha}}^{i}\cap\R\times(\kappa_\alpha)^{<\omega} = 
\mathcal{A}_{\mouse^{\gamma_\beta}}^{i}\cap\R\times(\kappa_\alpha)^{<\omega}\]
for all \ $i\le n$. \ Because \ $\kappa_\alpha<\kappa_\beta$, \ we conclude that \
$\ov{\mouse^{\gamma_\alpha}}\in
\ov{M^{\gamma_\beta}}$.
\end{proof}

For each ordinal \ $\alpha$, \ let \ $A_\alpha\subseteq\R\times(\kappa_\alpha)^{<\omega}$ \ be some canonical coding of the
predicates \[\mathcal{A}_{\mouse^{\gamma_\alpha}}^{1}\cap
\ov{M^{\gamma_\alpha}},\, \mathcal{A}_{\mouse^{\gamma_\alpha}}^{2}\cap\ov{M^{\gamma_\alpha}},\,
\dots,\,\mathcal{A}_{\mouse^{\gamma_\alpha}}^{n}\cap\ov{M^{\gamma_\alpha}}\] so that  \
$\ov{M^{\gamma_\alpha}}=L_{\kappa_\alpha}(\R)[A_\alpha]$. \ Note that \
$\R\in\ov{M^{\gamma_\alpha}}$ \ and the predicate \ $A_\alpha$ \ can be thought of as a class in the
structure \ $\ov{\mouse^{\gamma_\alpha}}$.

\begin{clam}\label{claimthree} For ordinals \ $\alpha$, \ we have that \ $\ov{\mouse^{\gamma_\alpha}}\models
\ZF + V\eq L(\R)[A_\alpha]$, \ and \
$\Theta^{\ov{\mouse^{\gamma_\alpha}}}=\Theta^{\ov{\mouse^{\gamma}}}<\kappa$.
\end{clam}

\begin{proof}[Proof of Claim  \ref{claimthree}] Let \ $\alpha$ \ be any ordinal. Clearly, \
$\ov{\mouse^{\gamma_\alpha}} \models V\eq L(\R)[A_\alpha]$. \ We now show that \
$\ov{\mouse^{\gamma_\alpha}} \models
\ZF$. \ Since \ $\ov{\mouse^{\gamma_\alpha}}\models\Rplus$, \ we only need to check that it satisfies the
separation, collection and power set axioms. To see the \ $\ov{\mouse^{\gamma_\alpha}}$ \ satisfies
separation and collection, let \
$\beta>\alpha$ \ be any ordinal. By \ Claim \ref{claimone} and Claim \ref{claimtwo}, we have that 
\newcounter{puff}
\begin{list}{{\upshape (\roman{puff})}}
{\setlength{\labelwidth}{20pt}
\usecounter{puff}}
\item $\ov{\mouse^{\gamma_\alpha}} \prec_\omega \ov{\mouse^{\gamma_\beta}}$,
\item $\ov{\mouse^{\gamma_\alpha}}\in \ov{M^{\gamma_\beta}}$ \ and
$\ov{M^{\gamma_\alpha}}\in \ov{M^{\gamma_\beta}}$.
\end{list}
The conditions (i) and (ii) easily imply that \ $\ov{\mouse^{\gamma_\alpha}}$ \ satisfies separation and
collection.

Finally, we prove that \ $\ov{\mouse^{\gamma_\alpha}}$ \ satisfies the power set axiom.\footnote{
As noted previously, our proof is modeled after the proof presented in \cite[see p. 89]{Dodd}.
Dodd's proof, however, uses the axiom of choice. Our proof explicitly shows that the axiom of choice
is not necessary.} Consider the inner model \ $N=L(\mouse^{\gamma_\alpha+1})$. \ Let \
$\Q\in N$  
\ be the partial order in Definition \ref{col_to_omega}. Now let \ $G$ \ be \ $\Q$--generic over \ $N$. \ Thus, \
$N[G]\models\ZFC$ \ and, by absoluteness, everything that we have proven about \ $\mouse^{\gamma_\alpha+1}$ \ holds in \
$N[G]$; \ in particular, \ $\ov{\mouse^{\gamma_\alpha}}\in N[G]$. \ Now we shall show that \
\[N[G]\models\text{``$\ov{\mouse^{\gamma_\alpha}}$ satisfies the power set axiom''}\] and thus, again by
absoluteness, \
$\ov{\mouse^{\gamma_\alpha}}\models \text{``the power set axiom''}$.
\ We will now work in \ $N[G]$ \ (until the end of this paragraph).  Let \ $a\in\ov{M^{\gamma_\alpha}}$. \
Let \ $\theta$ \ be a regular cardinal such that \ $\theta>\abs{\pow(a)}$. \ Then \ $\kappa_\theta=\theta$ \ and
therefore, \
$\abs{\pow(a)\cap\ov{M^{\gamma_\theta}}}<\theta$. \ Because \
$\ov{M^{\gamma_\theta}}=J_{\kappa_\theta}^{\ov{\mouse^{\gamma_\theta}}}(\underline{\R})
=J_{\theta}^{\ov{\mouse^{\gamma_\theta}}}(\underline{\R})$
\ (the Jensen hierarchy of sets relatively constructible above \ $\R^N$ \ in
the predicates \ $\mathcal{A}_{\mouse^{\gamma_\theta}}^{1}\cap \ov{M^{\gamma_\theta}},
\dots,\mathcal{A}_{\mouse^{\gamma_\theta}}^{n}\cap\ov{M^{\gamma_\theta}}$), \ we conclude that \ 
$\pow(a)\cap\ov{M^{\gamma_\theta}}\subseteq
J_{\gamma}^{\ov{\mouse^{\gamma_\theta}}}(\underline{\R})$ \ for some \ $\gamma<\theta=\kappa_\theta$. \
Thus, \ $\ov{M^{\gamma_\theta}}\models\text{``$\pow(a)$ exists''}$. \ Since 
\ $\pi_{\alpha\theta}(a)=a$, \ we see that \ $\ov{M^{\gamma_\alpha}}\models\text{``$\pow(a)$ exists''}$. 

Therefore, 
\ $\ov{\mouse^{\gamma_\alpha}} \models \ZF$. \ Claim
\ref{claimzero} implies that
\ $\Theta^{\ov{\mouse^{\gamma_\alpha}}}=\Theta^{\ov{\mouse^{\gamma}}}<\kappa$.
\end{proof}

\begin{clam}\label{claimfour} For all ordinals \ $\alpha<\beta$, \ 
$\ov{M^{\gamma_\alpha}}=H_{\kappa_\alpha}^{\ov{\mouse^{\gamma_\beta}}}$.\footnote{Recall
Definition
\ref{trans}.} \ Thus, \ $\ov{M^{\gamma_\alpha}}$ \ is uniformly definable over \
$\ov{\mouse^{\gamma_\beta}}$ \ from
\ $\kappa_\alpha$.
\end{clam}

\begin{proof}[Proof of Claim \ref{claimfour}] Let \ $\alpha<\beta$. \ Since \
$\ov{M^{\gamma_\alpha}}=H^{\mouse^{\gamma_\alpha}}_{\kappa_\alpha}$, \ it follows that \
$\ov{M^{\gamma_\alpha}}\subseteq H_{\kappa_\alpha}^{\ov{\mouse^{\gamma_\beta}}}$ \ because for
each \
$a\in H^{\mouse^{\gamma_\alpha}}_{\kappa_\alpha}$ \ the elements that witness this fact are in \
$\ov{M^{\gamma_\alpha}}\subseteq \ov{M^{\gamma_\beta}}$. \ To show that \
$H_{\kappa_\alpha}^{\ov{\mouse^{\gamma_\beta}}}\subseteq \ov{M^{\gamma_\alpha}}$, \ let \ $a\in
H_{\kappa_\alpha}^{\ov{\mouse^{\gamma_\beta}}}$. \ Because \ 
$\ov{\mouse^{\gamma_\alpha}}\prec_\omega
\ov{\mouse^{\gamma_\beta}}$ \ (by Claim \ref{claimone}), \ one concludes that \
$\ov{\mouse^{\gamma_\beta}}\models
\text{``$\kappa_\alpha$ is a limit cardinal''}$\footnote{Here, cardinality is to be interpreted as in Definition
\ref{R-card}.}. \ Therefore, \ $\ov{\mouse^{\gamma_\beta}}\models a\in H_\gamma$ \ for some \
$\gamma<\kappa_\alpha$. \ Since \ $\ov{\mouse^{\gamma_\alpha}}\prec_\omega
\ov{\mouse^{\gamma_\beta}}$,
\ we see that \
$H_{\gamma}^{\ov{\mouse^{\gamma_\beta}}}=H_{\gamma}^{\ov{\mouse^{\gamma_\beta}}}$. \ Thus, \
$a\in\ov{M^{\gamma_\alpha}}$.
\end{proof}

\begin{clam}\label{uniformly} For all ordinals \ $\alpha<\beta$, \ $\mouse^{\gamma_\alpha}\in
\ov{M^{\gamma_\beta}}$ \ and
\ $\mouse^{\gamma_\alpha}$ \ is uniformly definable over \ $\ov{\mouse^{\gamma_\beta}}$ \ from \
$\kappa_\alpha$.
\end{clam}

\begin{proof}[Proof of Claim \ref{uniformly}] Let \ $\alpha<\beta$. \ Claim \ref{claimtwo} (with its proof) and Claim
\ref{claimthree}, together with Theorem 2.11\footnote{This is the model extension theorem (see Theorem \ref{MET}
of this paper).} and Corollary 2.12 of \cite{Cfsrm}, imply the desired conclusions with a minor
extension of \cite[Theorem 2.11]{Cfsrm}. Recall that \
$\overline{\mouse^{\gamma_\alpha}}=\left({\mouse^{\gamma_\alpha}}\right)^{n}$ \ and \
$\ov{\mouse^{\gamma_\alpha}}=\left({\mouse^{\gamma_\alpha}}\right)^{n+1}$ \ where \
$n=n(\mouse^{\gamma_\alpha})$. \ We note that \ $\overline{\mouse^{\gamma_\alpha}}$ \ is a
model of \ $\mathcal{T}^{n}$, \ whereas \ $\ov{\mouse^{\gamma_\beta}}$
\ is not a model of \ $\mathcal{T}^{n+1}$ \ (since \ $\ov{\mouse^{\gamma_\beta}}$ \ ``is not above
$\kappa$'', \ it is not a premouse). However, the measure of \ $\overline{\mouse^{\gamma_\beta}}$ \ is coded by the
$\Sigma_1$--master code of \ $\overline{\mouse^{\gamma_\beta}}$. \ So we can define a set of \ $\Pi_2$ \ axioms \ 
${\wht{\mathcal{T}}}^{n+1}$ \ which is like \ $\mathcal{T}^{n+1}$ \ except that the axioms in \
${\wht{\mathcal{T}}}^{n+1}$ \ describe the measure \ $\mu$ \ ``in the master code.'' By slightly generalizing
Theorem 2.11 of \cite{Cfsrm} (and its proof) to include the ``single step'' for building  \
$\overline{\mouse^{\gamma_\alpha}}$ \ from \ $\ov{\mouse^{\gamma_\alpha}}$, \ we see that \
$\mouse^{\gamma_\alpha}$ \ can be constructed from \ $\ov{M^{\gamma_\alpha}}$. \ Now since \
$\ov{M^{\gamma_\alpha}}\in \ov{M^{\gamma_\beta}}$, \ this uniform construction can be carried
out in \ $\ov{M^{\gamma_\beta}}$ \ because it is a model of \ $\ZF$. \ Thus, the conclusion of the claim holds.
\end{proof}

\begin{clam}\label{claimfive} Let \ $\alpha\le\beta$ \ be ordinals. Then \  $\{\,\kappa_\gamma : \alpha\le\gamma < \beta\,\}$ \
is a set of order \ $\Sigma_\omega(\ov{\mouse^{\gamma_\beta}},\R\times\kappa_{\alpha})$ \ indiscernibles.
\end{clam}

\begin{proof}[Proof of Claim \ref{claimfive}] Let \ $\alpha\le\beta$. \ Then  \  $\{\,\kappa_\gamma : \alpha\le\gamma <
\beta\,\}$ \ is a set of order
\ $\Sigma_1(\mouse^{\gamma_\beta+1},\{\,\pi_{\alpha\beta}(a) : a\in M^{\gamma_{\alpha}+1}\,\})$ \ indiscernibles by
Corollary 2.14(1) of \cite{Crcm}. Since \ $\ov{\mouse^{\gamma_{\alpha}}}\in \mouse^{\gamma_{\alpha}+1}$ \
and \
$\ov{\mouse^{\gamma_\beta}}=\pi_{\alpha\beta}(\ov{\mouse^{\gamma_{\alpha}}})$, \ the conclusion
of the claim holds.
\end{proof}

\begin{clam}\label{claim5plus} Let \ $\alpha\le\beta$ \ be ordinals. \ Let \
$X_\alpha\subseteq\R\times\kappa_\alpha$ \ be 
\[\Sigma_\omega(\ov{\mouse^{\gamma_{(\alpha+i)}}},\R\cup\{\nu,\kappa_\alpha,\dots,\kappa_{\alpha+(i-1)}\}),\]
where \ $\nu<\kappa_\alpha$ \ and \ $i\in\omega$. \ Let \ $X_\beta$ \ have the same 
$\Sigma_\omega(\ov{\mouse^{\gamma_{(\beta+i)}}},\R\cup\{\nu,\kappa_\beta,\dots,\kappa_{\beta+(i-1)}\})$ \
definition  (using the same $\nu$). \ Then \ $X_\alpha\in \mouse^{\gamma_{\alpha}+1}$, \ $X_\beta\in
\mouse^{\gamma_{\beta}+1}$ \ and \ $\pi_{\alpha\beta}(X_\alpha)=X_\beta$.
\end{clam}

\begin{proof}[Proof of Claim \ref{claim5plus}] Let \ $\alpha\le\beta$, \ $i\in\omega$ \ and let \
$X_\alpha\subseteq\R\times\kappa_\alpha$ \ be 
\[\Sigma_\omega(\ov{\mouse^{\gamma_{(\alpha+i)}}},\R\cup\{\nu,\kappa_\alpha,\dots,\kappa_{\alpha+(i-1)}\}).\] 
Hence, \ $X_\alpha$ \ is \
$\Sigma_0(\mouse^{\gamma_{(\alpha+i)}+1},\R\cup\{\ov{\mouse^{\gamma_{(\alpha+i)}}},
\nu,\kappa_\alpha,\dots,\kappa_{\alpha+(i-1)}\})$. \ Since 
\[\pi_{\alpha,\alpha+i}(\ov{\mouse^{\gamma_{\alpha}}})=\ov{\mouse^{\gamma_{(\alpha+i)}}},\] and \
$\pi_{\alpha,\alpha+i}$ \ fixes all the elements in \ $\R\cup\{\R,\nu\}$, \ the proof of Theorem 2.13 of \cite{Crcm}
shows that there is a rudimentary function \ $G$ \ such that for all \ $b\in\R\times\kappa_\alpha$, \
\[b\in X_\alpha \iff \mouse^{\gamma_{\alpha}+1} \models G(b)\in\mu_i\]
where \  $\mu_i$ \ is uniformly ``rudimentary over $\mouse^{\gamma_{\alpha}+1}$'', \ that is, there is a
rudimentary function \ $\mathcal{F}$ \ such that 
\[a\in\mu_i \iff \nouse \models \mathcal{F}(a)\ne\emptyset\]
for all \ $a\in \nouse$, \ 
where \ $\nouse$ \ is any real premouse (see subsection \ref{realpremice} for the definition of $\mu_i$). 

We note that the
definition of \ $G$ \ depends only on the {\it definition\/} of \ $X_\alpha$ \ and the {\it parameters\/} from \
$\R\cup\{\R,\nu\}$ \ used in the definition. \ Consequently,  for all \ $b\in\R\times\kappa_\beta$, \
\[b\in X_\beta \iff \mouse^{\gamma_{\beta}+1} \models G(b)\in\mu_i.\]
Thus, \ $X_\alpha\in \mouse^{\gamma_{\alpha}+1}$, \ $X_\beta\in \mouse^{\gamma_{\beta}+1}$. \
Since \ $\pi_{\alpha\beta}\colon\mouse^{\gamma_{\alpha}+1}\mapsigma{1}\mouse^{\gamma_{\beta}+1}$, \ it follows that \ 
$\pi_{\alpha\beta}(X_\alpha)=X_\beta$.
\end{proof}

Let \ $\alpha$ \ be an ordinal and \ $k\in\omega$ \ where \ $k> 0$. \ Let \ $H_{\alpha}^k =
H_{\kappa_\alpha^+}^{\ov{\mouse^{\gamma_{(\alpha+k)}}}}$, \ where \ $\kappa_\alpha^+$ \ denotes the
successor of \ $\kappa_\alpha$ \ in \ $\ov{\mouse^{\gamma_{(\alpha+k)}}}$ \ (applying Definition \ref{realcard}). 
One can check that \ $H_{\alpha}^1=H_{\alpha}^m$ \ for all integers \ $m>0$ \ (see the proof of Claim \ref{claimfour}). 
Now, define
\[K^0_\alpha =\left\{a\in H_{\alpha}^1 : \{a\} \in
\Sigma_\omega\!\left(\ov{\mouse^{\gamma_{\alpha}}},
\R\cup\kappa_\alpha\right)\right\}\]
and for integers \ $m> 0$ \ define
\[K^m_\alpha = \left\{a\in H_{\alpha}^1 : \{a\} \in
\Sigma_\omega\!\left(\ov{\mouse^{\gamma_{(\alpha+m)}}},
\R\cup\kappa_\alpha\cup\{\kappa_\alpha,\dots\kappa_{\alpha+(m-1)}\}\right)\right\}.\]

\begin{clam}\label{claimsix} Let \ $\alpha$ \ be an ordinal and \ $m\in\omega$. \ Then
\newcounter{cntone}
\newcounter{cnttwo}
\newcounter{cntthree}
\newcounter{cntfour}
\newcounter{cntfive}
\newcounter{cntsix}
\begin{list}{{\upshape (\alph{puff})}}
{\setlength{\labelwidth}{20pt}
\usecounter{puff}}
\item \setcounter{cntone}{\value{puff}} $K^0_\alpha=\ov{M^{\gamma_{\alpha}}}$,
\item \setcounter{cnttwo}{\value{puff}}$K^m_\alpha\in H_{\alpha}^1$,
\item \setcounter{cntthree}{\value{puff}}$K^m_\alpha\in K^{m+1}_\alpha$,
\item \setcounter{cntfour}{\value{puff}}$\mouse^{\gamma_\alpha}\in K^1_\alpha$,
\item \setcounter{cntfive}{\value{puff}}$K^m_\alpha$ \ is transitive,
\item \setcounter{cntsix}{\value{puff}}$K^m_\alpha$ \ is closed under the ``rudimentary'' functions.
\end{list}
\end{clam}

\begin{proof}[Proof of Claim \ref{claimsix}] We prove (\alph{cntone})--(\alph{cntsix}) as follows:

\medskip
(\alph{cntone}) We need to show that 
\[\ov{M^{\gamma_{\alpha}}}= \left\{a\in H_{\alpha}^1 : \{a\} \in
\Sigma_\omega\!\left(\ov{\mouse^{\gamma_{\alpha}}},
\R\cup\kappa_\alpha\right)\right\}.\]
By Claim \ref{claimfour}, \
$\ov{M^{\gamma_\alpha}}=H_{\kappa_\alpha}^{\ov{\mouse^{\gamma_{(\alpha+1)}}}}$.
\  Thus, \ $\ov{M^{\gamma_\alpha}}\subseteq H_{\alpha}^1$. \ We note that \
$\OR^{\ov{\mouse^{\gamma_\alpha}}}=\kappa_\alpha$. \ Since there is a
\ $\Sigma_1(\ov{\mouse^{\gamma_\alpha}})$ \ function \ $f$ \ such that  
$f\colon \OR^{\ov{\mouse^{\gamma_\alpha}}} \times {\R} \maps{onto} \ov{M^{\gamma_\alpha}}$ \ (see
\cite[Corollary 1.8]{Crcm} and \cite[Lemma 1.4]{Cfsrm}), we conclude that \
$\ov{M^{\gamma_{\alpha}}}=K^0_\alpha$.

\medskip
(\alph{cnttwo}) We shall show that \ $K^m_\alpha\in H_{\alpha}^1$, \ where
\[K^m_\alpha = \left\{a\in H_{\alpha}^1 : \{a\} \in
\Sigma_\omega\!\left(\ov{\mouse^{\gamma_{(\alpha+m)}}},
\R\cup\kappa_\alpha\cup\{\kappa_\alpha,\dots\kappa_{\alpha+(m-1)}\}\right)\right\}.\]
Let \ $k>m$. \ As noted before, \ $H_{\alpha}^1=H_{\kappa_\alpha^+}^{\ov{\mouse^{\gamma_{(\alpha+m)}}}}
=H_{\kappa_\alpha^+}^{\ov{\mouse^{\gamma_{(\alpha+k)}}}}$. \ By Claim \ref{claimtwo} \
$\ov{\mouse^{\gamma_{(\alpha+m)}}}\in \ov{M^{\gamma_{(\alpha+k)}}}$. \ 
Because \ $\ov{\mouse^{\gamma_{(\alpha+k)}}}\models \ZF$ \ by Claim \ref{claimthree}, it follows that 
\ $K^m_\alpha\in \ov{M^{\gamma_{(\alpha+k)}}}$. \ In addition, 
\ $\ov{M^{\gamma_{(\alpha+k)}}}\models \abs{K^m_\alpha}\le\kappa_\alpha$\footnote{Here, the cardinal
$\abs{K^m_\alpha}$ is to be interpreted as in Definition
\ref{R-card}.} and from the point of view of \ $\ov{M^{\gamma_{(\alpha+k)}}}$ \ 
\[K^m_\alpha= \left\{a\in H_{\kappa_\alpha^+} : \{a\} \in
\Sigma_\omega\!\left(\ov{\mouse^{\gamma_{(\alpha+m)}}},
\R\cup\kappa_\alpha\cup\{\kappa_\alpha,\dots\kappa_{\alpha+(m-1)}\}\right)\right\}.\]
Since \ $\ov{M^{\gamma_{(\alpha+k)}}}$  \ is a model of \ $\ZF +
V\eq L(\R)[A_{\alpha+k}]$, \ Theorem \ref{card} implies that \
$\ov{M^{\gamma_{(\alpha+k)}}}\models K^m_\alpha\in H_{\kappa_\alpha^+}$ \ (see \cite[Lemma 6.4, pp.
131-132]{Kunen}).  Hence, \ $K^m_\alpha\in H_{\alpha}^1$.

\medskip
(\alph{cntthree}) We will prove that \ $K^m_\alpha\in K^{m+1}_\alpha$. \ By \ (\alph{cnttwo}), \ $K^m_\alpha\in H_{\alpha}^1$.
\ So we must show that \ $\{K^m_\alpha\}\in \Sigma_\omega\!\left(\ov{\mouse^{\gamma_{(\alpha+k)}}},
\R\cup\kappa_\alpha\cup\{\kappa_\alpha,\dots\kappa_{(\alpha+m)}\}\right)$, \ where \ $k=m+1$. \ By Claim \ref{uniformly}, \
$\mouse^{\gamma_{\alpha+m}}\in \ov{M^{\gamma_{(\alpha+k)}}}$ \ and is definable in \
$\ov{\mouse^{\gamma_{(\alpha+k)}}}$ \ from \ $\kappa_{\alpha+m}$.
\ Also, \ $H_\alpha^1$ \ is definable in \ $\ov{\mouse^{\gamma_{(\alpha+k)}}}$ \ from \ $\kappa_\alpha$. \ It
follows,  again because \ $\ov{\mouse^{\gamma_{(\alpha+k)}}}\models\ZF$, \ that \ 
$\{K^m_\alpha\}\in \Sigma_\omega\!\left(\ov{\mouse^{\gamma_{(\alpha+k)}}},
\R\cup\kappa_\alpha\cup\{\kappa_\alpha,\dots\kappa_{\alpha+m}\}\right)$. \ Thus, \ $K^m_\alpha\in K^{m+1}_\alpha$.

\medskip
(\alph{cntfour}) We must show that  \[\mouse^{\gamma_\alpha}\in
K^{1}_\alpha = \left\{a\in H_{\alpha}^1 :
\{a\} \in \Sigma_\omega\!\left(\ov{\mouse^{\gamma_{(\alpha+1)}}},
\R\cup\kappa_\alpha\cup\{\kappa_\alpha\}\right)\right\}.\]
By Claim \ref{uniformly}, \ $\mouse^{\gamma_\alpha}\in \ov{M^{\gamma_{(\alpha+1)}}}$ \ and is definable over
\ $\ov{\mouse^{\gamma_{(\alpha+1)}}}$ \ from \ $\kappa_\alpha$. \ As noted in (\arabic{plist_recall2}) (just
above Claim \ref{claimzero}), \ $\rho_{\mouse^{\gamma_\alpha}}^{n+1} = \kappa_\alpha$. \ So Corollaries 1.32 and 2.13 of
\cite{Cfsrm} imply that there is a \ $\boldface{\Sigma}{n+1}(\mouse^{\gamma_\alpha})$ \ function \ $f$ \ such that  \
$f\colon \kappa_\alpha \times {\R} \maps{onto} M^{\gamma_\alpha}$. \ This implies that \ $\mouse^{\gamma_\alpha}\in
H_{\alpha}^1$. \ Therefore, \ $\mouse^{\gamma_\alpha}\in K^{1}_\alpha$.

\medskip
(\alph{cntfive}) By (\alph{cntone}) we see that \ $K^0_\alpha$ \ is transitive. So we now consider the case \ $m>0$ \ and prove
that \ $K^m_\alpha$ \ is transitive. Let \ $b\in a\in K^m_\alpha$. \ We must show that \ $b\in K^m_\alpha$. \ 
Since \ $a\in K^m_\alpha$, \ $\{a\} \in \Sigma_\omega\!\left(\ov{\mouse^{\gamma_{(\alpha+m)}}},
\R\cup\kappa_\alpha\cup\{\kappa_\alpha,\dots\kappa_{\alpha+(m-1)}\}\right)$. \ In addition, \
$\ov{\mouse^{\gamma_{(\alpha+m)}}}\models a\in H_{\kappa_\alpha^+}$. \ It follows that \
$\ov{\mouse^{\gamma_{(\alpha+m)}}}\models b\in H_{\kappa_\alpha^+}$ \ and, furthermore, it follows that
there is a function \ $f\in \ov{M^{\gamma_{(\alpha+m)}}}$ \ such that \ $f\colon
\kappa_\alpha\times\R\maps{onto} a$. \  Because there is a \ $\Sigma_1(\ov{\mouse^{\gamma_{(\alpha+m)}}})$ \
function \ $F$ \ such that  
\ $F\colon [\OR^{\ov{\mouse^{\gamma_{(\alpha+m)}}}}]^{<\omega} \times
{\R} \maps{onto}\ov{M^{\gamma_{(\alpha+m)}}}$, \ one can assume that the function \ $f$ \ is definable in \
$\ov{\mouse^{\gamma_{(\alpha+m)}}}$ \ where the definition of \ $f$ \ uses only the parameters \
$\kappa_\alpha$, \
$a$ \ and some real \ $y$.\footnote{Let $q$ \ be the $\le_{\BK}$--least such that `$f=F(q,y)$ has the desired
property'.} \ Now let \ $\eta\in\kappa_\alpha$ \ and \
$w\in\R$ \ be such that \ $b=f(\eta,w)$. \ Thus, \ $b$ \ is a definable element in \
$\ov{\mouse^{\gamma_{(\alpha+m)}}}$
\ where the definition uses only the parameters \ $\kappa_\alpha, a, y, \eta, w$. \ Hence, \ $b\in
K^m_\alpha$.

\medskip
(\alph{cntsix}) Finally, we note that \ $K^m_\alpha$ \ is closed under the ``rudimentary'' functions. This holds
because the rudimentary functions are definable and increase rank by a finite ordinal.
This completes the proof of Claim \ref{claimsix}. \end{proof}

\begin{clam}\label{claimlast} Let \ $\alpha$ \ be an ordinal and let \
$\mouse^{\gamma_\alpha+1}=(M^{\gamma_\alpha+1},\R,\kappa_\alpha,\mu_\alpha)$. \  For all \ $m\ge 0$, 
\ $\mu_\alpha\cap K^m_\alpha\in K^{m+2}_\alpha$.
\end{clam}

\begin{proof}[Proof of Claim \ref{claimlast}] Let \ $X=\{a\in K^m_\alpha : a\subseteq \kappa_\alpha\}$. \ By (\alph{cntthree})
and (\alph{cntsix}) of Claim \ref{claimsix}, it follows that \ $X\in K^{m+1}_\alpha$. \ Thus, there is a function \
$f_\alpha\in \ov{M^{\gamma_{(\alpha+k)}}}$ \ such that \ $f_\alpha\colon \kappa_\alpha\times\R\maps{onto}
X$, \ where \
$k=m+1$. \  Since there is a
\ $\Sigma_1(\ov{\mouse^{\gamma_{(\alpha+k)}}})$ \ function \ $F$ \ such that \ 
$F\colon [\OR^{\ov{\mouse^{\gamma_{(\alpha+k)}}}}]^{<\omega} \times
{\R}\maps{onto}\ov{M^{\gamma_{(\alpha+k)}}}$,
\ one can assume that the function \ $f_\alpha$ \ is definable in \ $\ov{\mouse^{\gamma_{(\alpha+k)}}}$
\ where the definition of \ $f_\alpha$ \ uses only the parameters \ $\kappa_\alpha$, \ $X$ \ and some real \ $y$. \ 
Thus, \ $f_\alpha\in
\Sigma_\omega(\ov{\mouse^{\gamma_{(\alpha+k)}}},\R\cup\{\nu,\kappa_\alpha,\dots,\kappa_{\alpha+m}\})$.
\ Hence, Claim \ref{uniformly} also implies that  
\begin{equation*}
f_\alpha\in\Sigma_\omega(\ov{\mouse^{\gamma_{(\alpha+k+1)}}},\R\cup\{\nu,\kappa_\alpha,\dots,\kappa_{\alpha+(m+1)}\}).\tag{$*$}
\end{equation*}
Because \ $f_\alpha$ \ can be coded as a subset of \ $\kappa_\alpha\times\R$, \ Claim
\ref{claim5plus} implies that \ $f_\alpha\in M^{\gamma_{\alpha}+1}$, \ $f_{\alpha+1}\in
M^{\gamma_{(\alpha+1)}+1}$ \ and \ $\pi_{\alpha{(\alpha+1)}}(f_\alpha)=f_{\alpha+1}$, \ where 
\begin{equation*} f_{\alpha+1}\in
\Sigma_\omega(\ov{\mouse^{\gamma_{(\alpha+k+1)}}},\R\cup\{\nu,\kappa_\alpha,\dots,\kappa_{\alpha+(m+1)}\})
\tag{$**$}\end{equation*}
has the same definition as that of \ $f_\alpha$. \ 
So, for any \
$\eta\in\kappa_\alpha\times\R$, 
\[f_\alpha(\eta)\in\mu_\alpha\iff\kappa_\alpha\in
\pi_{\alpha{(\alpha+1)}}(f_\alpha)(\eta)\iff \kappa_\alpha\in f_{\alpha+1}(\eta).\]
Hence, \ $\mu_\alpha\cap X = \{f_\alpha(\eta) : \kappa_\alpha\in f_{\alpha+1}(\eta) \land
\eta\in\kappa_\alpha\times\R\}$.
\ Now, because \ $\mu_\alpha\cap X$ \ can be coded as a subset of \ $\kappa_\alpha\times\R$, \ properties ($*$) and
($**$) imply that \ $\mu_\alpha\cap X\in K^{m+2}_\alpha$. \ Therefore, \ $\mu_\alpha\cap K^m_\alpha\in K^{m+2}_\alpha$.
\end{proof}

\begin{clam}\label{clast} Let \ $\alpha$ \ be an ordinal. Then \ $H_{\kappa_\alpha}^{\mouse^{\gamma_\alpha+1}}=
H_{\kappa_\alpha}^{\mouse^{\gamma_\alpha}}$.
\end{clam}

\begin{proof}[Proof of Claim \ref{clast}] Let \ $\alpha$ \ be an ordinal. Recall that \ $H^1_\alpha=
H_{\kappa_\alpha^+}^{\ov{\mouse^{\gamma_{(\alpha+1)}}}}$.
 \ Define \ $K_\alpha=\bigcup\limits_{m\in\omega}K^m_\alpha$. \ Since \ $K^m_\alpha\subseteq
H^1_\alpha$ \ for each \ $m\in\omega$, \ we see that \ $K_\alpha\subseteq H^1_\alpha$. \ 
Claim \ref{claimlast} and  Claim \ref{claimsix}(\alph{cntfive})(\alph{cntsix}) imply that \ $K_\alpha$ \ is closed under
the ``rudimentary in \ $\mu_\alpha$'' \ functions. Thus, \ $M^{\gamma_\alpha+1}\subseteq K_\alpha$. \ Hence, \
$M^{\gamma_\alpha+1}\subseteq H_{\kappa_\alpha^+}^{\ov{\mouse^{\gamma_{(\alpha+1)}}}}$. \ It follows  that \ 
$H_{\kappa_\alpha}^{\mouse^{\gamma_\alpha+1}}\subseteq
H_{\kappa_\alpha}^{\ov{\mouse^{\gamma_{(\alpha+1)}}}}$. \ Because \
$\ov{M^{\gamma_\alpha}}=H^{\mouse^{\gamma_\alpha}}_{\kappa_\alpha}$, \ Claim \ref{claimfour} implies that \
$H_{\kappa_\alpha}^{\mouse^{\gamma_\alpha+1}}\subseteq H_{\kappa_\alpha}^{\mouse^{\gamma_\alpha}}$. \ Therefore,
$H_{\kappa_\alpha}^{\mouse^{\gamma_\alpha+1}}= H_{\kappa_\alpha}^{\mouse^{\gamma_\alpha}}$.
\end{proof}

Claim \ref{clast} implies that \ $H_{\kappa}^{\mouse^{\gamma+1}}=
H_{\kappa}^{\mouse^{\gamma}}$ \ and this completes the proof of Lemma \ref{lemmatwo}.\end{proof}

\begin{corollary}\label{nobigger} Let \ $\mouse=(M,\R,\kappa,\mu)$ \ be  an iterable real premouse. Let \ $\gamma <
\widehat{\OR}^{\mouse}$. \ Assume that \
${\mouse^\gamma}$ \ is acceptable and critical. Let \ $n=n(\mouse^\gamma)$. \
If \ $\rho_{\mouse^\gamma}^{n+1} = \kappa$, \ then \
$\rho_{\mouse^\gamma}^{m} = \kappa$ \ for all \ $m> n$.
\end{corollary}

\begin{proof}
This follows from Claim \ref{claimthree} of the above proof (in addition, see Lemma 1.31 and Corollary 2.38 of
\cite{Cfsrm}).
\end{proof}

\begin{lemma} \label{lemma3} Let \ $\mouse=(M,\R,\kappa,\mu)$ \ be an iterable real premouse. Then
\begin{list}{{\upshape (\roman{puff})}}
{\setlength{\labelwidth}{20pt}
\usecounter{puff}}
\item $\mouse^\lambda=(J_\lambda(\R),\R)$ \ is acceptable for all $\lambda\le\kappa$,
\item $\rho^{}_{\mouse^\kappa}=\omega\rho^{}_{\mouse^\kappa}=\kappa$,
\item $\mouse^{\kappa+1}$ \ is acceptable.
\end{list}
\end{lemma}

\begin{proof} For (i), one can show that \ $J_\lambda(\R)$ \ is acceptable (above the reals) for all \
$\lambda\le\kappa$ \ by generalizing, for example, the proof of Lemma 4.27 of \cite{Dodd}. For (ii) note that since \
$\mouse^{\kappa+1}$ \ is iterable, one can assume that \ $\kappa$ \ is a \ $\Lr$--cardinal. For any \ $\xi<\kappa$, \ a
standard  ``$\Lr$ condensation argument'' \ shows that any \ $\boldface{\Sigma}{1}(J_\kappa(\R))$ \ subset of \
$\R\times\xi$ \ is an element of \ $J_\kappa(\R)$. \ It follows that \
$\rho^{}_{J_\kappa(\R)}=\omega\rho^{}_{J_\kappa(\R)}=\kappa$. \ Finally, a proof of (iii) can be obtained as a special
case of the proof of Lemma \ref{lemmatwo}.
\end{proof}

\begin{lemma}\label{lemma4} Let \ $\mouse=(M,\R,\kappa,\mu)$ \ be a real premouse. Let \ $\gamma$ \ be such
that \ $\kappa<\gamma<\widehat{\OR}^{\mouse}$. \ Suppose that \ $\mouse^\gamma$ \ is acceptable.  If \
$\rho_{\mouse^\gamma}^{n} > \delta\ge\kappa$ \ for all \ $n\in\omega$, \ then
\[M^{\gamma+1}\cap\pow(\R\times\delta)= \boldface{\Sigma}{\omega}(\mouse^{\gamma})\cap\pow(\R\times\delta)=
M^{\gamma}\cap\pow(\R\times\delta).\]
\end{lemma}

\begin{proof} Let \ $\mu^{\gamma+1}=M^{\gamma+1}\cap\mu$. \
Recalling Definition \ref{preddef}, we prove the following claim.
\begin{claiim} $\mu^{\gamma+1}$ \ is predictable.  
\end{claiim}

\begin{proof}[Proof of Claim]
Let \ $\chi(v_0,v_1,\dots,v_k)\in\Sigma_\omega$ \ be an \ $\Lng$
\ formula.  We will define  an \ $\Lng$ \ formula \ $\psi(v_1,\dots,v_k)\in\Sigma_\omega$ \ and show that 
\begin{equation}
B_{a_1,\dots,a_k}\in\mu^{\gamma+1}\iff \mouse^\gamma\models\psi(a_1,\dots,a_k)\label{stupid}
\end{equation}
for all \ $a_1,\dots,a_k\in M^\gamma$, \ where \ $B_{a_1,\dots,a_k}=\{\varkappa\in\kappa : \mouse^\gamma\models
\chi(\varkappa,a_1,\dots,a_k)\}$. \  Since \ $\mouse^\gamma$ \ is acceptable and \ $\rho_{\mouse^\gamma}^{n} >
\delta\ge\kappa$ \ for all \ $n\in\omega$, \  Corollary \ref{gammaproj} and Theorem \ref{kappasound} imply that \
$\gamma_{\mouse^\gamma}^n=\rho_{\mouse^\gamma}^n$ \  for all \ $n\in\omega$. \ So, \ $\gamma_{\mouse^\gamma}^{n} >
\delta\ge\kappa$ \ for all \ $n\in\omega$ (recall Definition
\ref{defproj}). It follows that \ $B_{a_1,\dots,a_k}\in M^\gamma$ \ for all \ $a_1,\dots,a_k\in M^\gamma$.
Let \ $\psi(v_1,\dots,v_k)$ \ be the formula 
\[(\exists A)\left[\,(\forall\varkappa\in\kappa)(\varkappa\in A \leftrightarrow \chi(\varkappa,v_1,\dots,v_k))\land
A\in\mu\,\right]\] and let \ $a_1,\dots,a_k\in M$. \  Since \ $B_{a_1,\dots,a_k}\in M^\gamma$, \ (\ref{stupid}) now follows.
This completes the proof of the Claim.\end{proof} 
 
Lemma \ref{predlemma} and the above Claim imply that $\pow(M^\gamma)\cap M^{\gamma+1}\subseteq
\boldface{\Sigma}{\omega}(\mouse^\gamma)$.
\ Thus, \ $M^{\gamma+1}\cap\pow(\R\times\delta)=
\boldface{\Sigma}{\omega}(\mouse^{\gamma})\cap\pow(\R\times\delta)$. \ In addition, because \
$\gamma_{\mouse^\gamma}^{n} > \delta$ \ for all \ $n\in\omega$, \ it follows that \
$\boldface{\Sigma}{\omega}(\mouse^{\gamma})\cap\pow(\R\times\delta)=M^{\gamma}\cap\pow(\R\times\delta)$.
\end{proof}

\begin{lemma}\label{lemma4.5} Let \ $\mouse=(M,\R,\kappa,\mu)$ \ be an iterable real premouse. Let \ $\gamma$ \ be such
that \ $\kappa<\gamma<\widehat{\OR}^{\mouse}$. \ Assume that \
${\mouse^\gamma}$ \ is real mouse.  If \ $\rho_{\mouse^\gamma}^{n} > \delta$ \ for all \ $n\in\omega$, \
then
\begin{equation}
M^{\gamma+1}\cap\pow(\R\times\delta)= \boldface{\Sigma}{\omega}(\mouse^{\gamma})\cap\pow(\R\times\delta)=
M^{\gamma}\cap\pow(\R\times\delta).\label{bigdeal}\end{equation}
\end{lemma}

\begin{proof} Since \ $\mouse$ \ is a mouse, we have that \ $\rho_{\mouse^\gamma}^{n+1} \le
\kappa < \rho_{\mouse^\gamma}^n$ \ where \ $n=n({\mouse^\gamma})$. \ Hence, \ $\delta< \rho_{\mouse^\gamma}^{n+1} \le
\kappa < \rho_{\mouse^\gamma}^n$. \ There two cases to consider.

\smallskip
\noindent{\sc Case 1:} $\rho_{\mouse^\gamma}^{n+1} = \kappa$. \ 
Lemma \ref{lemmatwo} \ implies that \ $H^{\mouse^\gamma}_{\kappa}=H^{\mouse^{\gamma+1}}_{\kappa}$. \ Because \
$\delta<\kappa$, \ we conclude that equation (\ref{bigdeal}) holds.

\smallskip
\noindent{\sc Case 2:} $\rho_{\mouse^\gamma}^{n+1} < \kappa$. \ 
Since \ $\mouse^\gamma$ \ is mouse and \ $\rho_{\mouse^\gamma}^{n+1} < \kappa$, \ 
Lemma \ref{lemma1} implies that \ $\pow(M^\gamma)\cap M^{\gamma+1}\subseteq
\boldface{\Sigma}{\omega}(\mouse^\gamma)$. \ Thus, \ $M^{\gamma+1}\cap\pow(\R\times\delta)=
\boldface{\Sigma}{\omega}(\mouse^{\gamma})\cap\pow(\R\times\delta)$. \
In addition, because \ $\rho_{\mouse^\gamma}^{n} > \delta$ \ for all \ $n\in\omega$, \  Theorem \ref{gamma} implies, again
because \ $\mouse^\gamma$ \ is mouse, that \ $\gamma_{\mouse^\gamma}^n=\rho_{\mouse^\gamma}^n$ \  for all \ $n\in\omega$. \ So,
\ $\gamma_{\mouse^\gamma}^{n} > \delta$ \ for all \ $n\in\omega$. \ Hence,
\ $\boldface{\Sigma}{\omega}(\mouse^{\gamma})\cap\pow(\R\times\delta)= M^{\gamma}\cap\pow(\R\times\delta)$. \ Therefore,
equation (\ref{bigdeal}) holds.
\end{proof}

\begin{definition} Let \ $\mouse=(M,\R,\kappa,\mu)$ \ be an iterable real premouse. \ Let \ $\gamma < \widehat{\OR}^{\mouse}$ \
and let \  $\delta<{\OR}^{\mouse^{\gamma}}$. \ We say that \ $\mouse^{\gamma+1}$ \ is {\it acceptable at\/}  \ $\delta$
\ provided that, if \ $\mathcal{P}(\delta\times{\R})\cap M^{\gamma+1}\nsubseteq M^\gamma$, \ then for each \
$u\in M^{\gamma+1}$ \ there is an \ $f\in M^{\gamma+1}$ \ such that
\ $f=\langle f_{\xi} : \delta\le \xi<{\OR}^{\mouse^\gamma}\rangle$ \
and
\begin{equation*}f_{\xi}\colon \xi\times\R \maps{onto}
\{\xi\}\,\cup\,(\mathcal{P}(\xi\times\R)\cap u).
\end{equation*} 
\end{definition} 

\begin{lemma}\label{lemma5} Let \ $\mouse=(M,\R,\kappa,\mu)$ \ be an iterable real premouse. Let \ $\gamma$ \ be such
that \ $\kappa<\gamma<\widehat{\OR}^{\mouse}$. \ Assume that \
${\mouse^\gamma}$ \ is acceptable.  Then  \  $\mouse^{\gamma+1}$ \ is acceptable at \
$\delta$ \ whenever \ $\kappa\le\delta<{\OR}^{\mouse^{\gamma}}$.
\end{lemma}

\begin{proof} Let \ $\kappa\le\delta<{\OR}^{\mouse^{\gamma}}$ \ and assume that there is an \ $a\subseteq
\delta\times{\R}$ \ such that \ $a\in M^{\gamma+1}\setminus M^\gamma$.
\ Lemma \ref{lemma4} implies that \  $\rho_{\mouse^{\gamma}}^{m}\le\delta$ \ for some \ $m\ge 1$. \ Let \ $k$ \ be 
such that $\rho_{\mouse^{\gamma}}^{k+1}\le\delta<\rho_{\mouse^{\gamma}}^{k}$. \ Let \ $\lambda=\textup{max}\{\kappa,
\rho_\mouse^{k+1}\}\le\delta$. \ Lemma \ref{SCL} implies that there is a \ $\boldface{\Sigma}{k+1}(\mouse^{\gamma})$ \ function
\ $F\colon \lambda\times\R \maps{onto} [\omega\gamma]^{<\omega}\times\R$, \ since \ $\omega\gamma=\OR^{\mouse^\gamma}$. 
\ Now let \ $u\in M^{\gamma+1}$ \ be nonempty. So, \
$u\in S_{\omega\gamma+n}(\R)$ 
\ for some \ $n\in\omega$, \ where \ $S_{\omega\gamma+n}(\R)=S_{\omega\gamma+n}^{\mouse^{\gamma+1}}(\underline{\R})$. \
Lemma 1.7 of
\cite{Crcm} states that there is a function \ $g\in M^{\gamma+1}$
\ such that
\ $g\colon[\omega\gamma+n]^{<\omega}\times\R\maps{onto}S_{\omega\gamma+n}(\R)$ \ and consequently, one can easily
construct a function \
$h\in M^{\gamma+1}$ \ such that \ $h\colon[\omega\gamma]^{<\omega}\times\R\maps{onto}S_{\omega\gamma+n}(\R)$. \ Hence \
$h\circ F\colon \lambda\times\R \maps{onto} S_{\omega\gamma+n}(\R)$. \ So one  can  produce in \ $\mouse^{\gamma+1}$ \  a
function \ $f'\colon \lambda\times\R \maps{onto} u$. \ Using \ $f'$,
\ one can easily define the desired sequence of functions \ $\langle f_{\xi} : \delta\le \xi<{\OR}^{\mouse^\gamma}\rangle\in
M^{\gamma+1}$ \ such that \ $f_{\xi}\colon \xi\times\R \maps{onto} \{\xi\}\cup(\mathcal{P}(\xi\times\R)\cap u)$. \  For
example, given a desired \ $\xi>\delta\ge\lambda$ \ define
\begin{equation*}
f_{\xi}(\alpha, y)=
\begin{cases}
\xi,            &\text{if $\alpha\ge\lambda$;}\\
\xi,            &\text{if $\alpha<\lambda \,\land\, f'(\alpha, y)\notin\mathcal{P}(\xi\times\R)\cap u$;}\\
f'(\alpha, y),  &\text{if $\alpha<\lambda \,\land\, f'(\alpha, y)\in\mathcal{P}(\xi\times\R)\cap u$.}
\end{cases}
\end{equation*} 
This completes the proof.
\end{proof}

\begin{lemma}\label{lemma6} Let \ $\mouse=(M,\R,\kappa,\mu)$ \ be an iterable real premouse. Let \ $\gamma$ \ be such
that \ $\kappa<\gamma<\widehat{\OR}^{\mouse}$. \ Assume that \
${\mouse^\gamma}$ \ is acceptable.  Then for all \ $\delta<\kappa$, \  $\mouse^{\gamma+1}$ \ is acceptable at \
$\delta$.
\end{lemma}

\begin{proof} Let \ $\delta<\kappa$ \ and assume that there is an \ $a\subseteq
\delta\times{\R}$ \ such that \ $a\in M^{\gamma+1}\setminus M^\gamma$. \ 
Lemma \ref{lemma1}, Lemma \ref{lemmatwo} and Corollary \ref{nobigger} imply that there exists a \ $k\in\omega$ \ such that
\ $\rho_{\mouse^{\gamma}}^{k+1}\le\kappa<\rho_{\mouse^{\gamma}}^{k}$. \ Now let \ $u\in M^{\gamma+1}$ \ be nonempty.
Since \ $\mouse^\gamma$ \ is acceptable and \ $\rho_{\mouse^{\gamma}}^{k+1}\le\kappa<\rho_{\mouse^{\gamma}}^{k}$, \  we see
that \ $\mouse^\gamma$ \ is critical.  Theorem \ref{intismi} implies that \ $\mouse^\gamma$ \ is a real mouse and
thus, Lemma \ref{lemma4.5} implies that \  $\rho_{\mouse^{\gamma}}^{m+1}\le\delta$ \ for some \ $m\ge 0$. \ Let \ $m$ \
be the least such integer. Note that \ $m\ge k$. \ We need to define a sequence \ $\langle f_{\xi} : \delta\le\xi<
{\OR}^{\mouse^\gamma}\rangle\in M^{\gamma+1}$ \ such that \ $f_{\xi}\colon \xi\times\R \maps{onto}
\{\xi\}\cup(\mathcal{P}(\xi\times\R)\cap u)$. 

We first outline how we shall obtain the desired sequence. For ordinals \ $\alpha<\beta$ \ define the interval \
$[\alpha,\beta)=\{\xi : \alpha\le\xi<\beta\}$. \ If \ $m=k$ \ then it is sufficient to define our sequence on the
interval \  $[\rho_{\mouse^{\gamma}}^{k+1},{\OR}^{\mouse^\gamma})$. \ In this case, we shall define the sequence in two pieces:
first on the interval \ $[\kappa,{\OR}^{\mouse^\gamma})$ \ and then on the interval \ $[\rho_{\mouse^{\gamma}}^{k+1},\kappa)$.
\ These two sequences will be in \ $\mouse^{\gamma+1}$ \ and so their union will be the required sequence defined on the
entire interval \ $[\rho_{\mouse^{\gamma}}^{k+1},{\OR}^{\mouse^\gamma})$. \ If \ $m>k$, \ then we shall construct the required
sequence in ``finitely many pieces''.  Note that  \ $\rho_{\mouse^{\gamma}}^{m+1}\le\rho_{\mouse^{\gamma}}^{m}\le
\cdots\le \rho_{\mouse^{\gamma}}^{k+1}<\kappa<{\OR}^{\mouse^\gamma}$ \ is a finite increasing sequence. Without loss of 
generality, we shall assume that this finite sequence is strictly increasing, that is,
\[\rho_{\mouse^{\gamma}}^{m+1}<\rho_{\mouse^{\gamma}}^{m}<
\cdots< \rho_{\mouse^{\gamma}}^{k+1}<\kappa<{\OR}^{\mouse^\gamma}.\] 
So, in addition to the above two intervals, we shall construct sequences in \ $\mouse^{\gamma+1}$ \ on each of the intervals \
$[\rho_{\mouse^{\gamma}}^{i+1}, \rho_{\mouse^{\gamma}}^{i})$ \ for \ $m\le i \le k+1$. \ In the end, the union of this finite
set of sequences produces our required sequence defined on the entire interval
$[\rho_{\mouse^{\gamma}}^{m+1},{\OR}^{\mouse^\gamma})$.  \ Now we construct these sequences.

\smallskip
{\sc Case 1:} {\sl The interval \ $[\kappa,{\OR}^{\mouse^\gamma})$.} \ Using the argument in the proof of Lemma \ref{lemma5},
there is a sequence \ $\langle f_{\xi} : \kappa\le \xi<{\OR}^{\mouse^\gamma}\rangle\in M^{\gamma+1}$ \ such that \
$f_{\xi}\colon \xi\times\R \maps{onto} \{\xi\}\cup(\mathcal{P}(\xi\times\R)\cap u)$.

\smallskip
{\sc Case 2:} {\sl The interval \ $[\rho_{\mouse^{\gamma}}^{k+1},\kappa)$}. \ We need to construct a sequence \ $\langle
f_{\xi} :
\rho_{\mouse^{\gamma}}^{k+1}\le\xi< \kappa\rangle\in M^{\gamma+1}$ \ such that \
$f_{\xi}\colon \xi\times\R \maps{onto} \{\xi\}\cup(\mathcal{P}(\xi\times\R)\cap u)$. \ Thus we may now assume that \
$u\subseteq \mathcal{P}(\kappa\times\R)\cap M^{\gamma+1}$. \ Note that \
$\rho_{\mouse^{\gamma}}^{k+1}=\rho^{}_{\,\overline{\mouse^\gamma}}$ \ and \
$p_{\mouse^{\gamma}}^{k+1}=p^{}_{\,\overline{\mouse^\gamma}}$. \  Let \ $\core =
\core(\mouse^\gamma)$. \ Since \
$\mouse^\gamma$ \ is a real mouse, it is a mouse iterate of \ $\core$. \ Let
\[\premouseiterationover{\overline{\core}}{\alpha}{\beta}\]
be the premouse iteration of \ $\overline{\core}$ \ and let \ $\kappa_\alpha =
\overline{\pi}_{0\alpha}({\kappa}^{\core})$ \ for \ $\alpha\in \OR$. \ Also, we shall let \ $\overline{C}_\alpha$ \ denote
the domain of the structure \ $\overline{\core}_\alpha$. \ We make the following observations:
\be
\item $\overline{\core}_\kappa=\overline{\mouse^\gamma}$ \ and \ $\overline{\pi}_{0\kappa}({\kappa}^{\core})=\kappa$, \ by
Lemma \ref{deftwo}
\item $\pow(\R\times\kappa)\cap M^{\gamma+1}\subseteq
\boldface{\Sigma}{\omega}\!\left(\,\overline{\mouse^\gamma}\,\right)$, \ by 
Lemma 2.19 of
\cite{Cfsrm} and Lemma \ref{deftwo} (above)
\item $I_0\in\boldface{\Sigma}{\omega}\!\left(\,\overline{\mouse^\gamma}\,\right)$ \ where \ $I_0=\{\kappa_\alpha : \alpha
<\kappa\}$, \ by Lemma \ref{defone}
\item $\rho^{}_{\,\overline{\core}}=\rho^{}_{\,\overline{\core}_\alpha}=\rho^{}_{\,\overline{\mouse^\gamma}}$ \ for all \
$\alpha$, \ by Lemma \ref{relation1} and Theorem \ref{relationtwo}
\item $\overline{\pi}_{\alpha\beta}(p^{}_{\,\overline{\core}_\alpha})=p^{}_{\,\overline{\core}_\beta}$ \ for all \
$\alpha\le\beta$, \ by Lemma \ref{relation1}.
\ee
In addition, the proof of Lemma \ref{lemma5} can be used to show that \ $u$ \ can be coded by a
single subset of \ $\kappa\times\R$ \ in \ $\mouse^{\gamma+1}$. \ Thus, \
$u\subseteq\boldface{\Sigma}{m}\!\left(\,\overline{\mouse^\gamma}\,\right)$ \ for some fixed \ $m\in\omega$. \ 

Before we construct our desired sequence of functions \ $\langle f_{\xi} :
\rho_{\mouse^{\gamma}}^{k+1}\le\xi< \kappa\rangle\in M^{\gamma+1}$ \ such that \
$f_{\xi}\colon \xi\times\R \maps{onto} \{\xi\}\cup(\mathcal{P}(\xi\times\R)\cap u)$, \ we shall first discuss a method for 
constructing a single such function.
To do this, let \ $\xi$ \ be fixed such that
\ $\rho^{}_{\,\overline{\mouse^\gamma}}\le\xi<\kappa$ \ and let \ $A_\xi=\mathcal{P}(\xi\times\R)\cap u$. \
So, for each  \ $a\in A_\xi$, \ $a\subseteq\xi\times\R$. \  We will construct an individual function \
$f\in M^{\gamma+1}$
\ such that \ $f\colon\xi\times\R\maps{onto}A_\xi$. \  Let \
$\kappa_\alpha$ \ be the least such that \ $\xi\le\kappa_\alpha<\kappa$. \ It follows from Theorem \ref{relationtwo} and Lemma
2.19 of \cite{Cfsrm} that \ 
$a\in\boldface{\Sigma}{m}(\overline{\core}_\alpha)$ \ for each \ $a\in A_\xi$. \ Lemma 2.16 of \cite{Cfsrm} implies that \
\begin{equation}\overline{\core}_\alpha=\Hull_1^{\overline{\core}_\alpha}(\R\cup \kappa_\alpha\cup
\{p^{}_{\,\overline{\core}_\alpha}\}).\tag{$\star$}\end{equation} 

\begin{claiim}
$\overline{\core}_\alpha=\Hull_1^{\overline{\core}_\alpha}(\R\cup \xi\cup
\{p^{}_{\,\overline{\core}_\alpha}\})$.
\end{claiim} 

\begin{proof}
If \ $\xi=\kappa_\alpha$, \ then the Claim follows from ($\star$). If \ $\xi\ne\kappa_\alpha$, \ then
 \ $\rho^{}_{\,\overline{\core}}\le\xi<\kappa_\alpha$ \ (note that 
$\rho^{}_{\,\overline{\core}}=\rho^{}_{\,\overline{\mouse^\gamma}}$). Since \ $\kappa_\alpha$ \ was chosen to be the least such
that \ $\xi\le\kappa_\alpha$, \ it follows that either \
$\alpha=0$ \ or \ $\alpha$ \ is a successor ordinal. If \ $\alpha=0$ \ then the Claim follows from the definition of the core
of \ $\mouse^\gamma$, \ because \ $\xi\ge\rho^{}_{\,\overline{\mouse^\gamma}}$. \ Suppose now that \
$\kappa_\alpha=\kappa_{\eta+1}$ \ for some \ $\eta$. \ Hence, \ $\overline{\pi}_{\eta\alpha}\colon\overline{\core}_\eta\mapsigma{1}\overline{\core}_\alpha$ \  and \
$\overline{\pi}_{\eta\alpha}(p^{}_{\,\overline{\core}_\eta})=p^{}_{\,\overline{\core}_\alpha}$. \ Because \
$\overline{\core}_\eta=\Hull_1^{\overline{\core}_\eta}(\R\cup \kappa_\eta\cup \{p^{}_{\,\overline{\core}_\eta}\})$,  \
Lemma 2.8(3) of \cite{Crcm} implies that \ 
$\overline{\core}_\alpha=\Hull_1^{\overline{\core}_\alpha}(\R\cup \xi\cup
\{p^{}_{\,\overline{\core}_\alpha}\})$ \ since \ $\kappa_\eta<\xi$. \ This completes the proof of the Claim.
\end{proof}

Thus, \ $\overline{\core}_\alpha\cong\Hull_1^{\overline{\mouse^\gamma}}(\R\cup \xi\cup
\{p^{}_{\,\overline{\mouse^\gamma}}\})$ \ via the embedding \
$\overline{\pi}_{\alpha\kappa}\colon\overline{\core}_\alpha\mapsigma{1} \overline{\mouse^\gamma}$,
\ because \ $\overline{\core}_\kappa=\overline{\mouse^\gamma}$. \ Hence, \ $a\in\boldface{\Sigma}{m}(\mathcal{H}_\xi)$
\ for each \ $a\in A_\xi$, \ where \ $\mathcal{H}_\xi=\Hull_1^{\overline{\mouse^\gamma}}(\R\cup \xi\cup
\{p^{}_{\,\overline{\mouse^\gamma}}\})$. \ In addition, \ $\mathcal{H}_\xi\in\mouse^{\gamma+1}$ \ and because there is a
\ $\Sigma_1$ \ Skolem function for \ $\mathcal{H}_\xi$, \ Lemma \ref{skolemimage} implies that there is a {\it
canonical\/} total function \ $h\in\mouse^{\gamma+1}$ \ such that  \ $h^{\prime\prime}(\R\times(\xi)^{<\omega})=H_\xi$,
\ where \ $H_\xi$ \ denotes the domain of the structure \ $\mathcal{H}_\xi$. \ Now, since \ $\xi\le\kappa_\alpha$ \ it
follows that \ $\xi\in\overline{C}_\alpha$. \ Hence, Lemma 1.4 of \cite{Cfsrm} implies that there is a function \ $g\in
H_\xi$ \ such that \ $g\colon\xi\times\R\maps{onto}(\xi)^{<\omega}$. \ Therefore,  the function \ $f'\in M^{\gamma+1}$ \
defined by \ $f'(\zeta,y)=h(y_1,g(\zeta,y_2))$ \ is such that \ $f'\colon\xi\times\R\maps{onto}H_\xi$. \ Lemma
\ref{satisfaction} implies that the \ $\Sigma_m$ \ satisfaction relation over \ $\mathcal{H}_\xi$ \ is \
$\Sigma_m(\mathcal{H}_\xi)$. \ Let \ $\textup{Sat}=\textup{Sat}_m^3$ \ denote the \ $\Sigma_m$ \ satisfaction relation
over \ $\mathcal{H}_\xi$ \ ranging over formulae of three variables.  Define the function \ $\mathfrak{a}$ \ by \
$\mathfrak{a}(\zeta,y)=\{ v\in\R\times\xi : \textup{Sat}(y_1,v,f'(\zeta,y_2))\}$ \ for \ $y\in\R$ \ and \
$\rho^{}_{\,\overline{\mouse^\gamma}}\le\zeta<\kappa$.  \ Note that the function \ $\mathfrak{a}$ \ is in \
$\mouse^{\gamma+1}$. \ Now define \ $f$ \ by
\begin{equation}\label{pless}
f(\zeta, y)=
\begin{cases}
\xi,            &\text{if $\zeta<\rho^{}_{\,\overline{\mouse^\gamma}}$;}\\
\xi,            &\text{if $\rho^{}_{\,\overline{\mouse^\gamma}}\le\zeta<\kappa 
\,\land\, \mathfrak{a}(\zeta,y)\notin\mathcal{P}(\xi\times\R)\cap u$;}\\  \mathfrak{a}(\zeta,y),  &\text{if
$\rho^{}_{\,\overline{\mouse^\gamma}}\le\zeta<\kappa \,\land\, \mathfrak{a}(\zeta,y)\in\mathcal{P}(\xi\times\R)\cap u$.}
\end{cases}
\end{equation} 
It follows that \ $f\in\mouse^{\gamma+1}$ \ and \ $f\colon \xi\times\R \maps{onto}
\{\xi\}\cup(\mathcal{P}(\xi\times\R)\cap u)$.
\ This completes our discussion of the method used to construct the individual function \ $f$. 

We note that the above 
construction of \ $f$ \ depends only on \ $\xi$ \ and the choice of the function \ $g\in H_\xi$. \ One can choose, however,
such a function \ $g$ \ which is definable over \ $\mathcal{H}_\xi$ \ in {\it some\/}
real parameter as we shall now show. First let \ $g'\in \mouse^{\gamma+1}$ \ be such that \
$g'\colon\xi\times\R\maps{onto}(\xi)^{<\omega}$. \ Since there is a uniform \ $\Sigma_1(\mathcal{H}_\xi)$ \ function \
$F$ \ such that \ $F\colon [\OR^{\mathcal{H}_\xi}]^{<\omega} \times {\R} \maps{onto} H_\xi$, \ let \ $x\in\R$ \ and \ $b'$ \ 
be such that \ $g'=F(b',x)$. \ Let \ $b$ \ be the $\le_{\BK}$--least such that \ $g=F(b,x)$ \ is such that \ 
$g\colon\xi\times\R\maps{onto}(\xi)^{<\omega}$. \ Let \ $g_{\xi,x}$ \ denote the function \ $g$ \ obtained in this way.  Even though there always exists reals \ $x$ \ such that \ $g_{\xi,x}$ \ is ``defined'' (that
is, \ $g_{\xi,x}\colon\xi\times\R\maps{onto}(\xi)^{<\omega}$), \ there may be reals \ $y$ \ such that \  $g_{\xi,y}$
\ is not defined and hence, in this case we let \ $g_{\xi,y}(\gamma,z)=\emptyset$ \ for all \ $\gamma\in\xi$, $z\in\R$.
\ Define the function \ $f'_\xi\in M^{\gamma+1}$ \ by \ $f'_\xi(\zeta,y)=h(y_1,g_{\xi,y_3}(\zeta,y_2))$ \ and define the
function \ $\mathfrak{a}_\xi$ \ by \ $\mathfrak{a}_\xi(\zeta,y)=\{ v\in\R\times\xi : \textup{Sat}(y_1,v,f'_\xi(\zeta,y_2))\}$.
\  Now define the function \ $f_{\xi}$ \ as in
(\ref{pless}). It follows that \ $f_\xi\in\mouse^{\gamma+1}$ \ and \ $f_\xi\colon \xi\times\R \maps{onto}
\{\xi\}\cup(\mathcal{P}(\xi\times\R)\cap u)$. \ Since the function \ $f_\xi$ \ is uniformly definable over \ $\mathcal{H}_\xi$,
\ one can obtain a sequence
\  $\langle f_{\xi} : \rho_{\mouse^{\gamma}}^{k+1}\le\xi< \kappa\rangle\in M^{\gamma+1}$ \ such that \
$f_{\xi}\colon \xi\times\R \maps{onto} \{\xi\}\cup(\mathcal{P}(\xi\times\R)\cap u)$.

\smallskip
{\sc Case 3:} {\sl The interval \ $[\rho_{\mouse^{\gamma}}^{i+1}<\rho_{\mouse^{\gamma}}^{i})$ \ for \ $k+1\le i \le m$}. \  
Let \ $\core = \core(\mouse^\gamma)$ \ and assume the notation that was presented in Case 2.  \ Note that
\ $\rho_{\mouse^{\gamma}}^{i+1}=\rho_{\,\overline{\mouse^\gamma}}^{\,j+1}=\rho_{\,\overline{\core}}^{\,j+1}$ \ and \
$p_{\mouse^{\gamma}}^{i+1}=p_{\,\overline{\mouse^\gamma}}^{\,j+1}=\overline{\pi}_{0\kappa}(p_{\,\overline{\core}}^{\,j+1})$ \
whenever \ $k+1\le i \le m$, \ where  \ $j=i-(k+1)$. \ Let such an \ $i$ \ and \ $j=i-(k+1)$ \ be fixed. 
We need to obtain a sequence \ 
$\langle f_{\xi} : \rho_{\mouse^{\gamma}}^{i+1}\le\xi< \rho_{\mouse^{\gamma}}^{i}\rangle\in M^{\gamma+1}$ \ such that \
$f_{\xi}\colon \xi\times\R \maps{onto}
\{\xi\}\cup(\mathcal{P}(\xi\times\R)\cap u)$. \ Thus we may assume that \ $u\subseteq 
\mathcal{P}(\rho_{\mouse^{\gamma}}^{i}\times\R)\cap M^{\gamma+1}$ \ and, as in Case 2, 
 \ $u\subseteq\boldface{\Sigma}{m}\!\left(\,\overline{\mouse^\gamma}\,\right)$ \ for some fixed \ $m\in\omega$. \  
Theorem \ref{relationtwo} and Lemma 2.19 of
\cite{Cfsrm} imply that  \ $u\subseteq\boldface{\Sigma}{m}(\overline{\core})$.

We first construct the sequence \ $\langle f_{\xi}:
\rho_{\,\overline{\core}}^{\,j+1}\le\xi<\rho_{\,\overline{\core}}^{\,j}
\rangle$ \ such that \ $f_{\xi}\colon \xi\times\R \maps{onto}
\{\xi\}\cup(\mathcal{P}(\xi\times\R)\cap u)$. \ Afterwards,  we will show that this sequence is in \
$\mouse^{\gamma+1}$. \  Now, Lemma 2.34 and Corollary 1.32 of \cite{Cfsrm} imply since there is a \
$\boldface{\Sigma}{\omega}(\overline{\core})$ \ total function \ $h$ \ such that  \
$h^{\prime\prime}(\R\times(\rho_{\,\overline{\core}}^{\,j+1})^{<\omega})=\overline{C}$, \ where \ $\overline{C}$ \ is
the domain of the structure \ $\overline{\core}$. \ Lemma 1.4 of
\cite{Cfsrm} implies that there is a function \ $g\in \overline{C}$ \ such that \
$g\colon\rho_{\,\overline{\core}}^{\,j+1}\times\R\maps{onto}(\rho_{\,\overline{\core}}^{\,j+1})^{<\omega}$.
\ Therefore,  the function \ $f'\in \boldface{\Sigma}{\omega}(\overline{\core})$, \ defined by \
$f'(\zeta,y)=h(y_1,g(\zeta,y_2))$, \ is such that \
$f'\colon\rho_{\,\overline{\core}}^{\,j+1}\times\R\maps{onto}\overline{C}$. \ Lemma \ref{satisfaction} implies that the \
$\Sigma_m$ \ satisfaction relation over \ $\overline{\core}$ \ is \
$\Sigma_m(\overline{\core})$. \ Let \ $\textup{Sat}=\textup{Sat}_m^3$ \ denote the \ $\Sigma_m$ \
satisfaction relation over \ $\overline{\core}$ \ ranging over formulae of three variables.
Let \ $\xi$ \ be such that \ $\rho_{\,\overline{\core}}^{\,j+1}\le\xi<\rho_{\,\overline{\core}}^{\,j}$. \ For each \ $y\in\R$ \
and
\ $\zeta<\rho_{\,\overline{\core}}^{\,j+1}$ \ define \ $\mathfrak{a}_\xi(\zeta,y)=\{ v\in\R\times\xi :
\textup{Sat}(y_1,v,f'(\zeta,y_2))\}$. \ Note that the function \ $\mathfrak{a}_\xi$ \ is onto \
$\mathcal{P}(\xi\times\R)\cap\boldface{\Sigma}{m}(\overline{\core})$. \ Define \ $f_\xi$ \ by
\begin{equation}
f_\xi(\zeta, y)=
\begin{cases}
\xi,            &\text{if $\zeta\ge\rho_{\,\overline{\core}}^{\,j+1}$;}\\
\xi,            &\text{if $\zeta<\rho_{\,\overline{\core}}^{\,j+1} 
\,\land\, \mathfrak{a}_\xi(\zeta,y)\notin\mathcal{P}(\xi\times\R)\cap u$;}\\  \mathfrak{a}_\xi(\zeta,y),  &\text{if
$\zeta<\rho_{\,\overline{\core}}^{\,j+1} \,\land\, \mathfrak{a}_\xi(\zeta,y)\in\mathcal{P}(\xi\times\R)\cap u$.}
\end{cases}
\end{equation} 
for each \ $y\in\R$ \ and \ $\zeta\in\xi$. \ Clearly  \ $f_\xi\colon \xi\times\R \maps{onto}
\{\xi\}\cup(\mathcal{P}(\xi\times\R)\cap u)$. \ Thus, we have constructed a sequence \ 
$\langle f_{\xi} : \rho_{\,\overline{\core}}^{\,j+1}\le\xi<\rho_{\,\overline{\core}}^{\,j}
\rangle$ \ such that \ $f_{\xi}\colon \xi\times\R \maps{onto}
\{\xi\}\cup(\mathcal{P}(\xi\times\R)\cap u)$.

Now we shall observe that above sequence is, in fact, in \ $\mouse^{\gamma+1}$. \ First note that any rudimentarily closed
transitive structure that contains \ $u$, $\overline{\core}$ \ as elements, also contains the above sequence as an element. 
Recall, by definition of the core, that \ $\overline{\core}=\Hull_1^{\overline{\core}}(\R\cup \rho^1_{\,\overline{\core}}\cup
\{p^1_{\,\overline{\core}}\})$. \ Thus, \ $\overline{\core}\cong\Hull_1^{\overline{\mouse^\gamma}}(\R\cup 
\rho^1_{\,\overline{\mouse^\gamma}}\cup
\{p^1_{\,\overline{\mouse^\gamma}}\})$  \ via the embedding \
$\overline{\pi}_{0\kappa}\colon\overline{\core}\mapsigma{1} \overline{\mouse^\gamma}$,
\ since \ $\overline{\core}_\kappa=\overline{\mouse^\gamma}$, \
$\rho^1_{\,\overline{\core}}=\rho^1_{\,\overline{\mouse^\gamma}}$ \ and \
$\overline{\pi}_{0\kappa}(p^1_{\,\overline{\core}})=p^1_{\,\overline{\mouse^\gamma}}$.   \ Because \
$\overline{\mouse^\gamma}\in\mouse^{\gamma+1}$, \ it follows that the  sequence \ 
$\langle f_{\xi} : \rho_{\,\overline{\core}}^{\,j+1}\le\xi<\rho_{\,\overline{\core}}^{\,j}\rangle=
\langle f_{\xi} : \rho_{\mouse^{\gamma}}^{i+1}\le\xi< \rho_{\mouse^{\gamma}}^{i}\rangle$ \ is in \
$\mouse^{\gamma+1}$. 
\ This completes the proof of Case 3. Hence, the proof of the Lemma is now complete.\end{proof}

We can now prove
Theorem \ref{GenDJ}.

\begin{proof}[Proof of Theorem \ref{GenDJ}]
Suppose that \ $\mouse$ \ is an iterable real premouse.  One proves
by induction on \ $\gamma \le \widehat{\OR}^{\mouse}$ \ that \ $\mouse^\gamma$ \ is acceptable above the reals. 
Lemma \ref{lemma3} implies that \ $\mouse^\gamma=(J_\gamma(\R),\R)$ \ is acceptable for all
$\gamma\le\kappa$. \ Lemmas \ref{lemma3}, \ref{lemma5} and \ref{lemma6} ensure that \
$\mouse^{\gamma+1}$ \ is acceptable, assuming that \ $\mouse^\gamma$ \ is acceptable for \ $\gamma\ge\kappa$. \ If \
$\gamma$ \ is a limit ordinal and \ $\mouse^\delta$ \ is acceptable for all \ $\delta<\gamma$, \ it then follows
directly that \ $\mouse^\gamma$ \ is acceptable.
\end{proof}

Lemma \ref{lemma1} and Lemma \ref{lemmatwo} are essential ingredients for the proof of the above theorem. These lemmas will also allow us to solve a problem for constructing scales in \ $\Kr$. \ Let \ $\mouse$ \ be an iterable real premouse and suppose that an {\it arbitrary\/} set \ $A$ \ is constructed in \  $M^{\gamma+1}\setminus M^{\gamma}$, \
where
\ $\kappa^\mouse\le\gamma<{\OR}^\mouse$. \ Since we are using the measure \ $\mu^\mouse$ \ to construct the new set \ $A$ \ in
\ $M^{\gamma+1}$, \ it is in fact possible that \ $A\in
M^{\gamma+1}\setminus\boldface{\Sigma}{\omega}(\mouse^\gamma)$.\footnote{``after all, the measure must be unpredictable
somewhere'' \cite[p. 88]{Dodd}.} \ However, Lemma \ref{lemma1} and Lemma \ref{lemmatwo} imply that this cannot happen
{\it when  \ $A$ \ is a set of reals\/}, as we will show in our next lemma. This observation is important in
\cite{Part3} for our determination of whether or not \ $A$ \ has a scale of minimal complexity in \ $\Kr$ \ (recall
Question (Q) in the introduction).
\begin{lemma}\label{minicomp} Let \ $\mouse=(M,\R,\kappa,\mu)$ \ be an iterable real premouse. Let \ $\gamma <
\widehat{\OR}^{\mouse}$. \ For any set of reals \ $A$, \ if \ $A\in M^{\gamma+1}\setminus
M^{\gamma}$, \ then \ $A\in\boldface{\Sigma}{\omega}(\mouse^\gamma)$. 
\end{lemma}

\begin{proof} By Theorem \ref{GenDJ}, \ $\mouse^\gamma$ \ is acceptable. Now, assume that \ $A$ \ is a
set of reals such that \ $A\in M^{\gamma+1}\setminus
M^{\gamma}$. \ We will show that \ $A\in\boldface{\Sigma}{\omega}(\mouse^\gamma)$.
\begin{claiim} $\rho_{\mouse^\gamma}^{n} \le \kappa$ \ for some \ $n\in\omega$.
\end{claiim}

\begin{proof}[Proof of Claim] Assume for a contradiction that \ $\rho_{\mouse^\gamma}^{n} > \kappa$ \ for all \
$n\in\omega$. \ Lemma \ref{lemma4} implies that
\ $M^{\gamma+1}\cap\pow(\R)=M^{\gamma}\cap\pow(\R)$. \ This contradicts our assumption that \ $A\in M^{\gamma+1}\setminus
M^{\gamma}$. \ The proof of the claim is now complete.
\end{proof}

Since \ $\mouse^\gamma$ \ is acceptable, the Claim implies that \ $\mouse^\gamma$ \ is a real mouse. Let \ $n=n(\mouse)$ \ and
so, $\rho_{\mouse^\gamma}^{n} \le \kappa$. \ If \ $\rho_{\mouse^\gamma}^{n} = \kappa$, \ then 
Lemma \ref{lemmatwo} \ asserts that \ $H^{\mouse^\gamma}_{\kappa}=H^{\mouse^{\gamma+1}}_{\kappa}$. \ Thus, \ $A\in
M^\gamma$ \ which contradicts our assumption. \ From this contradiction we conclude that \ $\rho_{\mouse^\gamma}^{n} <
\kappa$. \  Lemma \ref{lemma1} now implies that \ $\pow(M^\gamma)\cap
M^{\gamma+1}\subseteq \boldface{\Sigma}{\omega}(\mouse^\gamma)$. \ Therefore, \
$A\in\boldface{\Sigma}{\omega}(\mouse^\gamma)$.
\end{proof}

By an argument similar to the one proving Lemma \ref{minicomp}, one can show the next lemma.
\begin{lemma}\label{sharp} Let \ $\mouse=(M,\R,\kappa,\mu)$ \ be an iterable real premouse. Let \ $\gamma <
\widehat{\OR}^{\mouse}$. \ If \ $H^{\mouse^\gamma}_{\kappa}\ne H^{\mouse^{\gamma+1}}_{\kappa}$, \ then \
$\mouse^\gamma$ \ is a real mouse and \ $\rho_{\mouse^\gamma}^{n+1} < \kappa$ \ where \ $n=n(\mouse^\gamma)$. 
\end{lemma}
\subsubsection*{A final note} Theorem \ref{GenDJ} establishes that an iterable {\it real\/} premouse \ $\mouse$ \
is acceptable above the reals. Recall that an iterable real premouse \ $\mouse$ \ contains all the reals, that is, \
$\R^\mouse=\R$. \ Our definition (see Definition \ref{premice}) of an iterable premouse \ $\nouse$, \ however,  only
requires that \
$\R^\nouse\subseteq\R$; \ and so, \ $\nouse$ \ need not contain {\it all\/} of the reals.\footnote{Such iterable 
premice exist in any generic extension of \ $\Kr$ \ which adds reals.} Similarly, our definition of
``acceptability above reals''  does not presume that the relevant structure contains all of the reals (see Definition
\ref{acceptable}). Our proof of Theorem
\ref{GenDJ} also does not require \ $\mouse$ \ to contain all of the reals. 
\begin{theorem} Suppose that \ $\mouse$ \ is an iterable premouse.  Then \ $\mouse$ \ is
acceptable above the reals.
\end{theorem}

\providecommand{\bysame}{\leavevmode\hbox to3em{\hrulefill}\thinspace}
\providecommand{\href}[2]{#2}

\end{document}